\pdfoutput=1
\documentclass[final,3p,times,fleqn]{elsarticle}

\usepackage{amsmath,amssymb,amscd,mathtools}
\usepackage{mathrsfs}
\usepackage{enumitem}
\usepackage[colorinlistoftodos,prependcaption,textsize=scriptsize]{todonotes}

\usepackage{booktabs,cellspace}

\usepackage{graphicx}
\usepackage{xcolor}
\usepackage{algorithm}
\usepackage{algpseudocode}
\usepackage{diagbox}

\usepackage{lineno}
\modulolinenumbers[5]
\usepackage[colorlinks=true, allcolors=blue]{hyperref}
\usepackage{tikz}

\usepackage{pgfplots}
\usepgfplotslibrary{patchplots}
\usepackage{xparse}

\definecolor{mycolor1}{rgb}{0.00000,0.44700,0.74100}%
\definecolor{mycolor2}{rgb}{0.85000,0.32500,0.09800}%
\definecolor{mycolor3}{rgb}{0.92900,0.69400,0.12500}%
\definecolor{mycolor4}{rgb}{0.49400,0.18400,0.55600}%
\definecolor{mycolor5}{rgb}{0.46600,0.67400,0.18800}%
\definecolor{mycolor6}{rgb}{0.30100,0.74500,0.93300}%
\definecolor{mycolor7}{rgb}{0.63500,0.07800,0.18400}%

\usetikzlibrary{matrix,backgrounds,calc,arrows}
\pgfdeclarelayer{myback}
\pgfsetlayers{myback,background,main}

\usepackage{subcaption}
\usepackage{etoolbox}
\newcommand{\pder}[3][1]{\frac{\partial^{\ifstrequal{#1}{1}{}{#1}}#2}{\partial#3}}
\newcommand{\norm}[2][]{\left\|#2\right\|_{#1}}
\newcommand{\abs}[2][]{\left|#2\right|_{#1}}
\newcommand{\scal}[3][]{\int_{#1}#2\cdot #3}
\newcommand{\ah}[3][]{a_h^{#1}\ifstrempty{#2#3}{}{\left(#2,#3\right)}}
\newcommand{\vemstab}[3][]{S_h^{#1}\ifstrempty{#2#3}{}{\left(#2,#3\right)}}
\newcommand{\lebl}[2][2]{L^{#1}\ifstrempty{#2}{}{\left(#2\right)}} % Lebesgue spaces (default L^2): L^2(\Omega) = \lebl{\Omega}
\newcommand{\sobh}[3][]{H^{#2}_{#1}\ifstrempty{#3}{}{\left(#3\right)}} % Sobolev spaces H^1_0(\Omega) = \sobh[0]{1}{\Omega}
\newcommand{\proj}[3][0]{\Pi^{#1}_{#2}\ifstrempty{#3}{}{\left(#3\right)}} % VEM projectors
\newcommand{\II}[2][h]{\mathcal{I}_{#1}\ifstrempty{#2}{}{\left(#2\right)}}
\newcommand{\displ}{\tensorOne{u}}
\newcommand{\stress}{\tensorTwo{\sigma}}
\newcommand{\strain}{\tensorTwo{\varepsilon}}
\newcommand{\Th}{\mathcal{T}} % tessellation of Omega
\newcommand{\Poly}[2]{\mathbb{P}_{#1}\ifstrempty{#2}{}{\left(#2\right)}} % space of polynomials
\newcommand{\Monom}[2]{\mathcal{M}_{#1}\ifstrempty{#2}{}{\left(#2\right)}} % space of scaled monomials

\newcommand{\tensorOne}[1]{\boldsymbol{#1}}
\newcommand{\tensorTwo}[1]{\boldsymbol{#1}} % Ideally Greek letter

\renewcommand{\Vec}[1]{%
  \ifcat\noexpand#1\relax % check if the argument is a control sequence (Greek)
    \boldsymbol{#1}% Greek
  \else
    \mathbf{#1}% single character
  \fi
}
\newcommand{\Mat}[1]{#1}

\newcommand{\blkVec}[1]{\Vec{#1}}
\newcommand{\blkMat}[1]{\Vec{#1}}

\newcommand{\sub}[1]{_{#1}}
\newcommand{\transpose}{^{\intercal}}

\newdefinition{rmk}{Remark}

\journal{Computer Methods in Applied Mechanics and Engineering}

\newcommand{\logLogSlopeTriangle}[5]
{
    \pgfplotsextra
    {
        \pgfkeysgetvalue{/pgfplots/xmin}{\xmin}
        \pgfkeysgetvalue{/pgfplots/xmax}{\xmax}
        \pgfkeysgetvalue{/pgfplots/ymin}{\ymin}
        \pgfkeysgetvalue{/pgfplots/ymax}{\ymax}

        \pgfmathsetmacro{\xArel}{#1}
        \pgfmathsetmacro{\yArel}{#3}
        \pgfmathsetmacro{\xBrel}{#1+#2}
        \pgfmathsetmacro{\yBrel}{\yArel}
        \pgfmathsetmacro{\xCrel}{\xArel}
        \pgfmathsetmacro{\lnxB}{\xmin*(1-(#1-#2))+\xmax*(#1-#2)} % in [xmin,xmax].
        \pgfmathsetmacro{\lnxA}{\xmin*(1-#1)+\xmax*#1} % in [xmin,xmax].
        \pgfmathsetmacro{\lnyA}{\ymin*(1-#3)+\ymax*#3} % in [ymin,ymax].
        \pgfmathsetmacro{\lnyC}{\lnyA+#4*(\lnxA-\lnxB)}
        \pgfmathsetmacro{\yCrel}{\lnyC-\ymin)/(\ymax-\ymin)} % THE IMPROVED EXPRESSION WITHOUT 'DIMENSION TOO LARGE' ERROR.
        \coordinate (A) at (rel axis cs:\xArel,\yArel);
        \coordinate (B) at (rel axis cs:\xBrel,\yBrel);
        \coordinate (C) at (rel axis cs:\xCrel,\yCrel);

        \draw[#5]   (A)-- node[pos=0.5,anchor=north] {\footnotesize 1}
                    (B)-- 
                    (C)-- node[pos=0.5,anchor=east] {\footnotesize{#4}}
                    cycle;
    }
}
\begin{document}

\title{Hybrid mimetic finite-difference and virtual element formulation for coupled poromechanics}

\begin{frontmatter}
\author[polito]{Andrea~Borio\corref{cor1}}
\ead{andrea.borio@polito.it}
\author[total]{Fran\c cois~P.~Hamon}
\ead{francois.hamon@total.com}
\author[aeed]{Nicola~Castelletto}
\ead{castelletto1@llnl.gov}
\author[aeed]{Joshua~A.~White}
\ead{white230@llnl.gov}
\author[aeed]{Randolph~R.~Settgast}
\ead{settgast1@llnl.gov}
% \author[stanford]{Hamdi~A.~Tchelepi}
% \ead{tchelepi@stanford.edu}

\cortext[cor1]{Corresponding author.}
\address[polito]{Dipartimento di Scienze Matematiche, Politecnico di Torino, Torino, 10129, Italy} 
\address[total]{Total E\&P Research and Technology, Houston, TX 77002, USA}
\address[aeed]{Atmospheric, Earth, and Energy Division, Lawrence Livermore National Laboratory, Livermore, CA 94550, USA}
% \address[stanford]{Department of Energy Resources Engineering, Stanford University, Stanford, CA 94305, USA}

\begin{abstract}
We present a hybrid mimetic finite-difference and virtual element formulation for coupled single-phase poromechanics on unstructured meshes.
The key advantage of the scheme is that it is convergent on complex meshes containing highly distorted cells with arbitrary shapes.
We use a local pressure-jump stabilization method based on unstructured macro-elements to prevent the development of spurious pressure modes in incompressible problems approaching undrained conditions.
A scalable linear solution strategy is obtained using a block-triangular preconditioner designed specifically for the saddle-point systems arising from the proposed discretization.
The accuracy and efficiency of our approach are demonstrated numerically on two-dimensional benchmark problems.
\end{abstract}


\end{frontmatter}

\allowdisplaybreaks

% \linenumbers

\section{Introduction}
\label{sec:introduction}

Modeling hydro-mechanical coupling is essential to accurately simulate a wide range of subsurface processes involving fluid flow and mechanical deformation, such as oil and gas recovery \cite{zoback2010reservoir}, geological CO$_2$ storage \cite{rutqvist2012geomechanics}, and geothermal energy production \cite{majer2007induced}.
In life sciences, this coupling also plays a central role in the modeling of bone deformation \cite{cowin1999bone} and blood-vessel interaction in hemodynamics \cite{badia2009coupling}.
%
%Numerical simulations are often used to investigate solid-fluid interaction.
Studies involving the numerical solution of Biot's equations of poroelasticity on a computational mesh are routinely used to investigate these processes.
In most subsurface applications, generating a mesh that faithfully represents the structure of the porous medium is a difficult task.
Geological formations often exhibit a high heterogeneity characterized by stratigraphic layering and the presence of faults and fractures.
%
%Preserving a faithful representation of the porous medium in the mesh generation step naturally yields unstructured polyhedral meshes that conform to complex geological features \cite{mallison2014unstructured}.
As a result, there is strong interest in using unstructured polyhedral meshes that conform to the complex geological features of the porous medium \cite{mallison2014unstructured}.
To preserve accuracy and reduce computational cost, it is appealing to solve both the porous flow problem and the mechanical problem on the same mesh.
%
%This is challenging, since practical unstructured grids typically include highly distorted cells with large aspect ratios and arbitrary shapes.
%
In this work, we address this issue by developing a robust numerical scheme and fully coupled solution strategy for the displacement-velocity-pressure formulation of Biot's equations on arbitrary polyhedral meshes.

In recent years, considerable efforts have been invested in the development of stable and convergent numerical schemes for general second-order elliptic problems on polyhedral meshes \cite{da2014mimetic,Beirao2015b,di2019hybrid,droniou2018gradient}.
In poromechanical simulations, the mechanical problem has traditionally been solved with the Finite-Element Method (FEM).
To overcome the mesh restrictions imposed by the standard FEM and handle arbitrary cell shapes, generalized finite-element discretizations have been proposed for polyhedral meshes, such as the polyhedral finite-element method (see \cite{manzini2014new} and references therein).
For completeness, we mention here that discontinuous Galerkin methods, hybrid high-order methods \cite{botti2019hybrid}, and finite-volume approaches \cite{nordbotten2016stable,terekhov2020cell} have also been designed to solve the mechanical component of Biot's equations on polyhedral meshes.
Recently, using concepts inspired by Mimetic Finite Difference (MFD) methods, the introduction of the Virtual Element Method (VEM) provided a variational framework to construct a consistent and stable scheme on arbitrary polyhedral meshes \cite{Beirao2015a,Beirao2015b}.
An attractive feature that distinguishes VEM from other generalized FEMs is that the assembly of the VEM discrete equations does not require  knowledge of the analytical expression of the (non-polynomial) basis functions, which allows a generic and robust treatment of complex cell geometries.
For the same reason, hanging nodes can be dealt with in a simple and convenient fashion.
The scheme's properties have been studied extensively \cite{cangiani2017posteriori,berrone2017orthogonal,brenner2017some,BBBPSsupg,BBapost}, and have been assessed numerically in a wide range of physical simulations, including geomechanical applications \cite{andersen2017virtual,coulet2019fully,BBBPS,BBS}.
In this work, we take advantage of the appealing features discussed above, using a low-order VEM relying on vertex-based degrees of freedom (dofs) to approximate displacements in the discrete mechanical problem.

For the flow problem, a large number of numerical schemes for polyhedral meshes have been proposed \cite{lie2019introduction}, and a complete review is out of the scope of the present work.
We focus here on locally conservative, low-order, linear schemes.
Pioneering linear finite-volume approaches have relied on a cell-centered discretization, such as Multi-Point Flux Approximation (MPFA) schemes \cite{edwards1998finite,aavatsmark2002introduction}.
However, linear schemes combining cell-centered dofs with face-centered dofs---such as the family of Hybrid Mimetic Mixed (HMM) methods \cite{aarnes2008multiscale,droniou2010unified,abushaikha2020fully}---or with vertex-based dofs---such as the Vertex Approximate Gradient (VAG) scheme \cite{eymard2012small}---have been shown to be convergent on polyhedral meshes.
%
%Nonlinear finite-volume methods have also been designed (e.g., \cite{schneider2018monotone}).
%
A recent numerical comparison of some of these schemes with nonlinear methods can be found in \cite{klemetsdal2020comparison}, and a general analysis framework is provided by \cite{droniou2018gradient,bonaldi2020gradient}.
Motivated by the fact that VEM and mimetic schemes share a common theoretical basis, we make a natural choice---not studied before in a poromechanical context to the best of our knowledge---to couple VEM with a low-order hybrid MFD discretization that approximates pressure with cell-centered dofs and velocity with face-centered dofs.
In the hybrid formulation, the MFD scheme also involves face-based Lagrange multipliers, which allows a cell-wise assembly of the discrete flow equations.
Static condensation is used to eliminate the velocities during assembly and reduce the size of the linear systems.

The proposed hybrid MFD-VEM numerical scheme uses the lowest order Virtual Element space for the displacement field and a piecewise-constant interpolation for pressure.
This combination of displacement and pressure approximation spaces does not intrinsically satisfy the Ladyzhenskaya-Babu{\v{s}}ka-Brezzi (LBB) \textit{inf-sup} stability condition \cite{babuvska1971error,brezzi1974existence}.
For incompressible problems approaching undrained conditions, this results in the development of spurious modes in the pressure field for specific mesh topologies \cite{da2010mimetic}.
To circumvent this issue, the scheme is stabilized using a local pressure-jump method inspired by \cite{camargo2020macroelement,frigo2020efficient} and adapted here to unstructured meshes.
We design a fully implicit, fully coupled solution strategy for the stabilized discretization.
The workhorse of the algorithm is a block-triangular preconditioner constructed specifically for the linear systems arising from the hybrid MFD-VEM scheme using a methodology presented in \cite{frigo2020efficient}.
Using this, we demonstrate three key features of our numerical framework applied to arbitrary polygonal meshes with highly distorted cells.
First, the scheme is convergent upon space-time refinement and matches the analytical solution of the well-known benchmark problems \cite{Mandel1953}.
Second, the scheme exhibits an accurate behavior for incompressible problems with small time steps.
We note in particular that no stabilization is needed on a large class of arbitrary polygonal meshes \cite{talischi2014polygonal}, and that local stabilization effectively damps spurious pressure modes on those meshes triggering the instability.
Third, the linear solution strategy is robust and scalable for all the polygonal meshes considered here.
These encouraging results will be extended to three-dimensional cases in future work.

This paper is organized as follows.
We review the strong and weak forms of Biot's poroelasticity equations in Section~\ref{sec:model_problem}.
In Section~\ref{sec:numerical_model}, we introduce the hybrid MFD-VEM scheme to solve the initial boundary value problem.
We also introduce the local pressure-jump method used to stabilize the scheme.
Section~\ref{sec:solution_strategy} is dedicated to the presentation of the linear solution algorithm for the coupled systems.
In Section~\ref{sec:numerical_examples}, we demonstrate the accuracy of the proposed numerical scheme on polygonal meshes and we illustrate the scalability of the solution algorithm.
The paper ends with a few concluding remarks regarding future work.

\section{Model problem}
\label{sec:model_problem}

% mainfile: main.tex
\subsection{Initial-boundary value problem in strong form}

We consider a displacement-velocity-pressure formulation of Biot's poroelasticity equations \cite{biot1941general,coussy2004poromechanics} in a two-dimensional domain $\Omega \in \mathbb{R}^2$.
Let $\mathcal{I} = (0, T)$ denote the time interval.
The three-field strong form of the initial-boundary value problem (IBVP) consists of a linear momentum balance equation, a mass balance equation, and Darcy's law.
Using an \emph{excess pressure} formulation, the displacement $\tensorOne{u} : \overline{\Omega} \times \mathcal{I} \rightarrow \mathbb{R}^2$, the Darcy velocity $\tensorOne{q} : \overline{\Omega} \times \mathcal{I} \rightarrow \mathbb{R}^2$, and the excess pore pressure $p : \overline{\Omega} \times \mathcal{I} \rightarrow \mathbb{R}$ satisfy:
\begin{linenomath}
\begin{subequations}
\begin{align}
  -\text{div} \, \tensorTwo{\sigma} ( \tensorOne{u}, p )  
  &=
  \tensorOne{b} & &\mbox{ in } \Omega \times \mathcal{I}
  & &\mbox{(momentum balance)}, \label{eq:momentumBalance}\\
  \tensorTwo{\kappa}^{-1} \cdot \tensorOne{q} + \nabla p
  &= 
  \tensorOne{0} & &\mbox{ in } \Omega \times \mathcal{I}
  & &\mbox{(Darcy's law)}, \label{eq:darcy}  \\   
  \dot{\zeta} ( \tensorOne{u}, p ) + \text{div} \, \tensorOne{q} &=
  0 & &\mbox{ in } \Omega \times \mathcal{I}
  & &\mbox{(mass balance)}. \label{eq:massBalance}
\end{align}\label{eq:IBVP}\null
\end{subequations}
\end{linenomath}
In the momentum balance equation, the total Cauchy stress tensor is $\tensorTwo{\sigma} (\tensorOne{u}, p) = \tensorTwo{\sigma}^{\prime}( \tensorOne{u} ) - \alpha p \tensorTwo{I} $, where $\tensorTwo{\sigma}^{\prime}( \tensorOne{u} )$ is the effective stress tensor, $\alpha$ is Biot's coefficient, and $\tensorTwo{I}$ is the second-order unit tensor.
In this work, we consider isotropic linear elastic materials.
Hence, the effective stress can be defined as $\tensorTwo{\sigma}^{\prime}( \tensorOne{u} ) = 2G \strain( \tensorOne{u} )+ \lambda \, \text{trace} (\strain( \tensorOne{u} )) $, with $\strain( \tensorOne{u} ) = \frac{1}{2} (\nabla \tensorOne{u} + \nabla \tensorOne{u} ^T) $ the linearized strain tensor, and $\lambda$ and $G$ the Lam\'{e} parameters of the porous medium. The vector $\tensorOne{b}(\tensorOne{x},t)$ denotes body forces.
In Darcy's law, $\tensorTwo{\kappa}$ denotes the intrinsic permeability tensor of the porous medium divided by the fluid viscosity, which is assumed constant.
Using the superposed dot, $\dot{()}$, to denote a derivative with respect to time, the fluid increment \cite{biot1941general} is defined as $\dot{\zeta} ( \tensorOne{u}, p ) = \alpha \, \text{div} \, \dot{\tensorOne{u}} + S_{\varepsilon} \dot{p}$, where $S_{\varepsilon}$ is the constrained specific storage coefficient, i.e. the inverse of Biot's modulus. 

To define the boundary conditions in the mechanical and flow problems, the domain boundary $\Gamma$ is decomposed as $\Gamma = \overline{\Gamma_u \cup \Gamma_{\sigma}}$ and $\Gamma = \overline{\Gamma_p \cup \Gamma_q}$ such that $\Gamma_u \cap \Gamma_{\sigma} = \emptyset$ and $\Gamma_p \cap \Gamma_q = \emptyset$.
We prescribe the following boundary conditions to, respectively, the displacement, total Cauchy stress tensor, Darcy velocity, and excess pore pressure fields:
\begin{linenomath}
\begin{subequations}
\begin{align}
  \tensorOne{u}  
  &= \overline{\tensorOne{u}} & \mbox{ on } \Gamma_u \times \mathcal{I}, \label{eq:displacement_bc}\\
  \tensorTwo{\sigma} \cdot \tensorOne{n}
  &= \overline{\tensorOne{t}} & \mbox{ on } \Gamma_{\sigma} \times \mathcal{I}, \label{eq:traction_bc}\\
  \tensorOne{q} \cdot \tensorOne{n}
  &= \overline{q} & \mbox{ on } \Gamma_{q} \times \mathcal{I}, \label{eq:flux_bc}\\
  p
  &= \overline{p} & \mbox{ on } \Gamma_{p} \times \mathcal{I}, \label{eq:pressure_bc}
\end{align}\label{eq:boundary_conditions}\null
\end{subequations}
\end{linenomath}
where the (space- and time-dependent) boundary values on the right-hand sides are denoted by the notation $(\overline{\cdot})$.
In \eqref{eq:traction_bc}, $\overline{\tensorOne{t}}$ is the prescribed traction, and $\tensorOne{n}$ is the outward normal vector on $\Gamma$.

To complete the definition of the problem, we impose the following initial condition on the excess pore pressure field:
\begin{linenomath}
\begin{equation}
%\tensorOne{u}( \tensorOne{x}, 0 ) = \tensorOne{u}^0, \qquad p( \tensorOne{x}, 0 ) = p^0, \qquad \tensorOne{x} \in \overline{\Omega}. 
  p( \tensorOne{x}, 0 ) = p^0, \qquad \tensorOne{x} \in \overline{\Omega},
\end{equation}
\end{linenomath}
with $p^0$ its initial value.  Frequently, though not always, $p^0=0$.  The initial displacement $\tensorOne{u}^0$ and velocity $\tensorOne{q}^0$ are then computed such that Eqs. \eqref{eq:momentumBalance} and \eqref{eq:darcy}, respectively, are satisfied.

\subsection{Weak statement of the model problem}

In the weak statement of the model problem, we retain a mixed structure in which mass balance and Darcy's law are kept as separate equations.
This mixed formulation will be used in Section~\ref{sec:fully_discrete_coupled_scheme} to construct a consistent discretization of the flow problem based on a hybrid MFD scheme.
The weak form of the IBVP \eqref{eq:IBVP} involves the following spaces:
\begin{linenomath}
\begin{subequations}
\begin{align}
  \tensorOne{\mathcal{U}} &= \{ \tensorOne{u} \in \tensorOne{H}^1( \Omega ): \, \tensorOne{u}_{|\Gamma_u} = \overline{\tensorOne{u}} \},
  & & & \tensorOne{\mathcal{U}}_0 &= \{ \tensorOne{u} \in \tensorOne{H}^1( \Omega ): \, \tensorOne{u}_{|\Gamma_u} = \tensorOne{0} \}, 
  \\
  \tensorOne{\mathcal{Q}} &= \{ \tensorOne{q} \in \tensorOne{H}( \text{div}; \Omega ): \, \tensorOne{q} \cdot \tensorOne{n}_{|\Gamma_q} = \overline{q} \},
  & & & \tensorOne{\mathcal{Q}}_0 &= \{ \tensorOne{q} \in \tensorOne{H}( \text{div}; \Omega ): \, \tensorOne{q} \cdot \tensorOne{n}_{|\Gamma_q} = 0 \},
  \\
  \mathcal{P} &= L^2( \Omega ), & & & &
\end{align}\label{eq:sobolev_spaces}\null
\end{subequations}
\end{linenomath}
where $\tensorOne{H}^1( \Omega )$ and $\tensorOne{H}( \text{div}; \Omega )$ are the Sobolev spaces containing, respectively, the vector functions whose first derivatives belong to $L^2( \Omega )$, and the vector functions whose divergence is in $L^2( \Omega )$.
Using the notation $( \cdot, \cdot )_{\square}$ to denote the inner product on the functional space specified in the subscript, the weak form of the mixed IBVP \eqref{eq:IBVP} consists in finding the triplet $( \tensorOne{u}( t ), \tensorOne{q}(t), p(t) ) \in \tensorOne{\mathcal{U}} \times \tensorOne{\mathcal{Q}} \times \mathcal{P}$ that, for all $t \in \mathcal{I}$, satisfies:
\begin{linenomath}
\begin{subequations}
\begin{align}
  &( \tensorTwo{\sigma}^{\prime}, \strain ( \tensorOne{\eta} ) )_{[L^2(\Omega)]^{2 \times 2}}
  - ( \alpha p, \text{div} \, \tensorOne{\eta} )_{L^2(\Omega)}
  = (\tensorOne{b} , \tensorOne{\eta})_{[L^2(\Omega)]^2}
  + ( \overline{\tensorOne{t}}, \tensorOne{\eta} )_{[L^2(\Gamma_{\sigma})]^2}
  & & \forall \tensorOne{\eta} \in \tensorOne{\mathcal{U}}_0, \label{momentumBalance_weak_form} \\
  &( \tensorTwo{\kappa}^{-1} \cdot \tensorOne{q}, \tensorOne{\varphi} )_{[L^2(\Omega)]^2} - ( p, \text{div} \, \tensorOne{\varphi} )_{L^2(\Omega)} = - ( \overline{p}, \tensorOne{\varphi} \cdot \tensorOne{n} )_{L^2(\Gamma_p)} & & \forall \tensorOne{\varphi} \in \mathcal{\tensorOne{Q}}_0, \label{darcy_weak_form} \\
  &(\alpha \, \text{div} \, \dot{\tensorOne{u}}, \chi )_{L^2(\Omega)} + ( \text{div} \, \tensorOne{q}, \chi )_{L^2(\Omega)} + ( S_{\varepsilon} \dot{p}, \chi )_{L^2(\Omega)} = 0 & & \forall \chi \in L^2( \Omega ). \label{massBalance_weak_form}
\end{align}\label{eq:weak_form_IBVP}\null
\end{subequations}
\end{linenomath}
A detailed analysis of this problem can be found in \cite{lipnikov2003numerical}.
In the next section, we discretize the weak form \eqref{eq:weak_form_IBVP} on arbitrary polygonal meshes by applying a low-order VEM to the momentum balance equation and a hybrid MFD scheme to the mass balance and Darcy's law.

\section{Numerical model}
\label{sec:numerical_model}

\subsection{Polygonal meshes}
\label{sec:mesh_hierarchy}

To describe the arbitrary polygonal meshes considered in this work, we rely on a mesh structure consisting of the triplet $\{ \mathcal{T}, \mathcal{F}, \mathcal{V} \}$ (see also \cite{eymard2010discretization,nordbotten2016stable}).
In this triplet, $\mathcal{T}$ denotes the set of non-overlapping open polygonal cells in the mesh such that $\overline{\Omega} = \cup_{K \in \mathcal{T}} \overline{K}$.
 Let $\partial K = \overline{K} \setminus K$ be the boundary of cell $K \in \mathcal{T}$.
The second entry in the triplet, $\mathcal{F}$, is the set of mesh faces corresponding, in the two-dimensional case considered here, to one-dimensional hyperplanes of $\mathbb{R}^2$.
The set of mesh faces is decomposed into three non-intersecting subsets $\mathcal{F} = \mathcal{F}_{int} \cup \mathcal{F}_{p} \cup \mathcal{F}_{q}$, containing respectively the interior faces, the boundary faces located on $\Gamma_p$, and the boundary faces located on $\Gamma_q$.
For all $K \in \mathcal{T}$, $\mathcal{F}_K$ is a subset that contains the faces of cell $K$ such that $\partial K = \cup_{f \in \mathcal{F}_K} \overline{f}$.
Denoting by $\mathcal{T}_{f} = \{ K \in \mathcal{T}: \, f \in \mathcal{F}_K \}$ the subset consisting of the cell(s) adjacent to interface $f$, we require that the mesh is conforming in the sense that $\text{card}( \mathcal{T}_{f} ) = 2$ if $f \in \mathcal{F}_{int}$ and $\text{card}( \mathcal{T}_{f} ) = 1$ otherwise.
The third entry of the triplet, $\mathcal{V}$, is the set of mesh vertices.
The subset of the vertices of cell $K$ is $\mathcal{V}_K$, and the subset of the vertices of face $f$ is $\mathcal{V}_f$.

\begin{figure}
  \small
  
  \hfill
  \begin{subfigure}[t]{.3\linewidth}
    \begin{tikzpicture}
        \node[anchor=south west,inner sep=0] (image) at (0,0) {\includegraphics[width=\linewidth]{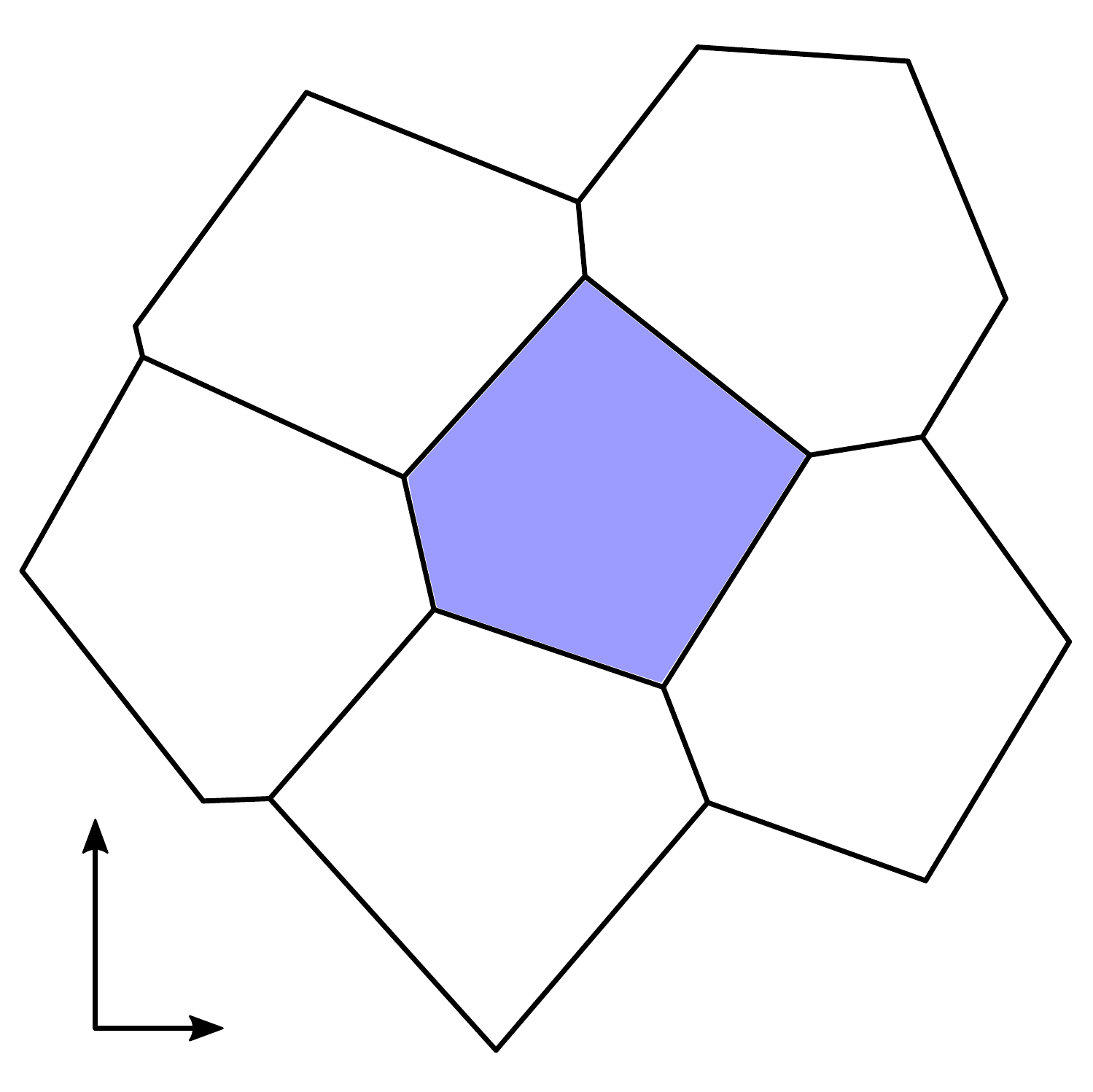}};
        \node [right] at (.19\linewidth,.05\linewidth) {$x$};
        \node [right] at (.01\linewidth,.24\linewidth) {$y$};
        \node [right] at (.475\linewidth,.525\linewidth) {$K$};
    \end{tikzpicture}
    \caption{}
  \end{subfigure}
  \hfill
  \begin{subfigure}[t]{.3\linewidth}
    \begin{tikzpicture}
        \node[anchor=south west,inner sep=0] (image) at (0,0) {\includegraphics[width=\linewidth]{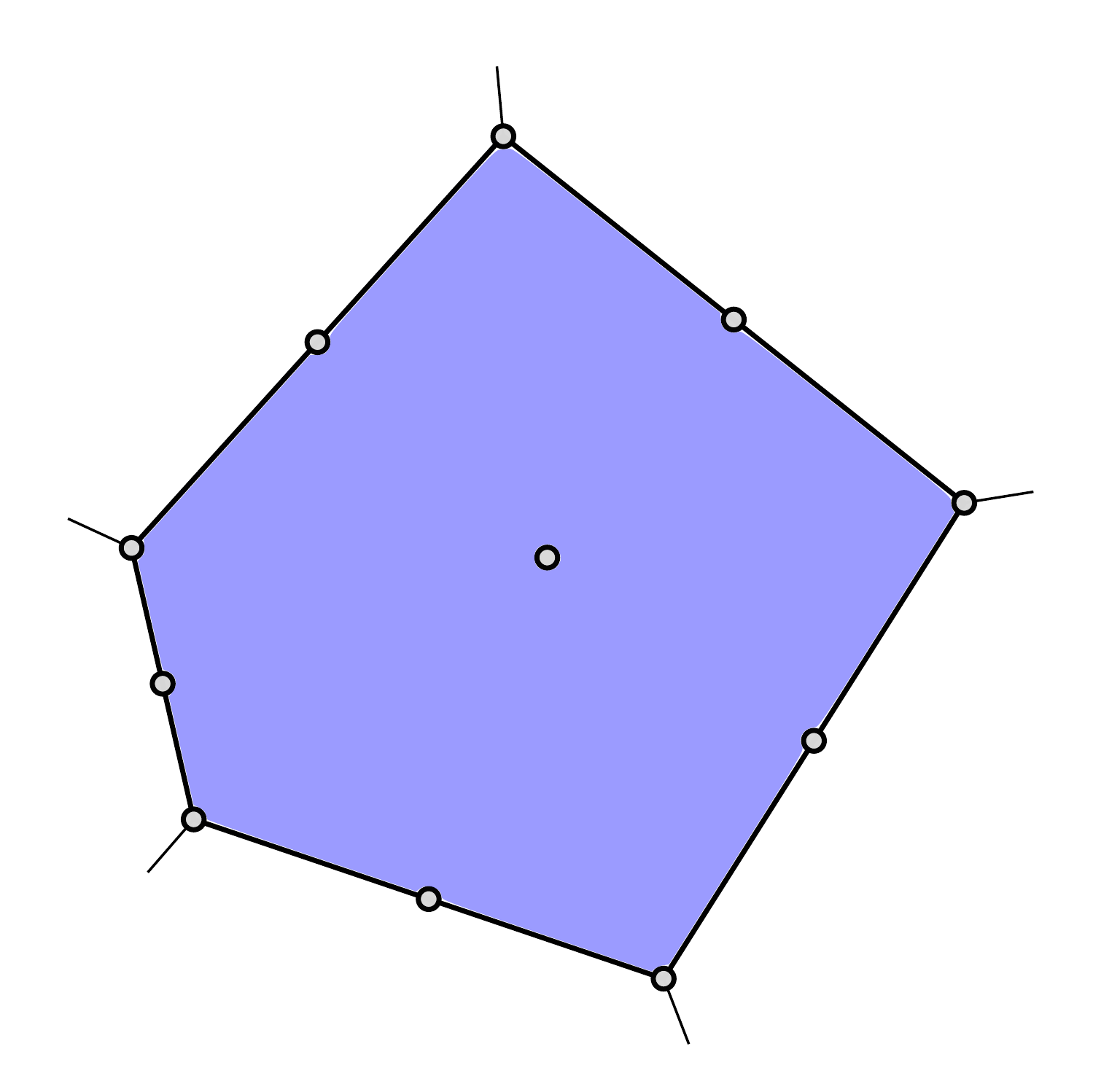}};
        \node [right] at (.5\linewidth,.525\linewidth) {$\tensorOne{x}_K$};
        \node [right] at (.05\linewidth,.55\linewidth) {$\tensorOne{x}_{v_i}$};
        \node [right] at (.15\linewidth,.200\linewidth) {$\tensorOne{x}_{v_j}$};
        \node [right] at (.61\linewidth,.095\linewidth) {$\tensorOne{x}_{v_k}$};
        \node [right] at (.83\linewidth,.580\linewidth) {$\tensorOne{x}_{v_l}$};
        \node [right] at (.46\linewidth,.880\linewidth) {$\tensorOne{x}_{v_m}$};
        \node [right] at (.04\linewidth,.370\linewidth) {$\tensorOne{x}_{f_i}$};
        \node [right] at (.35\linewidth,.130\linewidth) {$\tensorOne{x}_{f_j}$};
        \node [right] at (.75\linewidth,.310\linewidth) {$\tensorOne{x}_{f_k}$};
        \node [right] at (.67\linewidth,.730\linewidth) {$\tensorOne{x}_{f_l}$};
        \node [right] at (.18\linewidth,.720\linewidth) {$\tensorOne{x}_{f_m}$};
    \end{tikzpicture}
    \caption{}
  \end{subfigure}
  \hfill
  \begin{subfigure}[t]{.3\linewidth}
    \begin{tikzpicture}
        \node[anchor=south west,inner sep=0] (image) at (0,0) {\includegraphics[width=\linewidth]{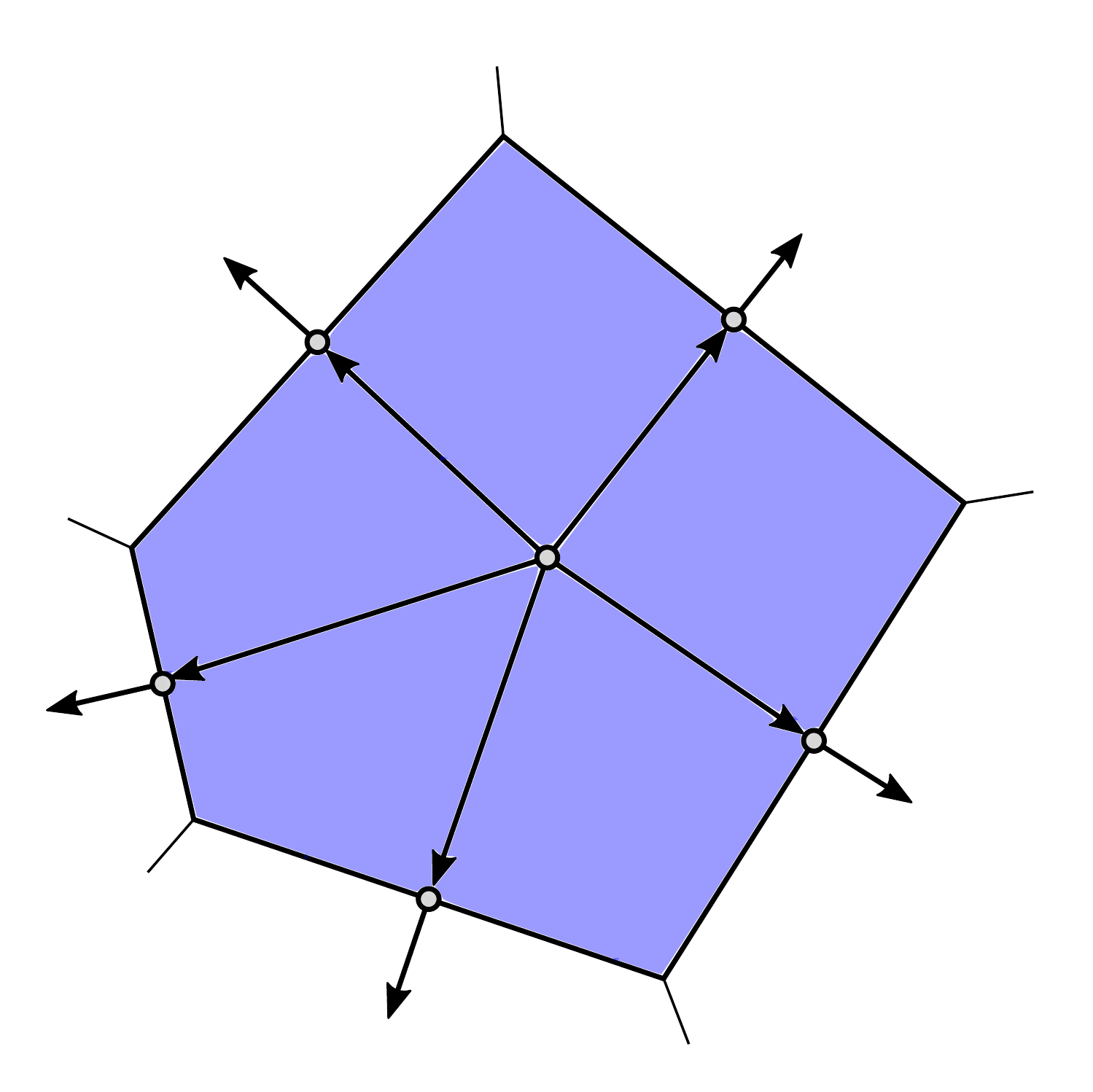}};
        \node [right] at (.00\linewidth,.290\linewidth) {$\tensorOne{n}_{K,f_i}$};
        \node [right] at (.36\linewidth,.080\linewidth) {$\tensorOne{n}_{K,f_j}$};
        \node [right] at (.75\linewidth,.220\linewidth) {$\tensorOne{n}_{K,f_k}$};
        \node [right] at (.68\linewidth,.820\linewidth) {$\tensorOne{n}_{K,f_l}$};
        \node [right] at (.18\linewidth,.800\linewidth) {$\tensorOne{n}_{K,f_m}$};
        \node [right] at (.25\linewidth,.380\linewidth) {$\tensorOne{c}_{K,f_i}$};
        \node [right] at (.44\linewidth,.330\linewidth) {$\tensorOne{c}_{K,f_j}$};
        \node [right] at (.60\linewidth,.420\linewidth) {$\tensorOne{c}_{K,f_k}$};
        \node [right] at (.55\linewidth,.550\linewidth) {$\tensorOne{c}_{K,f_l}$};
        \node [right] at (.36\linewidth,.620\linewidth) {$\tensorOne{c}_{K,f_m}$};
    \end{tikzpicture}
    \caption{}
  \end{subfigure}
  \hfill\null
    
  \caption{\label{fig:mesh_attributes}Mesh attributes.}
\end{figure}
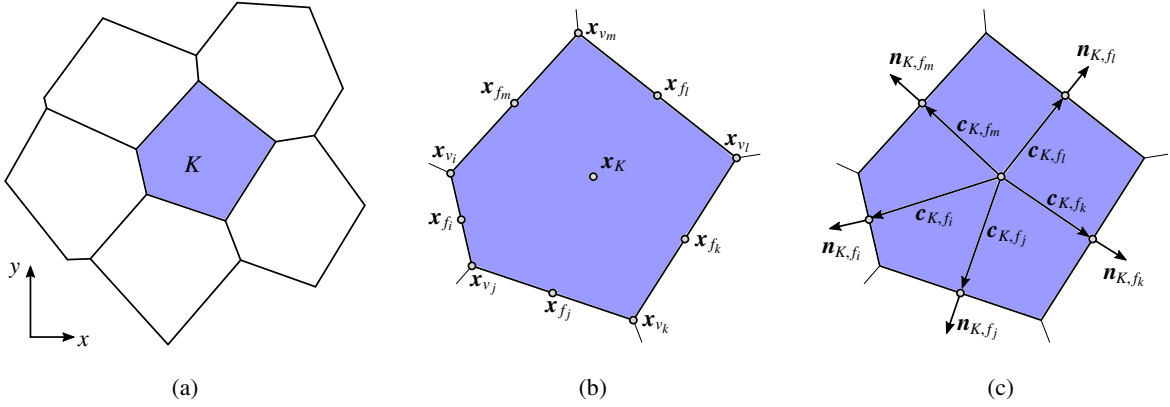

The notation $| \cdot |$ refers to the $d$-measure of a $d$-dimensional quantity (cell or face).
We denote the positions of the center of cell $K$, the center of face $f$, and the vertex $v$ by $\tensorOne{x}_K$, $\tensorOne{x}_f$, and $\tensorOne{x}_v$, respectively (see Fig. \ref{fig:mesh_attributes}).
The geometrical computations in the next sections involve the outward normal to face $f \in \mathcal{T}_K$ with respect to cell $K$, denoted by $\tensorOne{n}_{K,f}$, and the vector connecting $\tensorOne{x}_K$ to $\tensorOne{x}_f$ ($f \in \mathcal{F}_K$), denoted by $\tensorOne{c}_{K,f}$ (see Fig. \ref{fig:mesh_attributes}).
Using these notations, a mesh satisfies the $\tensorTwo{\kappa}$-orthogonality condition \cite{lie2019introduction} if:
\begin{linenomath}
\begin{equation}
  \tensorTwo{\kappa} \cdot \tensorOne{n}_{K,f} \parallel \tensorOne{c}_{K,f} \quad \forall f \in \mathcal{F}_K, \, \forall K \in \mathcal{T}.
  \label{eq:orthogonality_condition}
\end{equation}
\end{linenomath}
We make the regularity assumption that there exists $\gamma > 0$ such that for all $K \in \mathcal{T}$, $K$ is star-shaped with respect to a ball of radius larger than $\gamma h_K$, where $h_K$ is the cell diameter.

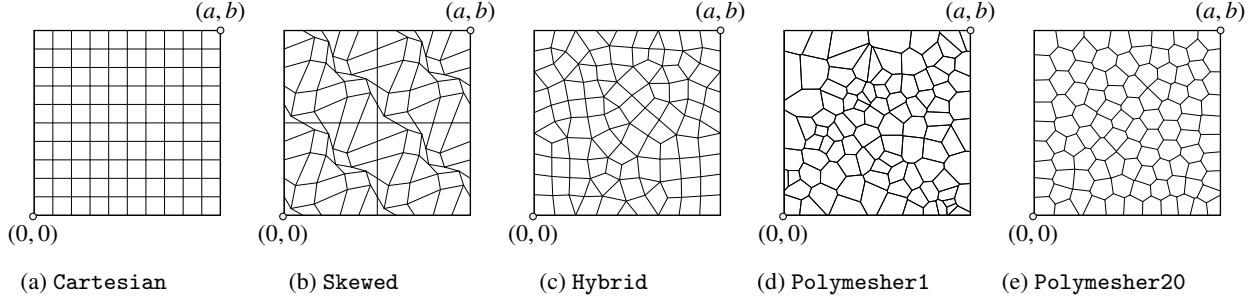
\begin{figure}  
    \small
    
    \begin{subfigure}{.15\linewidth}
      \begin{tikzpicture}[scale=1.0]
        \node[anchor=south west,inner sep=0] (image) at (0,0) {\includegraphics[width=\linewidth]{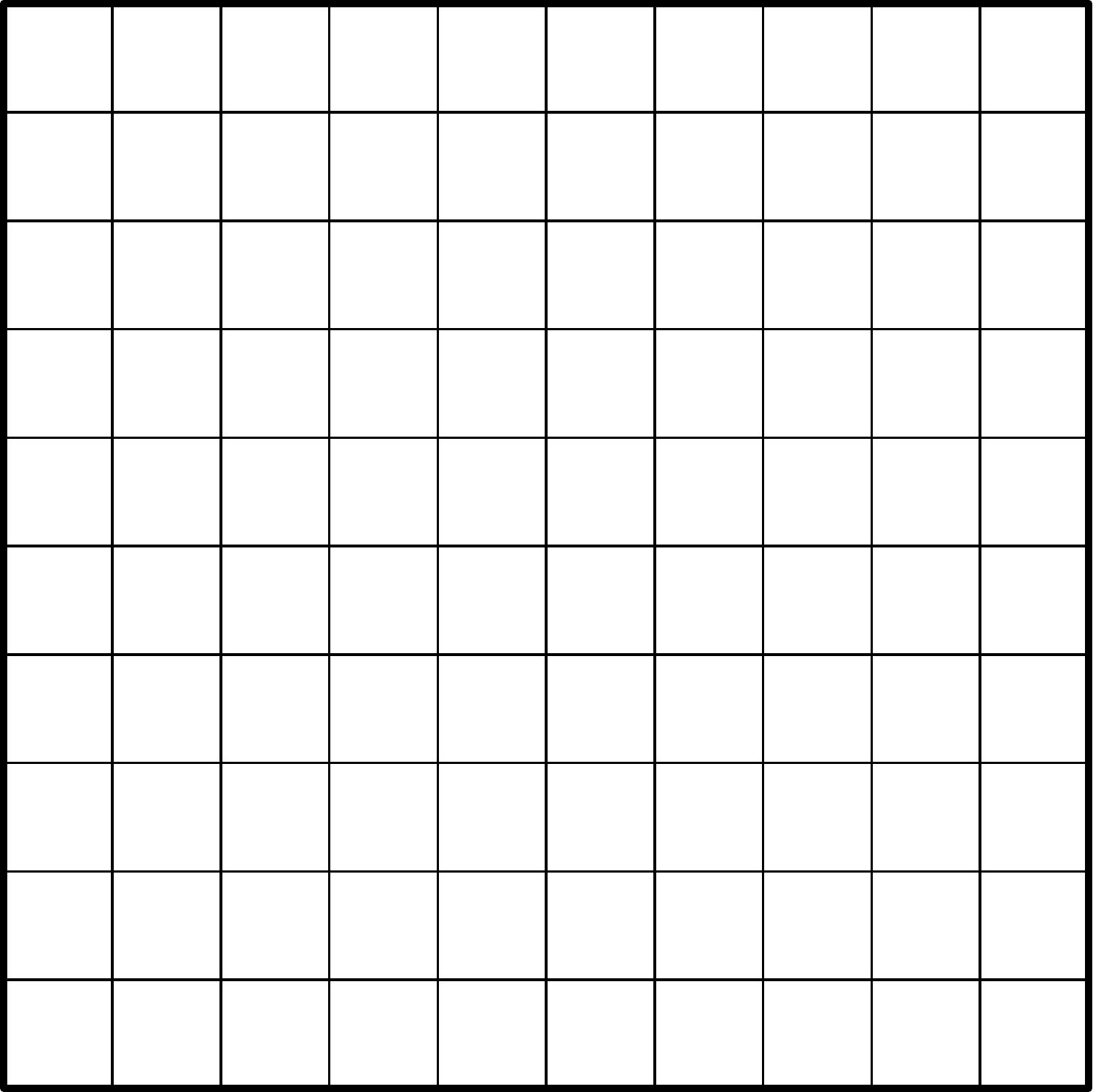}};
         \node[circle,draw=black, fill=gray!15, inner sep=0pt,minimum size=2.5pt] (b) at (0,0) {};
        \node [below] at (0,0) {$(0,0)$};
         \node[circle,draw=black, fill=gray!15, inner sep=0pt,minimum size=2.5pt] (b) at (\linewidth,\linewidth) {};
        \node [above] at (\linewidth,\linewidth) {$(a,b)$};
      \end{tikzpicture}
      \caption{\texttt{Cartesian}}
      \label{fig:meshes-Cart}
    \end{subfigure}
    \hfill    
    \begin{subfigure}{.15\linewidth}
      \begin{tikzpicture}[scale=1.0]
        \node[anchor=south west,inner sep=0] (image) at (0,0) {\includegraphics[width=\linewidth]{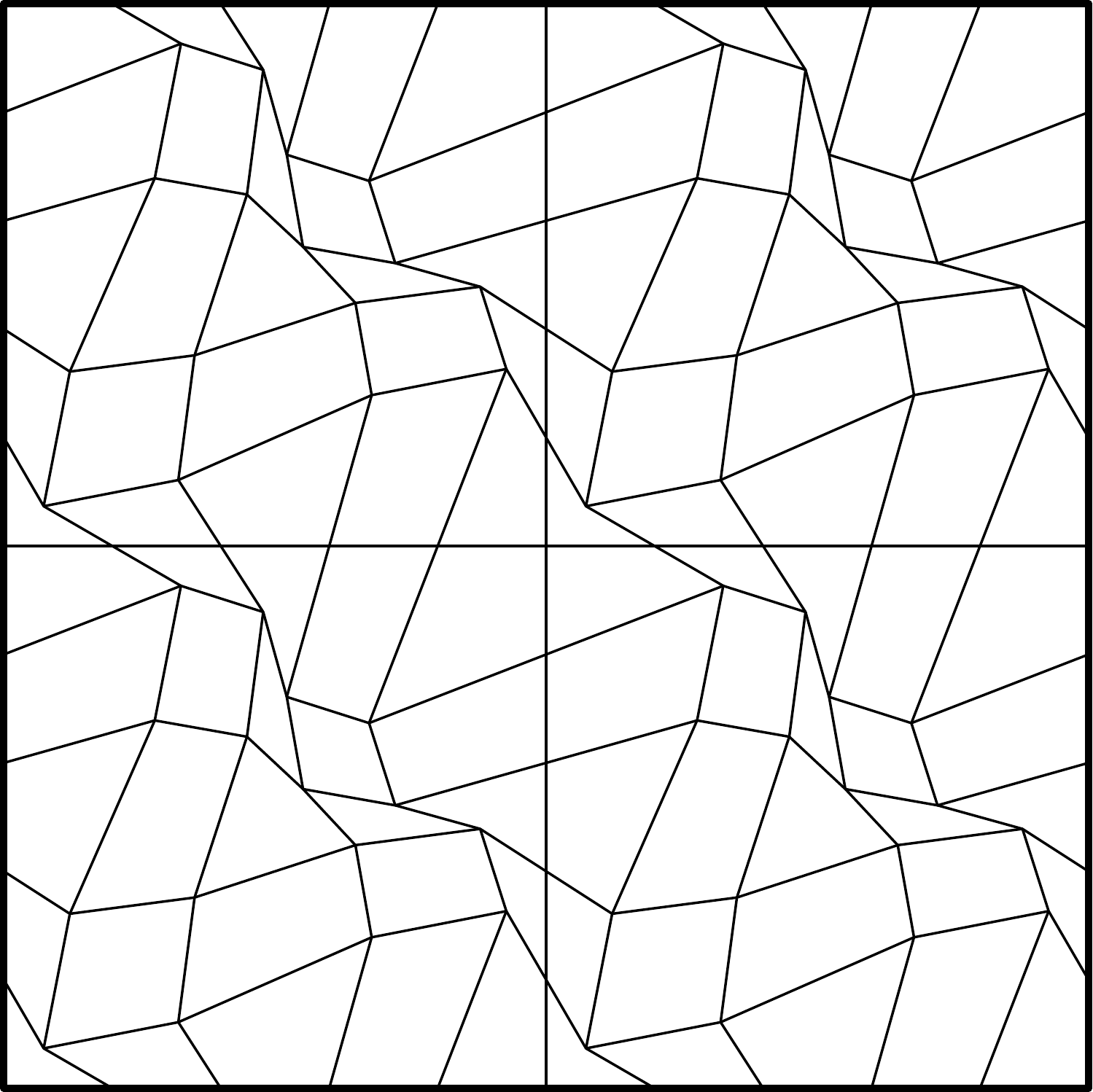}};
         \node[circle,draw=black, fill=gray!15, inner sep=0pt,minimum size=2.5pt] (b) at (0,0) {};
        \node [below] at (0,0) {$(0,0)$};
         \node[circle,draw=black, fill=gray!15, inner sep=0pt,minimum size=2.5pt] (b) at (\linewidth,\linewidth) {};
        \node [above] at (\linewidth,\linewidth) {$(a,b)$};
      \end{tikzpicture}
      \caption{\texttt{Skewed}}
      \label{fig:meshes-Skewed}
    \end{subfigure}  
    \hfill    
    \begin{subfigure}{.15\linewidth}
      \begin{tikzpicture}[scale=1.0]
        \node[anchor=south west,inner sep=0] (image) at (0,0) {\includegraphics[width=\linewidth]{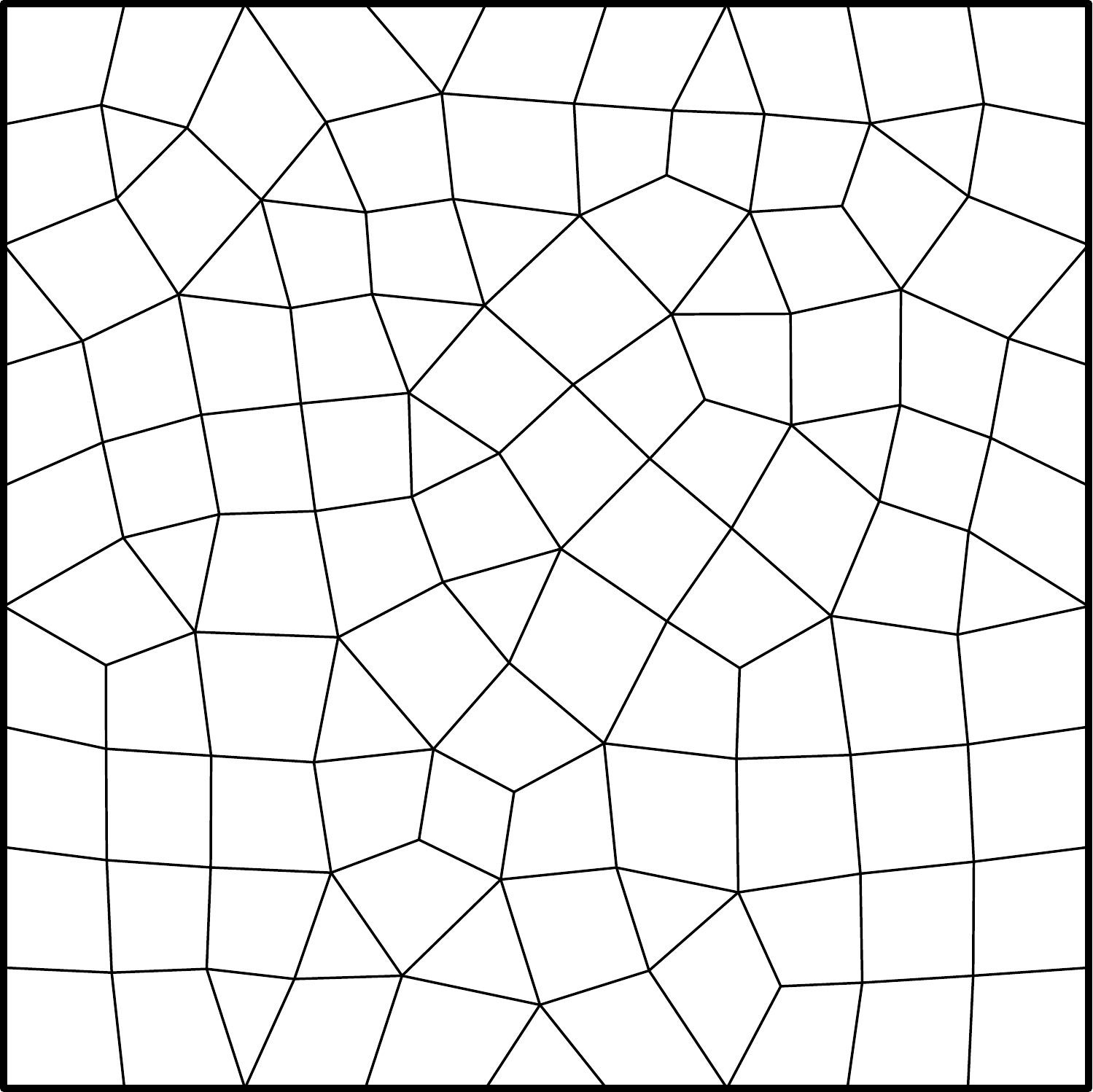}};
         \node[circle,draw=black, fill=gray!15, inner sep=0pt,minimum size=2.5pt] (b) at (0,0) {};
        \node [below] at (0,0) {$(0,0)$};
         \node[circle,draw=black, fill=gray!15, inner sep=0pt,minimum size=2.5pt] (b) at (\linewidth,\linewidth) {};
        \node [above] at (\linewidth,\linewidth) {$(a,b)$};
      \end{tikzpicture}
      \caption{\texttt{Hybrid}}
      \label{fig:meshes-Hybrid}
    \end{subfigure} 
    \hfill    
    \begin{subfigure}{.15\linewidth}
      \begin{tikzpicture}[scale=1.0]
        \node[anchor=south west,inner sep=0] (image) at (0,0) {\includegraphics[width=\linewidth]{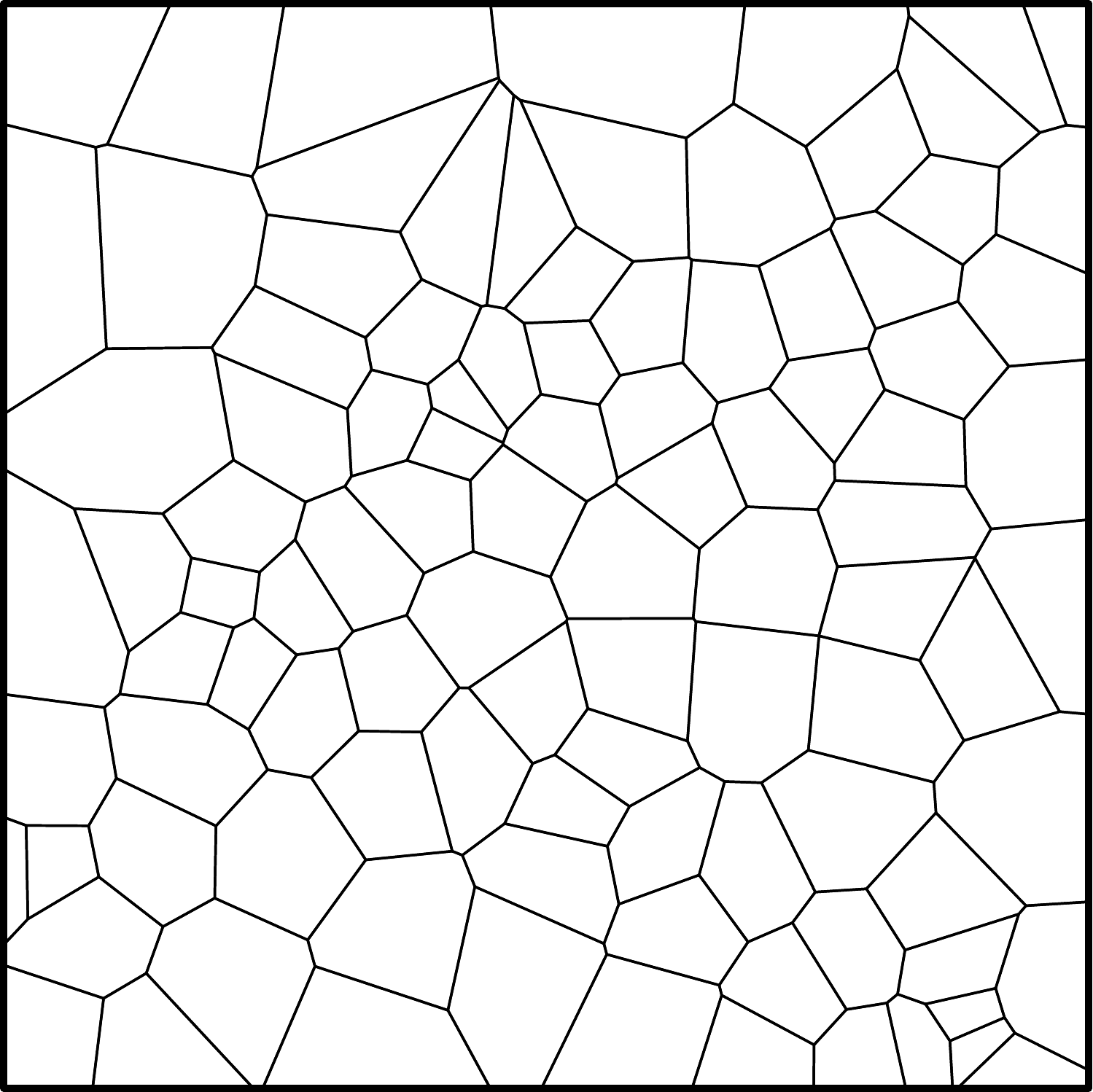}};
         \node[circle,draw=black, fill=gray!15, inner sep=0pt,minimum size=2.5pt] (b) at (0,0) {};
        \node [below] at (0,0) {$(0,0)$};
         \node[circle,draw=black, fill=gray!15, inner sep=0pt,minimum size=2.5pt] (b) at (\linewidth,\linewidth) {};
        \node [above] at (\linewidth,\linewidth) {$(a,b)$};
      \end{tikzpicture}
      \caption{\texttt{Polymesher1}}
      \label{fig:meshes-Poly1}
    \end{subfigure} 
    \hfill    
    \begin{subfigure}{.15\linewidth}
      \begin{tikzpicture}[scale=1.0]
        \node[anchor=south west,inner sep=0] (image) at (0,0) {\includegraphics[width=\linewidth]{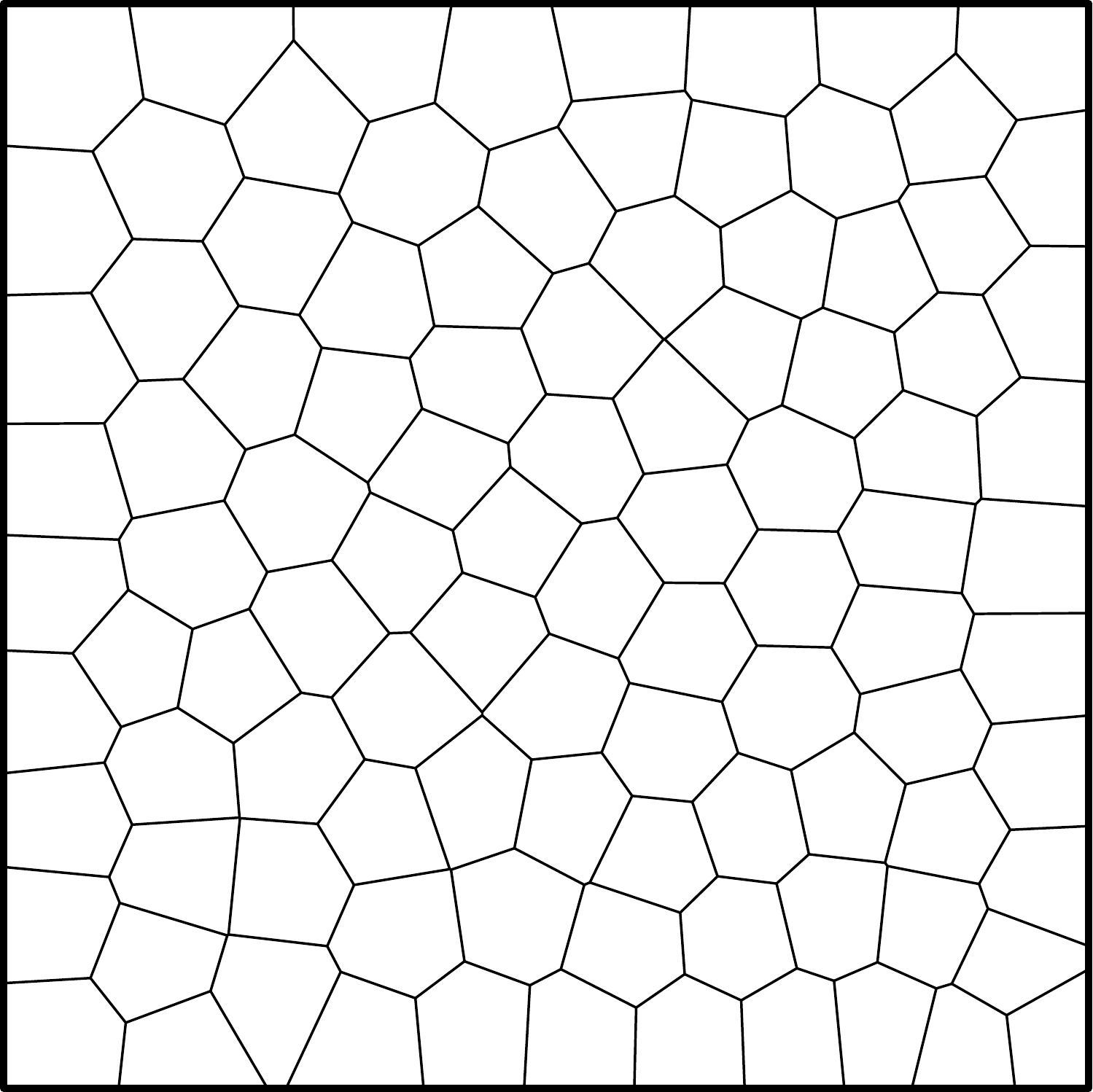}};
         \node[circle,draw=black, fill=gray!15, inner sep=0pt,minimum size=2.5pt] (b) at (0,0) {};
        \node [below] at (0,0) {$(0,0)$};
         \node[circle,draw=black, fill=gray!15, inner sep=0pt,minimum size=2.5pt] (b) at (\linewidth,\linewidth) {};
        \node [above] at (\linewidth,\linewidth) {$(a,b)$};
      \end{tikzpicture}
      \caption{\texttt{Polymesher20}}
      \label{fig:meshes-Poly20}
    \end{subfigure} 
    \hfill\null
  
  \caption{Representative meshes used in the numerical tests, with their respective tags.% \fr{I suggest that we place the meshes closer to the beginnging of the paper and that we include more refined versions. I also suggest that we replace the Cartesian mesh with a hybrid mesh here.}
  }
  \label{fig:mesh_illustration}
\end{figure}

Figure~\ref{fig:mesh_illustration} illustrates the different types of meshes considered in our numerical examples.
They include skewed meshes (\texttt{Skewed}) obtained by perturbing the vertices of a uniform Cartesian mesh (\texttt{Cartesian}) using the function:
\begin{linenomath}
\begin{equation}
g(x,y) =
\begin{pmatrix}
  x+ 0.07\sin\left(4\pi x\right)\cos\left(4\pi y+ \displaystyle \frac{\pi}{2}\right)
  \\
  y+ 0.07\sin\left(4\pi x\right)\cos\left(4\pi y+ \displaystyle \frac{\pi}{2}\right)
\end{pmatrix} \,.
\end{equation}
\end{linenomath}
The test suite also involves hybrid meshes (\texttt{Hybrid}) composed of triangles and quadrilaterals generated with GMSH \cite{geuzaine2009gmsh}, and arbitrary polygonal meshes (\texttt{Polymesher1} and \texttt{Polymesher20}) generated with respectively one and 20 smoothing steps of Polymesher \cite{Talischi2012309}.
We stress the fact that the meshes of types \texttt{Skewed}, \texttt{Hybrid}, \texttt{Polymesher1}, and \texttt{Polymesher20} do not satisfy, in general, the $\tensorTwo{\kappa}$-orthogonality condition \eqref{eq:orthogonality_condition}.

\subsection{Fully discrete coupled scheme}
\label{sec:fully_discrete_coupled_scheme}

In this section, we present the system of algebraic equations arising
from a hybrid MFD-VEM discretization of the weak form
\eqref{eq:weak_form_IBVP}.
To  simplify the presentation, we will momentarily delay the description of the construction of certain VEM and MFD operators, to first describe the overall form of the discrete equations.  A subsequent section will then revisit these important but subsidiary elements in detail.

To discretize the momentum balance equation
\eqref{momentumBalance_weak_form}, we consider a low-order VEM in
which the displacement field is approximated in the following
functional spaces:
\begin{linenomath}
\begin{subequations}
  \begin{align}
    \tensorOne{\mathcal{U}}_h &= \{ \tensorOne{u}_h \in
      C^0(\Omega) \times C^0(\Omega) \colon \, \tensorOne{u}_{h|\Gamma_u} = \overline{\tensorOne{u}}, \, 
      {\tensorOne{u}_h}_{|K} \in \mathcal{U}_{h,K}\times \mathcal{U}_{h,K}, \,
      \forall K \in \Th \}, \\ 
    %% \tensorOne{\mathcal{U}}_h &= \{ \tensorOne{u}_h \in
    %%                             [C^0(\Omega)]^2: \, \tensorOne{u}_{h|\Gamma_u} = \overline{\tensorOne{u}}, \,
    %%                             {\tensorOne{u}_h}_{|K} \in [ \mathcal{U}_{h,K} ]^2, \,
      %%                             \forall K \in \Th \}, \\
    \tensorOne{\mathcal{U}}_{h,0} &= \{ \tensorOne{u}_h \in
      C^0(\Omega) \times C^0(\Omega) \colon \, \tensorOne{u}_{h|\Gamma_u} = 0, \, 
      {\tensorOne{u}_h}_{|K} \in \mathcal{U}_{h,K}\times \mathcal{U}_{h,K}, \,
      \forall K \in \Th \}.      
    %% \tensorOne{\mathcal{U}}_{h,0} &= \{ \tensorOne{u}_h \in
    %%                                 [C^0(\Omega)]^2: \, \tensorOne{u}_{h|\Gamma_u} = 0, \,
    %%                                 {\tensorOne{u}_h}_{|K} \in [ \mathcal{U}_{h,K} ]^2, \,
    %%                                 \forall K \in \Th \}.
  \end{align}
  \label{eq:functional_spaces}\null
\end{subequations}
\end{linenomath}
The local space $\tensorOne{\mathcal{U}}_{h,K}$ is ``virtual'' in the
sense that the analytical expression of its basis functions is not
known and not needed for the construction of the scheme.
For the low-order VEM, a function of $\tensorOne{\mathcal{U}}_{h,K}$
is uniquely defined by its values at the vertices of $K$.
These vertex-based values are, for each component of the displacement,
the dofs of the VEM.
The VEM methodology provides two key operators acting on the
displacement dofs: a coercive bilinear form
$a_h( \cdot, \cdot ): \tensorOne{\mathcal{U}}_h \times
\tensorOne{\mathcal{U}}_h \rightarrow \mathbb{R}$ approximating the
$[L^2(\Omega)]^{2\times 2}$-inner product of Eq.~\eqref{momentumBalance_weak_form}, and a discrete divergence operator
discretizing the divergence term appearing in the $L^2(\Omega)$-inner
product of Eq.~\eqref{momentumBalance_weak_form}.
We denote this VEM divergence operator by
$\text{div}_{h}^{\textsc{vem}}: \tensorOne{\mathcal{U}}_h \rightarrow
\mathcal{P}_h$, where $\mathcal{P}_h$ is:
\begin{linenomath}
\begin{equation}
  \mathcal{P}_h = \{ p_h \in
  L^2(\Omega): \, 
  p_{h|K} \in \mathbb{P}_0( K ), \,
  \forall K \in \Th \}.
\end{equation}
\end{linenomath}
The definition of the local virtual space
$\tensorOne{\mathcal{U}}_{h,K}$ as well as the construction of the VEM
operators $a_h( \cdot, \cdot )$ and $\text{div}_{h}^{\textsc{vem}}$ are
reviewed in Section~\ref{sec:vem}.
The right-hand side of \eqref{momentumBalance_weak_form}
involving body forces is computed by defining, for each $K\in\Th$,
$\tilde{\tensorOne{b}}_K = \frac{1}{\abs{K}} \int_{K} \tensorOne{b}$
(where the integral of a vector is intended to be performed
component-wise) and approximating the local right-hand side as
follows:
\begin{linenomath}
\begin{equation}
  \label{eq:bodyforcerhsapprox}
  (\tensorOne{b} , \tensorOne{\eta}_h )_{[L^2(K)]^2}\approx
  (\tilde{\tensorOne{b}}_K , \tilde{\tensorOne{\eta}}_{h,K} )_{[L^2(K)]^2} = 
  \abs{K}\tilde{\tensorOne{b}}_K \cdot \tilde{\tensorOne{\eta}}_{h,K}
  \quad  \forall \tensorOne{\eta}_h \in \tensorOne{\mathcal{U}}_h\,,\text{ where }
  \tilde{\tensorOne{\eta}}_{h,K} =
  \frac{1}{\abs{K}} \int_K \tensorOne{\eta}_{h} \,.
\end{equation}
\end{linenomath}
In the low-order MFD scheme, we adopt a purely discrete representation
of the solution fields for the velocity, pressure, and Lagrange
multiplier variables.
To represent the Darcy velocity field, we denote by $\mathscr{W}_h$
the set of one-sided face-based discrete fields.
In the hybrid formulation, a discrete field
$\Vec{w}_h = ( w_{K,f} )_{K \in \mathcal{T}, f \in \mathcal{F}_K} \in
\mathscr{W}_h$ contains one dof per boundary face, and two dofs per
interior face that are not necessarily equal.
Each dof approximates the average Darcy velocity over a face such that
\begin{linenomath}
\begin{equation}
  w_{K,f} \approx \frac{1}{|f|} \int_f \tensorOne{q} \cdot \tensorOne{n}_{K,f}, \qquad \forall f \in \mathcal{F}_K, \, \forall K \in \mathcal{T}.
\end{equation}
\end{linenomath}
To discretize the weak form of the Darcy equation
\eqref{darcy_weak_form}, the MFD scheme involves a discrete weighted
inner product
$[ \cdot, \cdot ]_{\mathscr{W}_h}: \mathscr{W}_h \times \mathscr{W}_h
\rightarrow \mathbb{R}$ used to approximate the
$[L^2(\Omega)]^2$-inner product of Eq.~\eqref{darcy_weak_form}.
For the discretization of the divergence term present in the $L^2(\Omega)$-inner
product of Eq.~\eqref{darcy_weak_form}, we also define a discrete
divergence operator,
$\text{div}_{h}^{\textsc{mfd}}: \mathscr{W}_h \rightarrow
\mathscr{P}_h$, where the set of cell-based discrete fields, $\mathscr{P}_h$, is isomorphic to the space of piecewise constant functions, $\mathcal{P}_h$.
The construction of these two MFD operators is detailed in
Section~\ref{sec:mhfv}.
%
% The MFD discretization involves the definition of a discrete
% divergence operator
% $\text{div}_h: \mathscr{W}_h \rightarrow \mathscr{P}_h$ such that:
%
%\begin{equation}
%\text{div}_h \Vec{w}_h = \bigg( \frac{1}{|K|} \sum_{f \in \mathcal{F}_K} |f| w_{K,f} \bigg)_{K \in \mathcal{T}} \in \mathscr{P}_h.
%\end{equation}
%
%A key element of the MFD is the construction of the discrete weighted
% inner product
% $[ \cdot, \cdot ]_{\mathscr{W}_h}: \mathscr{W}_h \times
% \mathscr{W}_h \rightarrow \mathbb{R}$ to approximate the first
% bilinear function in \eqref{massBalance_weak_form}.
%
% The principles guiding the construction of this inner product are
% similar to those employed to design $a_h( \cdot, \cdot )$ in VEM.

The pressure field is represented as a cell-based discrete field of
$\mathscr{P}_h$.
A discrete pressure solution
$\Vec{p}_h = ( p_{K} )_{K \in \mathcal{T}} \in \mathscr{P}_h$ is a
collection of dofs approximating the cell-based pressure averages,
i.e.,
\begin{linenomath}
\begin{equation}
  p_K \approx \frac{1}{|K|} \int_K p, \qquad \forall K \in \mathcal{T}.
\end{equation}
\end{linenomath}
Equivalently, the pressure solution can be viewed as a
piecewise-constant function of $\mathcal{P}_h$ using the decomposition
$p_h( \tensorOne{x} ) = (I^{\mathcal{P}_h} \Vec{p}_h)(\tensorOne{x})
= \sum_{K \in \mathcal{T}} p_{K} \chi_K( \tensorOne{x} )$, where
$I^{\mathcal{P}_h}: \mathscr{P}_h \rightarrow \mathcal{P}_h$ is the interpolation operator and $(\chi_K)_{K \in \mathcal{T}}$ is the canonical basis of $\mathcal{P}_h$.
The discrete counterpart of the $L^2( \Omega )$-inner product is
$[ \cdot, \cdot ]_{\mathscr{P}_h}: \mathscr{P}_h \times \mathscr{P}_h
\rightarrow \mathbb{R}$, classically defined as the sum of local inner
products:
\begin{linenomath}
\begin{equation}
  [ \Vec{p}_h, \Vec{\chi}_h ]_{\mathscr{P}_h} = \sum_{K \in \mathcal{T}} [ \Vec{p}_{h|K}, \Vec{\chi}_{h|K} ]_{\mathscr{P}_{h|K}} = \sum_{K \in \mathcal{T}} |K| p_K \chi_K.
\end{equation}
\end{linenomath}

In the hybrid formulation, the computation of the Darcy velocity
involves Lagrange multipliers that are represented as a face-based
discrete field.
We remark that, considering the set of face-based discrete fields,
$\mathscr{L}_h$, the field
$\Vec{\pi}_h = (\pi_f)_{f \in \mathcal{F}} \in \mathscr{L}_h$ contains
one dof per face that can be viewed as an approximation of the face
pressure average, i.e.,
\begin{linenomath}
\begin{equation}
  \pi_f \approx \frac{1}{|f|} \int_f p, \qquad \forall f \in \mathcal{F}.
\end{equation}
\end{linenomath}
The subset $\mathscr{L}_{h,0}$ contains the discrete fields
$\Vec{\pi}_h = ( \pi_f )_{f}$ satisfying $\pi_f = 0$ for
$f \in \mathcal{F}_p$.

We consider a fully implicit (backward-Euler) temporal discretization
of the coupled system.
The superscript $n$ denotes the time level at which the degrees of
freedom are evaluated.
Using the notations introduced above, the coupled problem in discrete
weak form reads:
given two functions $\{ \tensorOne{u}^0, p^0\}$ defining the initial
state, find the function
$\tensorOne{u}^n_h \in \tensorOne{\mathcal{U}}_h$ and the discrete
fields
$\{ \Vec{w}^n_h, \Vec{p}^n_h, \Vec{\pi}^n_h \} \in \mathscr{W}_h
\times \mathscr{P}_h \times \mathscr{L}_h$ such that for
$n \in \{1, \dots, N\}$:
\begin{linenomath}
\begin{subequations}
  \begin{align}
    & a_h( \tensorOne{u}^n_h, \tensorOne{\eta}_h ) - ( \alpha \, I^{\mathcal{P}_h} \Vec{p}^n_h, \text{div}_{h}^{\textsc{vem}} \, \tensorOne{\eta}_h )_{L^2(\Omega)} = \sum_{K\in\Th} \abs{K}\tilde{\tensorOne{b}}_K \cdot \tilde{\tensorOne{\eta}}_{h,K} + ( \overline{\tensorOne{t}}, \tensorOne{\eta}_h )_{[L^2(\Gamma_{\sigma})]^2}, & &\forall \tensorOne{\eta}_h \in \tensorOne{\mathcal{U}}_{h,0}, \label{discrete_momentum_equation}\\
    &[\Vec{w}^n_h, \Vec{\varphi}_h ]_{\mathscr{W}_h} - [ \Vec{p}^n_h, \text{div}_{h}^{\textsc{mfd}}  \, \Vec{\varphi}_h ]_{\mathscr{P}_h} + \sum_{K \in \mathcal{T}} \sum_{f \in \mathcal{F}} |f| \pi^n_f \varphi_{K,f} = 0,& &\forall \Vec{\varphi}_h \in \mathscr{W}_h, \label{discrete_darcy_equation} \\ 
    &( \alpha \, \text{div}_{h}^{\textsc{vem}} \, \tensorOne{u}^n_h, I^{\mathcal{P}_h} \Vec{\chi}_h )_{L^2(\Omega)} + \Delta t \, [ \text{div}_{h}^{\textsc{mfd}} \, \Vec{w}^n_h, \Vec{\chi}_h ]_{\mathscr{P}_h} + [ S_{\varepsilon} \Vec{p}^n_h, \Vec{\chi}_h ]_{\mathscr{P}_h} = s^{n-1}_h,&  &\forall \Vec{\chi}_h \in \mathscr{P}_h, \label{discrete_mass_balance_equation} \\
    &- \sum_{K \in \mathcal{T}} \sum_{f \in \mathcal{F}_K} |f| w^n_{K,f} \lambda_f  = - \sum_{f \in \mathcal{F}_q} |f| \overline{q} \lambda_f,&  &\forall \Vec{\lambda}_h \in \mathscr{L}_{h,0}, \label{discrete_face_constraints}
  \end{align}
  \label{eq:discrete_coupled_scheme}\null
\end{subequations}
\end{linenomath}
where the right-hand side of \eqref{discrete_mass_balance_equation}
is:
\begin{linenomath}
\begin{equation}
  s^{n-1}_h = ( \alpha \, \text{div}_{h}^{\textsc{vem}} \, \tensorOne{u}^{n-1}_h, I^{\mathcal{P}_h} \Vec{\chi}_h )_{L^2(\Omega)} + [ S_{\varepsilon} \Vec{p}^{n-1}_h, \Vec{\chi}_h ]_{\mathscr{P}_h}.
  \label{eq:discrete_right_hand_side_at_time_n_1}
\end{equation}
\end{linenomath}
In the hybridized system \eqref{eq:discrete_coupled_scheme}, the use
of Lagrange multipliers allows a local, cell-wise computation of
the one-sided face velocities in the discrete Darcy equation
\eqref{discrete_darcy_equation}.
To ensure that the hybrid scheme remains mass conservative, the set of
algebraic constraints \eqref{discrete_face_constraints} imposes the
continuity of the velocities at the mesh faces.
We stress the fact that, despite a relatively large number of dofs,
the hybridized system is amenable to static condensation, which is used to locally eliminate the one-sided face velocities during the
assembly.
The resulting algebraic system solved by the linear solver is discussed in Section~\ref{sec:solution_strategy}.

In the following sections, we focus on the terms of
\eqref{eq:discrete_coupled_scheme} that have not been fully defined
yet.
In Section~\ref{sec:vem}, we define the virtual space
$\tensorOne{\mathcal{U}}_{h}$ and review the construction of the VEM
operators $a_h( \cdot, \cdot )$ and $\text{div}_{h}^{\textsc{vem}}$.
In Section~\ref{sec:mhfv}, we show that a similar methodology is used
in MFD to form the operators $[ \cdot, \cdot ]_{\mathscr{W}_h}$ and
$\text{div}_{h}^{\textsc{mfd}}$.
For simplicity, we drop the superscript $n$ denoting the time level.

% mainfile: main.tex
\subsection{Local virtual space and VEM operators}
\label{sec:vem}

We start the section with the local virtual space introduced in
the low-order VEM \cite{Beirao2015a,Beirao2015b} to approximate the
displacement variable.
%
% Equation \eqref{momentumBalanceS} is discretized by the Virtual
% Element Method \cite{Beirao2015a,Beirao2015b}.
%
Let $K\in\Th$ be a polygon of the tessellation of $\Omega$.
Following \cite{Beirao2015a}, we define the following set of scaled
monomials on $K$:
\begin{linenomath}
  \begin{equation}
    \Monom{1}{K} = \left\{m_1( \tensorOne{x} ) = 1,\,
      m_2( \tensorOne{x} ) = \frac{x-x_K}{h_K},\,
      m_3( \tensorOne{x} ) = \frac{y-y_K}{h_K}\right\},
  \end{equation}
\end{linenomath}
where we recall that $\tensorOne{x}_K = (x_K,y_K)$ is the center of
$K$ and $h_K$ is the diameter of $K$.
We define the projection operator
$\proj[\nabla,K]{1}{} \colon \sobh{1}{K} \to \Poly{1}{K}$ such that,
$\forall \eta\in\sobh{1}{K}$,
\begin{linenomath}
  \begin{equation}
    \label{eq:defPiNabla}
    \begin{cases}
      \displaystyle \scal[K]{\nabla \proj[\nabla,K]{1}{} \eta}{\nabla
        m} = \scal[K]{\nabla \eta}{\nabla m}, & \forall
      m\in\Monom{1}{K},
      \\
      \displaystyle \int_{\partial K} \proj[\nabla,K]{1}{}\eta =
      \int_{\partial K} \eta \,.
    \end{cases}
  \end{equation}
\end{linenomath}
Then, we define the following functional space:
\begin{linenomath}
  \begin{equation}
    \label{eq:defVh}
    \mathcal{U}_{h,K} = \left\{ \eta\in\sobh{1}{K}\colon \Delta \eta \in\Poly{1}{K},\,
      \eta\in C^0\left(\partial K\right),\,
      \eta_{|f}\in\Poly{1}{f}\;\forall f\in\mathcal{F}_K,\,
      \scal[K]{\eta}{m} = \scal[K]{\proj[\nabla,K]{1}{}\eta}{m}\;
      \forall m\in\Monom{1}{K} \right\},
  \end{equation}
\end{linenomath}
where $\Poly{1}{K}$ is the set of polynomials of degree $\leq 1$
defined on $K$, and $\Poly{1}{f}$ the set of polynomials of degree
$\leq 1$ defined on the face $f$.

A function $\eta_h\in \mathcal{U}_{h,K}$ is completely defined by its
values at the vertices of $K$ \cite{Beirao2015b}.
Functions in $\mathcal{U}_{h,K}$ are called ``virtual'' because their
analytical expressions is not known. Instead, we know that they are
polynomials of degree $\leq 1$ on each edge. For each
$\eta_h\in\mathcal{U}_{h,K}$, its projection
$\proj[\nabla,K]{1}{}\eta_h$ is computable from the degrees of
freedom. Indeed, the right-hand sides of \eqref{eq:defPiNabla} can be
computed knowing only the analytical expression of $\eta_h$ on
$\partial K$. In particular, the first right-hand side is computed by
integrating by parts and applying Green's theorem. Exploiting the fact
that $\Delta m = 0$ $\forall m\in\Poly{1}{K}$, we get
\begin{linenomath}
  \begin{equation}
    \scal[K]{\nabla \eta_h}{\nabla m} = \scal[\partial K]{\eta_h}{\pder{m}{n}} \,.
  \end{equation}
\end{linenomath}
The last condition on the functions in $\mathcal{U}_{h,K}$ (see
Eq.~\eqref{eq:defVh}) allows us to compute the integral mean of
$\eta_h$ on $K$. Indeed,
\begin{linenomath}
  \begin{equation}
    \label{eq:vembasismean}
    \tilde{\eta}_{h,K} := \frac{1}{\abs{K}}\int_K \eta_h =
    \frac{1}{\abs{K}}\int_K \proj[\nabla,K]{1}{}\eta_h \,.
  \end{equation}
\end{linenomath}
Finally, the integral mean of the gradient of $\eta_h$ on $K$ is also
computable knowing only the analytical expression of $\eta_h$ on
$\partial K$:
\begin{linenomath}
  \begin{align}
    \widetilde{\partial}_{x,K} \eta_h := \frac{1}{\abs{K}}\int_K \pder{\eta_h}{x}
    &= -\frac{1}{\abs{K}} \int_{\partial K}\eta_h n_x \,,
    &
      \widetilde{\partial}_{y,K} \eta_h :=  \frac{1}{\abs{K}}\int_K \pder{\eta_h}{y}
    &= -\frac{1}{\abs{K}} \int_{\partial K}\eta_h n_y \,.
  \end{align}
\end{linenomath}
We denote
$\widetilde{\nabla}_K \eta_h := \begin{pmatrix} \widetilde{\partial}_{x,K} \eta_h \\
  \widetilde{\partial}_{y,K} \eta_h \end{pmatrix}$. We use the same
notation in the following for vectorial and tensorial functions, where
integral means are performed on each component.

With reference to the notation introduced in Section
\ref{sec:fully_discrete_coupled_scheme}, we discretize the
displacement $\displ^n_h$ defining the space
\begin{linenomath}
  \begin{equation}
    \tensorOne{\mathcal{U}}_h = \{ \tensorOne{\eta}_h \in
      C^0(\Omega) \times C^0(\Omega) \colon 
      {\tensorOne{\eta}_h}_{|K} \in \mathcal{U}_{h,K}\times \mathcal{U}_{h,K}, \,
      \forall K \in \Th \}.
  \end{equation}
\end{linenomath}
Furthermore, the discrete bilinear form of the problem is obtained by
first defining the following discretizations of the strains and
effective stresses, for each $K\in \Th$ and each $\tensorOne{\eta}_h =
\begin{pmatrix}
  \eta_{h,x}\\ \eta_{h,y}
\end{pmatrix}
\in \tensorOne{\mathcal{U}}_h$:
\begin{linenomath}
  \begin{align}
    \label{eq:defVemStrain}
    \widetilde{\strain}_K\left(\tensorOne{\eta}_h\right)
    &=\frac{1}{\abs{K}}\int_K\strain\left(\tensorOne{\eta}_h\right) =
      \frac{1}{2}\left( \widetilde{\nabla}_K\tensorOne{\eta}_h +
      \left( \widetilde{\nabla}_K\tensorOne{\eta}_h \right)\transpose\right),
    \\
    \label{eq:defVemStress}
    \widetilde{\stress}^\prime_K\left(\tensorOne{\eta}_h\right)
    &= 2G\widetilde{\strain}_K\left(\tensorOne{\eta}_h\right) +
      \lambda \text{trace}
      \left( \widetilde{\strain}_K\left(\tensorOne{\eta}_h\right) \right),
  \end{align}
\end{linenomath}
where $\widetilde{\nabla}_K\tensorOne{\eta}_h =
\begin{pmatrix}
  \left( {\widetilde{\nabla}_K\eta_{h,x}} \right)\transpose
  \\
  \left( {\widetilde{\nabla}_K\eta_{h,y}} \right)\transpose
\end{pmatrix}
$. Moreover, we define the discrete VEM divergence operator
$\text{div}_{h}^{\textsc{vem}}: \tensorOne{\mathcal{U}}_h \rightarrow
\mathcal{P}_h$ as
\begin{linenomath}
  \begin{equation}
    \text{div}_{h}^{\textsc{vem}}\tensorOne{\eta}_h =
    \text{trace}\left(\widetilde{\nabla}_K\tensorOne{\eta}_h\right), \quad
    \forall \tensorOne{\eta}_h \in \tensorOne{\mathcal{U}}_h,
  \end{equation}
\end{linenomath}
that is used in \eqref{discrete_momentum_equation} and
\eqref{eq:discrete_right_hand_side_at_time_n_1}.

To complete the definition of VEM discrete bilinear forms, following
\cite{Beirao2015a}, we define, for each $K\in\Th$, the continuous and
coercive bilinear form
$\ah[K]{}{}\colon
\tensorOne{\mathcal{U}}_h\times\tensorOne{\mathcal{U}}_h \to
\mathbb{R}$ such that
\begin{linenomath}
  \begin{equation}
    \label{eq:defVemBilForm}
    \ah[K]{\tensorOne{u}_h}{\tensorOne{\eta}_h} =
    \int_K \widetilde{\stress}^\prime_K\left(\tensorOne{u}_h\right) \colon
    \widetilde{\strain}_K\left(\tensorOne{\eta}_h\right) +
    \left(2\sup_{K}G\right)
    \vemstab[K]{\tensorOne{u}_h-\proj[\nabla,K]{1}{}\tensorOne{u}_h}
    {\tensorOne{\eta}_h-\proj[\nabla,K]{1}{}\tensorOne{\eta}_h},
    \quad \forall \tensorOne{u}_h, \tensorOne{\eta}_h \in \tensorOne{\mathcal{U}}_h.
  \end{equation}
\end{linenomath}
The bilinear form $\vemstab[K]{}{}$ is defined in order to be
computable from the degrees of freedom of the two functions involved
and to ensure the coercivity of $\ah[K]{}{}$. The most common choice,
and the one we make here, is
\begin{linenomath}
  \begin{equation}
    \label{eq:defVemStab}
    \vemstab[K]{\tensorOne{u}_h-\proj[\nabla,K]{1}{}\tensorOne{u}_h}
    {\tensorOne{\eta}_h-\proj[\nabla,K]{1}{}\tensorOne{\eta}_h}
    = \sum_{v \in \mathcal{V}_K}
    \left(
      \tensorOne{u}_h(\tensorOne{x}_v) -
      \proj[\nabla,K]{1}{\tensorOne{u}_h}(\tensorOne{x}_v)
    \right)
    \left(
      \tensorOne{\eta}_h(\tensorOne{x}_v) -
      \proj[\nabla,K]{1}{\tensorOne{\eta}_h}(\tensorOne{x}_v)
    \right), \,
  \end{equation}
\end{linenomath}
where we recall that $\mathcal{V}_K$ denotes the set of the vertices
of $K$. With the above definitions, we can build the discrete VEM
bilinear form, that is used in \eqref{discrete_momentum_equation}, as
\begin{linenomath}
  \begin{equation}
    \ah{\tensorOne{u}_h}{\tensorOne{\eta}_h} = \sum_{K\in\Th}
    \ah[K]{\tensorOne{u}_h}{\tensorOne{\eta}_h}, \quad \forall
    \tensorOne{u}_h, \tensorOne{\eta}_h \in \tensorOne{\mathcal{U}}_h.
  \end{equation}
\end{linenomath}
Moreover, Eq.~\eqref{eq:vembasismean} applied to each component of a test
function $\tensorOne{\eta}_h \in \tensorOne{\mathcal{U}}_h$, allows us
to compute the approximation of the right-hand side term involving the
body force as described in \eqref{eq:bodyforcerhsapprox}. Finally, we
note that, since the VEM basis functions are known to be polynomials on
faces, we can compute the right-hand side of
\eqref{discrete_momentum_equation} involving
$\overline{\tensorOne{t}}$ as we would do for a classical finite
element discretization.

\begin{rmk}
  \label{rmk:decomposition_consistency_stability}
  The bilinear form $\ah[K]{}{}$ defined by \eqref{eq:defVemBilForm}
  is composed of two terms: the first one (first addend in
  \eqref{eq:defVemBilForm}) accounts for the consistency of the scheme
  on polynomials of degree $1$, the second one (second addend in
  \eqref{eq:defVemBilForm}, defined by \eqref{eq:defVemStab}) accounts
  for coercivity. Indeed, whenever either
  $\tensorOne{\eta}_h\in[\mathbb{P}_1(K)]^2$ or
  $\tensorOne{u}_h \in [\mathbb{P}_1(K)]^2$,
  $\vemstab[K]{\tensorOne{u}_h-\proj[\nabla,K]{1}{}\tensorOne{u}_h}
  {\tensorOne{\eta}_h-\proj[\nabla,K]{1}{}\tensorOne{\eta}_h} = 0$,
  since the operator $\proj[\nabla,K]{1}{}$ is the identity for
  polynomials of degree $1$. Moreover, since the operators
  $\widetilde{\stress}^\prime_K$ and $\widetilde{\strain}_K$ are
  consistent on polynomials of degree $1$ (see \eqref{eq:defVemStrain}
  and \eqref{eq:defVemStress}), we have that $\ah[K]{}{}$ is exact on
  $[\mathbb{P}_1(K)]^2 \times [\mathbb{P}_1(K)]^2$. The same structure
  can be observed for the MFD bilinear form
  $[ \cdot\,, \cdot ]_{\mathscr{W}_h}$, that is described in Section
  \ref{sec:mhfv} (see \eqref{eq:generic_form_of_M}).
\end{rmk}

%%% Local Variables:
%%% mode: latex
%%% TeX-master: "main"
%%% End:

\subsection{MFD operators}
\label{sec:mhfv}

In this section, we focus on the two key discrete MFD operators---namely, the divergence operator and the weighted inner product---acting on the set of one-sided face-based fields $\mathscr{W}_h$ that appear in the discrete Darcy and mass balance equations \eqref{discrete_darcy_equation}-\eqref{discrete_mass_balance_equation}.
Let $\Vec{w}_h$, $\Vec{q}_h \in \mathscr{W}_h$ be two discrete fields representing one-sided face velocities.
The discrete divergence operator, $\text{div}_{h,MFD}: \mathscr{W}_h \rightarrow \mathscr{P}_h$, is defined by:
\begin{equation}
\text{div}_{h}^{\textsc{mfd}} \Vec{w}_h = \bigg( \frac{1}{|K|} \sum_{f \in \mathcal{F}_K} |f| w_{K,f} \bigg)_{K \in \mathcal{T}} \in \mathscr{P}_h.
\end{equation}
The weighted inner product $[ \cdot, \cdot ]_{\mathscr{W}_h}: \mathscr{W}_h \times \mathscr{W}_h \rightarrow \mathbb{R}$ is written as the sum of local inner products:
\begin{equation}
  [ \Vec{w}_h, \Vec{q}_h ]_{\mathscr{W}_h} = \sum_{K \in \mathcal{T}} [ \Vec{w}_{h|K}, \Vec{q}_{h|K} ]_{\mathscr{W}_{h|K}} = \sum_{K \in \mathcal{T}} \sum_{f \in \mathcal{F}_K} \sum_{f' \in \mathcal{F}_K} ( \Mat{M}_K )_{ff'} w_{K,f} q_{K,f}.
\end{equation}
The mimetic discretization framework provides a methodology to construct a symmetric positive definite matrix of size $\text{card}(\mathcal{F}_K) \times \text{card}(\mathcal{F}_K)$, denoted by $\Mat{M}_K$, that only depends on the permeability and geometric properties of cell $K$.
As in Section~\ref{sec:vem}, an expression for $\Mat{M}_K$ is obtained by imposing consecutively a consistency condition and a stability condition \cite{da2014mimetic}, which can then be used to demonstrate that the scheme is convergent for elliptic problems in mixed form \cite{brezzi2005convergence}.

For a low-order mimetic finite-difference scheme, the consistency condition states that the local inner product $[ \cdot, \cdot ]_{\mathscr{W}_{h|K}}$ must be exact when one of the two arguments is the projection on $\mathscr{W}_{h|K}$ of a constant function $\tensorOne{\kappa}_{|K} \cdot \nabla q$ with $q \in \mathbb{P}_1(K)$.
%
%Specifically, considering a function $\tensorOne{w} \in [L^2(K)]^2$ with a constant divergence on $K$ and a constant dot product $\tensorOne{w} \cdot \tensorOne{n}_{K,f}$ for all faces $f \in \mathcal{F}_{K}$, the inner product must be accurate in the following sense:
%\begin{equation}
%[ \Pi^{\mathscr{W}_{h|K}} ( \tensorOne{\kappa}_{|K} \cdot \nabla q ), \Pi^{\mathscr{W}_{h|K}} \tensorOne{w} ]_{\mathscr{W}_{h|K}} = \int_K \tensorOne{\kappa}^{-1}_{|K} \cdot ( \tensorOne{\kappa}_{|K} \cdot \nabla q \cdot \tensorOne{w} ), \quad \forall q \in \mathbb{P}_1( K ), \label{eq:consistency_condition}
%\end{equation}
%where the discrete projection operator $\Pi^{\mathscr{W}_{h|K}}: [L^2(K)]^2 \rightarrow \mathscr{W}_{h|K}$ is defined by $\Pi^{\mathscr{W}_{h|K}} \tensorOne{w} = \big( 1/|f| \int_f \tensorOne{w} \cdot \tensorOne{n}_{K,f} \big)_{f \in \mathcal{F}_K}$ for all $\tensorOne{w} \in [L^2(K)]^2$.
%
Setting $q$ to $q_1( \tensorOne{x} ) = x - x_K$ and $q_2( \tensorOne{x} ) = y - y_K$ in the consistency condition results in the following algebraic constraint on the entries of $\Mat{M}_K$ \cite{da2014mimetic}:
\begin{equation}
  \Mat{M}_K \Mat{N}_K = \Mat{R}_K.
  \label{eq:algebraic_consistency_condition}
\end{equation}
Let $f_i$ denote the $i$-th face in $\mathcal{F}_K$.
By construction, $\Mat{N}_K$ and $\Mat{R}_K$ are two full-rank matrices of size $\text{card}(\mathcal{F}_K) \times 2$ such that the $i$-th row of $\Mat{N}_K$ is $\tensorOne{n}\transpose_{K,f_i} \tensorOne{\kappa}$, and the $i$-th row of $\Mat{R}_K$ is $|f_i| \tensorOne{c}\transpose_{K,f_i}$.
Condition \eqref{eq:algebraic_consistency_condition} defines a family of consistent low-order mimetic finite-difference schemes, but does not provide a unique expression for $\Mat{M}_K$.
It is worth noting that the symmetric matrix $\frac{1}{|K|} \Mat{R}_K \tensorTwo{\kappa}^{-1} \Mat{R}\transpose_K$ satisfies \eqref{eq:algebraic_consistency_condition}, but is only positive semidefinite.
To enforce the SPD structure of $\Mat{M}_K$ and obtain a coercive bilinear form, a stabilization term is introduced.
Let $\Mat{C}_K$ be a matrix of size $\text{card} (\mathcal{F}_K) \times (\text{card}( \mathcal{F}_K) - 2)$ whose columns form a basis of $\text{ker}( \Mat{N}\transpose_K)$, and let $\Mat{U}_K$ be a SPD matrix of size $(\text{card} (\mathcal{F}_K) - 2) \times (\text{card}( \mathcal{F}_K) - 2)$.
As shown in \cite{brezzi2005family}, writing $\Mat{M}_K$ in the generic form:
\begin{equation}
  \Mat{M}_K = \frac{1}{|K|} \Mat{R}_K \tensorTwo{\kappa}^{-1}  \Mat{R}\transpose_K +  \gamma_K \widetilde{\Mat{U}}_K \quad \text{with} \quad \widetilde{\Mat{U}}_K := \Mat{C}_K \Mat{U}_K \Mat{C}\transpose_K,
  \label{eq:generic_form_of_M}
\end{equation}
yields a coercive local bilinear form that still satisfies the algebraic consistency condition \eqref{eq:algebraic_consistency_condition} since $\Mat{C}\transpose_K \Mat{N}_K = 0$.
In \eqref{eq:generic_form_of_M}, it is clear that the structure of the mimetic inner product is analogous to that of the VEM operator $\ah[K]{}{}$ as both bilinear forms can be split into a term ensuring consistency and a term enforcing stability (see Remark~\ref{rmk:decomposition_consistency_stability}).
Multiple definitions of the stabilization term in \eqref{eq:generic_form_of_M} have been proposed in previous work \cite{lie2019introduction}.
Here, we obtain an inner product satisfying the MFD stability condition by setting the scaling coefficient as:
\begin{equation}
  \gamma_K := \frac{1}{ \text{card}(\mathcal{F}_K) |K| } \text{trace}( \Mat{R}_K \tensorOne{\kappa}_{|K}^{-1} \Mat{R}\transpose_K ),
  \label{eq:definition_gamma}
\end{equation}
and by writing $\Mat{U}_K = ( \Mat{C}\transpose_K \Mat{C}_K )^{-1}$, which yields after some simplifications \cite{da2014mimetic}:
\begin{equation}
  \widetilde{\Mat{U}}_K = \Mat{I} - \Mat{N}_K ( \Mat{N}\transpose_K \Mat{N}_K )^{-1} \Mat{N}\transpose_K.
  \label{eq:definition_u_tilde}
\end{equation}
%Defining $\Mat{M}_K$ as above garantees that the resulting inner product satisfies the mimetic finite-difference stability condition, which states that there exists $(\sigma_{K})_{*} > 0$ and $(\sigma_K)^{*} > 0$ dependent on the rock and geometric properties of cell $K$ such that:
%\begin{equation}
%(\sigma_{K})_{*} \sum_{f \in \mathcal{F}_K} |f| w_{K,f} \leq [ \Vec{w}, \Vec{w} ]_{\mathscr{W}_{h|K}} \leq (\sigma_K)^{*} \sum_{f \in \mathcal{F}_K} |f| w_{K,f}, \quad \forall \Vec{w} \in \mathscr{W}_{h|K}. \label{eq:stability_condition}
%\end{equation}
%
%Provided that the mesh satisfies the assumptions of Section~\ref{sec:mesh_hierarchy}, equations \eqref{eq:generic_form_of_M}, \eqref{eq:definition_gamma}, and \eqref{eq:definition_u_tilde} yield a convergent mimetic finite-difference scheme for elliptic problems in mixed form \cite{brezzi2005convergence}.

\begin{rmk}
  We note that setting:
  \begin{equation}
  (\Mat{M}_{K})_{ff'} =
  \begin{dcases}
       \frac{|f| || \tensorOne{c}_{K,f} ||_2^2}{\tensorOne{n}_{K,f} \cdot \tensorTwo{\kappa}_{|K} \cdot \tensorOne{c}_{K,f}},
 &\text{if } f' = f \\
     0 &\text{otherwise,}
  \end{dcases} \label{tpfaHalfTrans}
  \end{equation}
  yields a matrix containing the inverses of the standard TPFA half-transmissibilities divided by $|f|^2$.
  This can be used to recover the cell-centered linear TPFA scheme after algebraic elimination of the face variables \cite{beaude2019combined}.
  However, it is well known that the linear TPFA scheme is not consistent when the mesh does not satisfy the $\tensorTwo{\kappa}$-orthogonality condition \eqref{eq:orthogonality_condition}.
\end{rmk}

\subsection{Local pressure-jump stabilization}
  \label{sec:stability}

  The stability of the coupled scheme described in
  \eqref{eq:discrete_coupled_scheme} is subject to the well-known
  Ladyzhenskaya-Babu{\v{s}}ka-Brezzi (LBB) \textit{inf-sup} stability
  condition \cite{babuvska1971error,brezzi1974existence}:
  \begin{equation}
    \exists \beta>0 \colon \inf_{p_h \in \mathcal{P}_h \setminus \{ 0 \} } \sup_{\tensorOne{v}_h \in \tensorOne{\mathcal{U}}_h \setminus \{ \tensorOne{0} \} } \frac{ ( \alpha p_h, \text{div}_{h}^{\textsc{vem}} \tensorOne{v}_h )_{L^2(\Omega)} }{ \norm[L^2(\Omega)]{p_h} \left\| \tensorOne{v}_h \right\|_{[L^2(\Omega)]^2} } \geq \beta.
    \label{lbb_stability_condition}
  \end{equation}
  The proposed scheme relies on approximation spaces for the
  displacement and pressure variables that do not, in
  general, satisfy the discrete \textit{inf-sup} condition given in
  \eqref{lbb_stability_condition}.
  Specifically, for incompressible solid and fluid constituents,
  i.e. $S_{\varepsilon} = 0$, in the presence of undrained
  conditions---resulting for instance from small time steps and/or low
  permeability---the instability may manifest itself by the presence
  of spurious modes in the pressure field (checkerboarding).

  It is worth noting that, for specific mesh topologies, the spurious
  pressure modes do not appear and no stabilization is required.
  In particular, the numerical examples of Section~\ref{sec:cantilever}
  indicate that the \texttt{Polymesher1} and \texttt{Polymesher20} mesh
  families do not exhibit any checkerboarding.
  These results are in agreement with the findings of \cite{da2009mimetic,da2010mimetic,talischi2014polygonal},
  in which equal-order interpolation schemes applied to the Stokes problem
  are shown to be stable whenever no vertex in the mesh is connected to
  more than three faces.

\begin{figure}
  \small
  
  \hfill
  \begin{subfigure}{.275\linewidth}
    \centerline{\includegraphics[width=\linewidth]{./figs/base_hybrid}}
    \caption{}
  \end{subfigure}
  \hfill
  \begin{subfigure}{.275\linewidth}
    \centerline{\includegraphics[width=\linewidth]{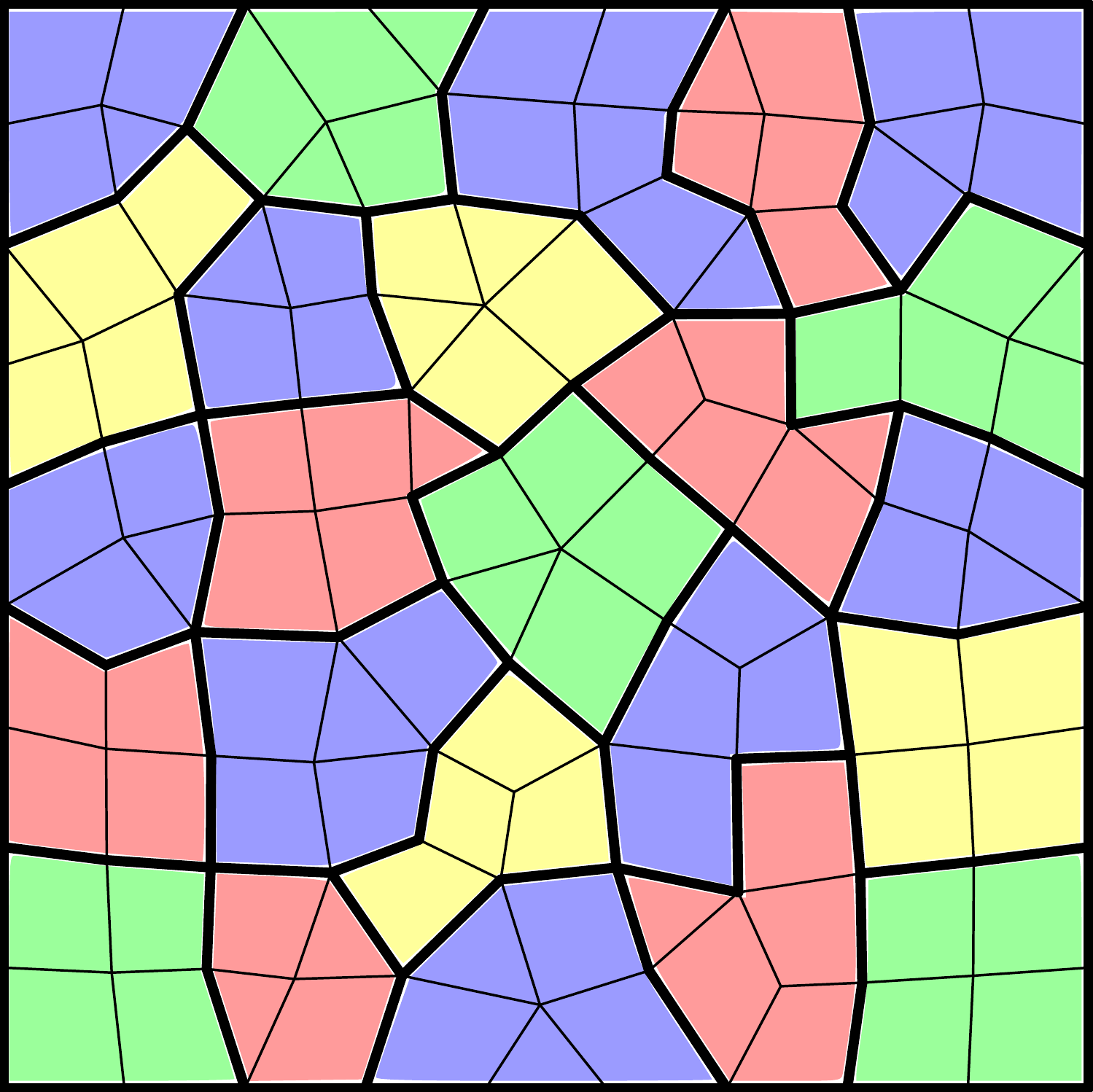}}
    \caption{}
  \end{subfigure}
  \hfill
  \begin{subfigure}{.275    \linewidth}
    %\centerline{\includegraphics[width=\linewidth]{./figs/macro_face}}
      \begin{tikzpicture}[scale=1.0]
        \node[anchor=south west,inner sep=0] (image) at (0,0) {\includegraphics[width=\linewidth]{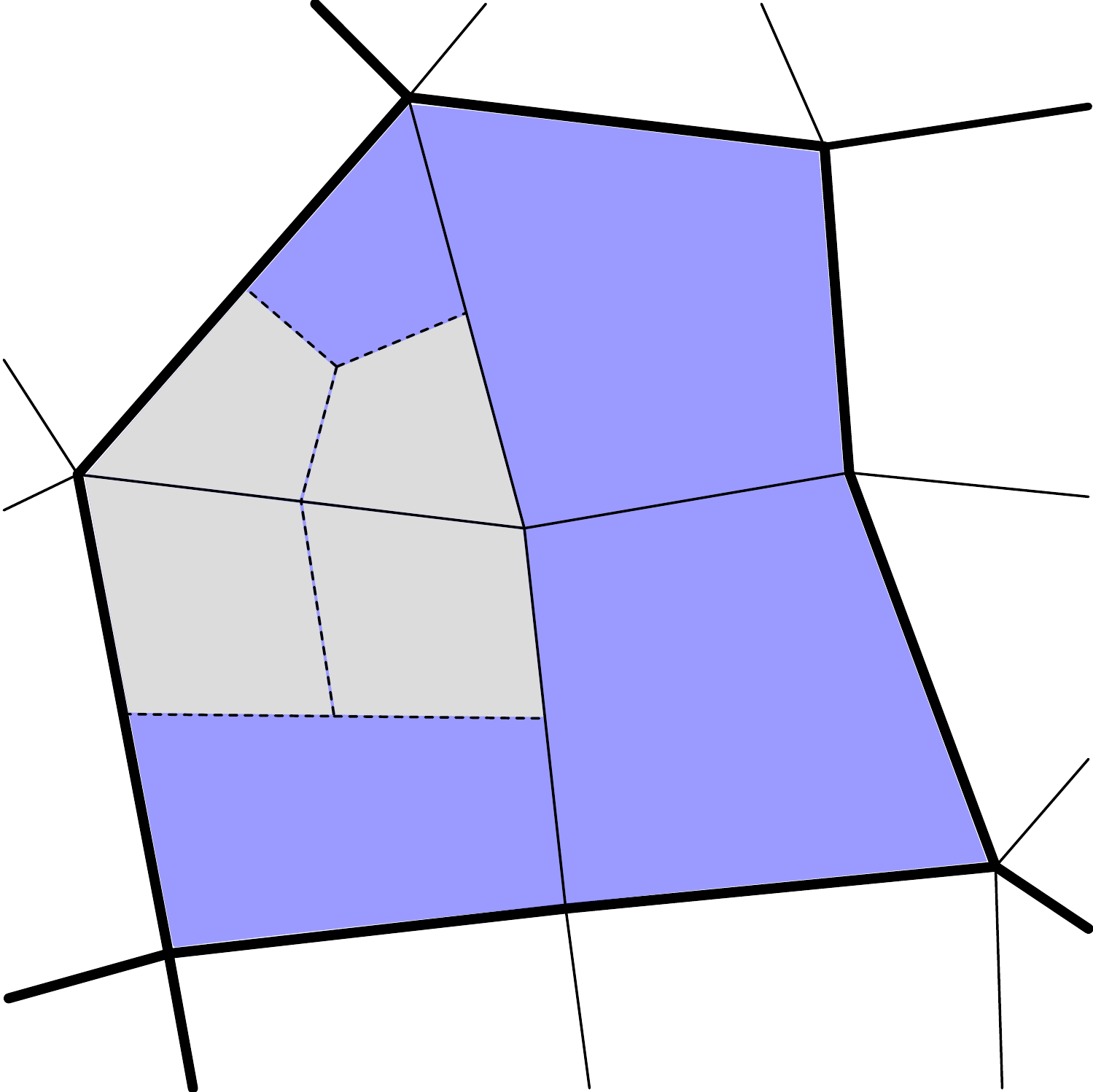}};
         \node[circle,draw=black, fill=gray!15, inner sep=0pt,minimum size=2.5pt] (b) at (0.0725\linewidth,0.5675\linewidth) {};
        \node [left] at (0.0725\linewidth,0.585\linewidth) {$v_2$};
         \node[circle,draw=black, fill=gray!15, inner sep=0pt,minimum size=2.5pt] (b) at (.482\linewidth,.5175\linewidth) {};
        \node [right] at (0.482\linewidth,0.48\linewidth) {$v_1$};
        \node [right] at (.25\linewidth,.25\linewidth) {$K_1$};
        \node [right] at (.6\linewidth,.35\linewidth) {$K_2$};
        \node [right] at (.525\linewidth,.7\linewidth) {$K_3$};
        \node [right] at (.28\linewidth,.775\linewidth) {$K_4$};
        \node [right] at (.11\linewidth,.44\linewidth) {$m_{1,2}$};
        \node [right] at (.31\linewidth,.42\linewidth) {$m_{1,1}$};
        \node [right] at (.3\linewidth,.59\linewidth) {$m_{4,1}$};
        \node [right] at (.14\linewidth,.62\linewidth) {$m_{4,2}$};
      \end{tikzpicture}    
    \caption{}
  \end{subfigure}
  \hfill\null
    
  \caption{\label{fig:macro_element_construction}Unstructured macro-element stabilization for the \texttt{Hybrid} mesh family: (a) mesh; (b) mesh partition into macro-elements; and (c) example of macro-element $E = K_1 \cup K_2 \cup K_3 \cup K_4$, made of three quadrilaterals ($K_1$, $K_2$, and $K_3$) and one triangle ($K_4$). In the stabilized scheme, we add artificial fluxes at the four (internal) faces adjacent to vertex $v_1$. Assuming that $f$ is the face defined by vertices $v_1$ and $v_2$, we compute $\Upsilon_f := m_{1,1} + m_{1,2} + m_{4,1} + m_{4,2}$ in Eq.~\eqref{eq:stabilization_term_in_mass_conservation} as the area of the region in gray.}
\end{figure}  
  
  To eliminate checkerboarding in the unstable mesh configurations---in
  this work, the \texttt{Cartesian}, \texttt{Skewed}, \texttt{Hybrid} mesh
  types---we rely on the local pressure-jump stabilization technique
  introduced originally for the Stokes problem \cite{kechkar1992analysis,silvester1994optimal}, and more specifically,
  we adapt the methodology used in \cite{camargo2020macroelement,frigo2020efficient}
  to unstructured meshes. 
  To achieve this, we partition the mesh into non-overlapping cell aggregates referred
  to as \emph{macro-elements} using the algorithm presented in \ref{sec:macro_element_construction}.
  We emphasize that these macro-elements are unstructured, in contrast to the traditional concept of structured macro-elements on a logically-nested grid.
  Considering the set of macro-elements $\mathcal{E}$, a macro-element $E \in \mathcal{E}$
  is constructed such that $E = \cup_{K \in E} K$ and
  $\cup_{E \in \mathcal{E}} E = \mathcal{T}$.
  We denote by $\mathcal{F}_{E, \textit{int}}$ the set of mesh faces that are
  in the interior of $E$.
  Following \cite{camargo2020macroelement,frigo2020efficient}, the local pressure-jump
  stabilization consists in introducing artificial fluxes at the internal faces
  $f \in \mathcal{F}_{E, \textit{int}}$ of each macro-element $E$ to prevent the development of spurious
  pressure modes.
  This is done by adding the following term in the mass balance equation \eqref{discrete_mass_balance_equation}:
  \begin{equation}
    J( \Vec{p}^n_h, \Vec{\chi}_h ) := \beta \sum_{E \in \mathcal{E}} \sum_{f \in \mathcal{F}_{E, \textit{int}}} \Upsilon_f [\![ \Vec{p}^n_h ]\!]_f [\![ \Vec{\chi}_h ]\!]_f, \qquad \Vec{p}^n_h, \Vec{\chi}_h \in \mathscr{P}_h, \label{eq:stabilization_term_in_mass_conservation}
  \end{equation}
  where the pressure jump for an internal face $f \in \mathcal{F}_{E, \textit{int}}$ is
  $[\![ \Vec{p}^n_h ]\!]_f := p^n_K - p^n_L$ ($K, L \in \mathcal{T}_f$ with $K < L$).
  %
  %For each internal face, we introduce a geometric coefficient, $\Upsilon_f$, representing a characteristic area computed as follows.
  We introduce a geometric coefficient, $\Upsilon_f$, representing a characteristic area associated with an internal face $f \in \mathcal{F}_{E, \textit{int}}$ and defined as follows.
  Considering cell $K \in E$, we associate each node $v \in \mathcal{V}_K$ with the measure, denoted by $m_{K,v}$, of the quadrilateral whose vertices (appropriately ordered) are $v$, the centers of the two faces adjacent to $v$ in $K$, and the centroid of $K$ (see Fig.~\ref{fig:macro_element_construction}).
  Using that, we define $\Upsilon_f$ as:
  \begin{equation}
    \Upsilon_f := \sum_{K \in \mathcal{T}_f} \sum_{v \in \mathcal{V}_f} m_{K,v}.
  \end{equation}
  The two-dimensional stabilization coefficient \cite{elman2014finite,frigo2020efficient}, denoted by $\beta$, is computed
  using Biot's coefficient and the Lam{\'e} parameters as:
  \begin{equation}
    \beta := \frac{\alpha^2}{4(2 G + \lambda)}.
  \end{equation}
  Since the artificial fluxes are only added at the internal faces of the macro-elements, the stabilized scheme is not cell-wise mass conservative but remains mass conservative at the level of the macro-elements---in the sense that the sum of the fluxes at the external faces of a macro-element is equal to zero for the incompressible setting considered here.
  The stability properties of the scheme and their impact on the performance of the iterative
  linear solver are assessed in Section~\ref{sec:cantilever}.

%\begin{figure}
%  \centering
%  \scalebox{1.25}{
%    \input{./face_characteristic_area.tex}
%  }
%  \caption{Macro-element $E = K_1 \cup K_2 \cup K_3 \cup K_4$ appearing for the \texttt{hybrid} mesh family, made of two quadrilaterals ($K_1$ and $K_2$) and two triangles ($K_3$ and $K_4$). In the stabilized scheme, we add artificial fluxes at the four (internal) faces adjacent to vertex $v_1$. Assuming that $f$ is the face in red in this schematic, we compute $\Upsilon_f := m_{1,1} + m_{1,2} + m_{4,1} + m_{4,2}$ in equation \eqref{eq:stabilization_term_in_mass_conservation} as the area of the region in blue.}
%  \label{fig:macro_element_construction_old}
%\end{figure}

\subsection{Solution strategy}
\label{sec:solution_strategy}

The discrete weak form \eqref{eq:discrete_coupled_scheme} produces a sequence of $4\times4$ block systems of algebraic equations of the type
\begin{linenomath}
\begin{align}
  \begin{bmatrix}
    \Mat{A}\sub{uu}
    & 0
    & -\Mat{A}\sub{up}
    & 0
    \\
    0
    & \Mat{A}\sub{ww}
    & -\Mat{A}\sub{wp}
    & -\Mat{A}\sub{w\pi}
    \\
    \Mat{A}\sub{up}\transpose
    & \Delta t \Mat{A}\sub{wp}\transpose
    & \overline{\Mat{A}}\sub{pp}
    & 0
    \\
    0  
    & \Mat{A}\sub{w \pi}\transpose
    & 0
    & 0
  \end{bmatrix}
  \begin{bmatrix}
    \Vec{u} \\ \Vec{w} \\ \Vec{p} \\ \Vec{\pi}
  \end{bmatrix}
  =
  \begin{bmatrix}
    \Vec{b}_{u} \\
    \Vec{b}_{w} \\
    \Vec{b}_{p} \\
    \Vec{b}_{\pi}
  \end{bmatrix},
  \label{eq:4x4_system}
\end{align}
\end{linenomath}
where
\begin{itemize}
  \item $\Mat{A}\sub{uu}$ is a symmetric positive definite (SPD) matrix corresponding to the elasticity block;
  \item $\Mat{A}\sub{ww}$ is a block-diagonal SPD matrix, having as many blocks as the number of cells, each one with a size equal to the number of faces of the corresponding cell; 
  \item $\overline{\Mat{A}}\sub{pp}$ consists of two contributions: (i) a diagonal matrix that depends cell-wise on the specific storage coefficient $S_{\varepsilon}$, and (ii) a local stabilization contribution (Eq. \eqref{eq:stabilization_term_in_mass_conservation}), if needed, whose sparsity pattern is a subset of that of a discrete (cell-centered) Laplace operator;
  \item $\Mat{A}\sub{up}$ and $\Mat{A}\sub{wp}$ are matrices generated by the inner products $( \alpha \, I^{\mathcal{P}_h} \Vec{\chi}^n_h, \text{div}_{h}^{\textsc{vem}} \, \tensorOne{\eta}_h )_{L^2(\Omega)}$ and $[ \Vec{\chi}^n_h, \text{div}_{h}^{\textsc{mfd}}  \, \Vec{\varphi}_h ]_{\mathscr{P}_h}$, respectively;
  \item $\Mat{A}\sub{w\pi}$ is a matrix whose columns correspond to a unique mesh interface, having two (respectively one) negative unit entries for interfaces belonging to $\mathcal{F}_{\text{int}}$ (respectively $\mathcal{F}_p \cup \mathcal{F}_q$).
\end{itemize}
Following standard practice \cite{Bof_etal08}, we take advantage of the block-diagonal structure of $\Mat{A}\sub{ww}$ by reducing \eqref{eq:4x4_system} via static condensation to a $3 \times 3$ linear system

\begin{linenomath}
\begin{align}
  \blkMat{A} \blkVec{x} &= \blkVec{b},
  &
  \blkMat{A} &= 
  \begin{bmatrix}
    \Mat{A}\sub{uu}
    & -\Mat{A}\sub{up}
    & 0
    \\
    \Mat{A}\sub{up}\transpose
    & \Mat{A}\sub{pp}
    & \Delta t \Mat{A}\sub{p \pi}
    \\
    0
    & \Mat{A}\sub{p \pi}\transpose
    & \Mat{A}\sub{\pi \pi}
  \end{bmatrix},
  &
  \blkVec{x} &=
  \begin{bmatrix}
    \Vec{u} \\ \Vec{p} \\ \Vec{\pi}
  \end{bmatrix},
  &
  \blkVec{b} &=
  \begin{bmatrix}
    \Vec{b}_{u} \\
    \Vec{b}_{p} - \Delta t \Mat{A}\sub{wp}\transpose \Mat{A}\sub{ww}^{-1} \Vec{b}\sub{w} \\
    \Vec{b}_{\pi} - \Mat{A}\sub{w \pi}\transpose \Mat{A}\sub{ww}^{-1} \Vec{b}\sub{w}
  \end{bmatrix},
  \label{eq:3x3_system}
\end{align}
\end{linenomath}
with
$\Mat{A}\sub{pp}
 = (\overline{\Mat{A}}\sub{pp}
 + \Delta t \Mat{A}\sub{wp}\transpose \Mat{A}\sub{ww}^{-1} \Mat{A}\sub{wp}$),
$\Mat{A}\sub{p \pi}
 = (\Mat{A}\sub{wp}\transpose \Mat{A}\sub{ww}^{-1} \Mat{A}\sub{w \pi})$, and
$\Mat{A}\sub{\pi \pi}
 = (\Mat{A}\sub{w \pi}\transpose \Mat{A}\sub{ww}^{-1} \Mat{A}\sub{w \pi})$.
We note that $\Mat{A}\sub{pp}$ shares the same sparsity pattern as $\overline{\Mat{A}}\sub{pp}$ since the matrix arising from the static condensation, i.e. $\Delta t \Mat{A}\sub{wp}\transpose \Mat{A}\sub{ww}^{-1} \Mat{A}\sub{wp}$, is diagonal.

The non-symmetric linear system \eqref{eq:3x3_system} is solved by a preconditioned Krylov subspace method. 
%%
%To warrant that mechanics and flow discrete equations have comparable magnitude, we apply a simple block-row and block-column scaling before entering the Krylov solver yielding
%%
%\begin{linenomath}
%\begin{align}
%  \blkMat{A} \blkVec{y} &= \blkVec{b},
%  &
%  \blkMat{A} &=
%  \blkMat{D}^{-1/2} \overline{\blkMat{A}} \blkMat{D}^{-1/2} =
%  \begin{bmatrix}
%    \Mat{A}\sub{uu}
%    & -\Mat{A}\sub{up}
%    & 0
%    \\
%    \Mat{A}\sub{up}\transpose
%    & \Mat{A}\sub{pp}
%    & \Delta t \Mat{A}\sub{p \pi},
%    \\
%    0
%    & \Mat{A}\sub{p \pi}\transpose
%    & \Mat{A}\sub{\pi \pi}
%  \end{bmatrix},  
%  &
%  \blkVec{y} &= \blkMat{D}^{1/2} \blkVec{x}
%  &
%  \blkVec{b} &= \blkMat{D}^{-1/2} \overline{\blkVec{b}},
%  \label{eq:3x3_system_scaled}
%\end{align}
%\end{linenomath}
%%
%where
%%
%\begin{linenomath}
%\begin{align}
%  \blkMat{D} &=
%  \begin{bmatrix}
%    \theta_u I & 0 & 0 \\
%    0 & \theta_p I & 0 \\
%    0 & 0 & \theta_{\pi} I
%  \end{bmatrix},
%  &
%  \theta_u &= || \Mat{A}\sub{uu} ||_{F},
%  &
%  \theta_p &= \theta_{\pi} = \left\| \begin{bmatrix} \Mat{A}\sub{pp} & \Delta t\Mat{A}\sub{p\pi} \\ \Mat{A}\sub{p\pi}\transpose & \Mat{A}\sub{\pi\pi} \end{bmatrix} \right\|_{F},
%  \label{eq:scaling}
%\end{align}
%\end{linenomath}
%%
%with $||\cdot||_F$ denoting the Frobenius norm of a matrix.
%%
It shares the same properties of the Biot system addressed in \cite{frigo2020efficient}, which is obtained based on a mixed hybrid finite element formulation of the same three-field formulation considered in this work.
Thus, we adopt the same block-triangular preconditioning strategy, used in conjunction with  a right-preconditioned generalized minimal residual (GMRES) method \cite{SaaSch86}.
Denoting by $\blkMat{P}^{-1}$ the preconditioning operator, we work with the modified system:
\begin{linenomath}
\begin{subequations}
  \begin{align}  
    &\blkMat{A} \blkMat{P}^{-1} \blkVec{y} = \blkVec{b}, \\
    &\blkVec{x} = \blkMat{P}^{-1} \blkVec{y}.
  \end{align}
\end{subequations}
\end{linenomath}
The preconditioner in factorized form reads:
\begin{linenomath}
\begin{align}
  \blkMat{P}^{-1}
%  = 
%  \begin{bmatrix}
%    \widetilde{\Mat{A}}\sub{uu} &
%    -\Mat{A}\sub{up} &
%    0 \\
%    0 &
%    \widetilde{\Mat{B}}\sub{pp} &
%    \Delta t \Mat{A}\sub{p\pi} \\
%    0 &
%    0 &
%    \widetilde{\Mat{C}}\sub{\pi \pi}
%  \end{bmatrix}^{-1}
  =
  \begin{bmatrix}
    \widetilde{\Mat{A}}\sub{uu}^{-1} & 0 &  0 \\ 
    0 & I & 0 \\
    0 & 0 & I   
  \end{bmatrix}
  \begin{bmatrix}
    I & \Mat{A}\sub{up} &  0 \\ 
    0 & I & 0 \\
    0 & 0 & I   
  \end{bmatrix}
  \begin{bmatrix}
    I & 0 &  0 \\ 
    0 & \widetilde{\Mat{B}}\sub{pp}^{-1} & 0 \\
    0 & 0 & I   
  \end{bmatrix}
  \begin{bmatrix}
    I & 0 &  0 \\ 
    0 & I & -\Delta t \Mat{A}\sub{p \pi} \\
    0 & 0 & I   
  \end{bmatrix}
  \begin{bmatrix}
    I & 0 &  0 \\ 
    0 & I & 0 \\
    0 & 0 & \widetilde{\Mat{C}}\sub{\pi \pi}^{-1}  
  \end{bmatrix},
  \label{eq:blk_precond}
\end{align}
\end{linenomath}
with $\widetilde{\Mat{A}}\sub{uu}^{-1}$ a suitable approximation of $\Mat{A}\sub{uu}^{-1}$, $\widetilde{\Mat{B}}\sub{pp}^{-1}$ a suitable approximation of the inverse of the first-level SPD Schur complement $\Mat{B}\sub{pp} = ( \Mat{A}\sub{pp} + \Mat{A}\sub{up}\transpose \Mat{A}\sub{uu}^{-1} \Mat{A}\sub{up})$, and  $\widetilde{\Mat{C}}\sub{\pi\pi}^{-1}$ a suitable approximation of the inverse of the second-level SPD Schur complement $\Mat{C}\sub{\pi\pi} = ( \Mat{A}\sub{\pi\pi} - \Delta t \Mat{A}\sub{p\pi}\transpose \Mat{B}\sub{pp}^{-1} \Mat{A}\sub{p \pi})$.
We define $\widetilde{B}\sub{pp}^{-1}$ considering the following expression:
\begin{linenomath}
\begin{align}
  \Mat{B}\sub{pp} &\approx \Mat{A}\sub{pp} + \text{diagm} (\Mat{A}\sub{up}\transpose \Mat{D}\sub{uu}^{-1} \Mat{A}\sub{up}), 
  \label{eq:sparse_schur1}
\end{align}
\end{linenomath}
where the so-called fixed-stress assumption  \cite{KimTchJua11,WhiCasTch16,CasWhiFer16} is used to introduce a sparse approximation of the triple product $\Mat{A}\sub{up}\transpose \Mat{A}\sub{uu}^{-1} \Mat{A}\sub{up}$.
In \eqref{eq:sparse_schur1}, $\text{diagm}( \cdot )$ is an operator constructing a diagonal matrix by extracting the diagonal entries of the input matrix, and $\Mat{D}\sub{uu} = \text{diagm}( \Mat{A}\sub{uu})$.
Based on \eqref{eq:sparse_schur1}, the action of $\widetilde{\Mat{B}}\sub{pp}^{-1}$ on a vector is always expressed by a single $\ell_1$-Jacobi iteration \cite{Baker11,VasYan14}.
Therefore, $\widetilde{\Mat{B}}\sub{pp}^{-1}$ is available explicitly and allows us to replace, for preconditioning purposes, the exact second-level Schur complement with 
\begin{linenomath}
\begin{align}
  \Mat{C}\sub{\pi\pi} \approx \Mat{A}\sub{\pi\pi} - \Delta t \Mat{A}\sub{p\pi}\transpose \widetilde{\Mat{B}}\sub{pp}^{-1} \Mat{A}\sub{p \pi}.
  \label{eq:sparse_schur2}
\end{align}
\end{linenomath}
The approximations considered for $\widetilde{\Mat{A}}\sub{uu}^{-1}$ and $\widetilde{\Mat{C}}\sub{\pi \pi}^{-1}$ are provided in Section \ref{sec:numerical_examples}.
For a detailed analysis of the features of preconditioner \eqref{eq:blk_precond} the reader is referred to \cite{frigo2020efficient}.

\section{Numerical examples}
\label{sec:numerical_examples}

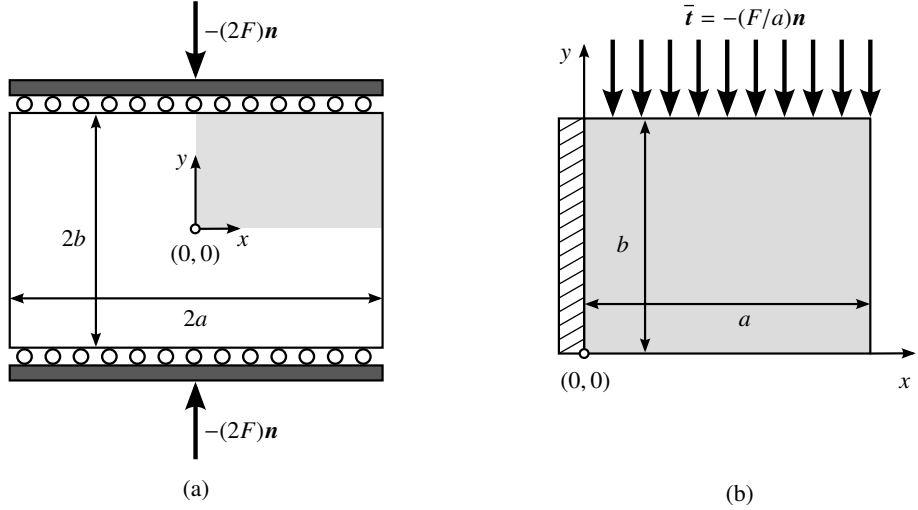
\begin{figure}[htbp]
  \small
  
  \hfill
  \begin{subfigure}{.3\linewidth}
      \begin{tikzpicture}
        \node[anchor=south west,inner sep=0] (image) at (0,0) {\includegraphics[width=\linewidth]{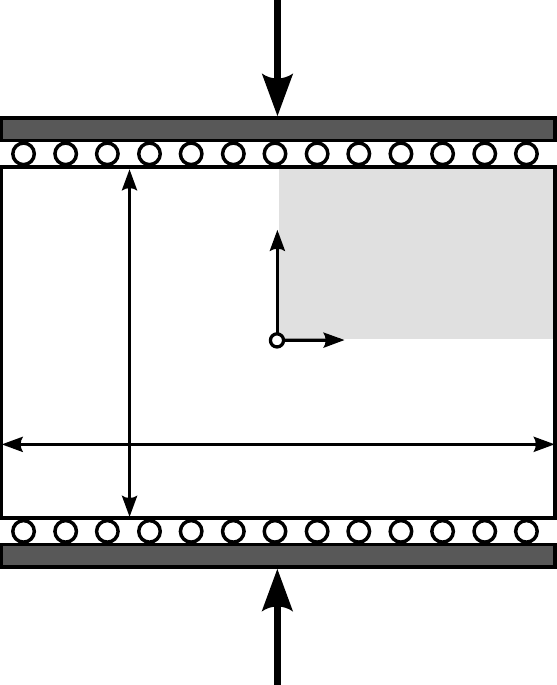}};
        \node[above] at (.5\linewidth,.5\linewidth) {$(0,0)$};
        \node[above] at (.625\linewidth,.55\linewidth) {$x$};
        \node[left] at (.5\linewidth,.8\linewidth) {$y$};
        \node[right] at (.5\linewidth,.075\linewidth) {$-(2F) \tensorOne{n}$};
        \node[right] at (.5\linewidth,1.15\linewidth) {$-(2F) \tensorOne{n}$};            
        \node[below] at (.5\linewidth,.425\linewidth) {$2a$}; 
        \node[above] at (.175\linewidth,.55\linewidth) {$2b$};       
      \end{tikzpicture}
    \caption{}
    \label{fig:mandel_sketch}
  \end{subfigure}
  \hfill
  \begin{subfigure}{.3\linewidth}
      \begin{tikzpicture}
        \node[anchor=south west,inner sep=0] (image) at (0,0) {\includegraphics[width=\linewidth]{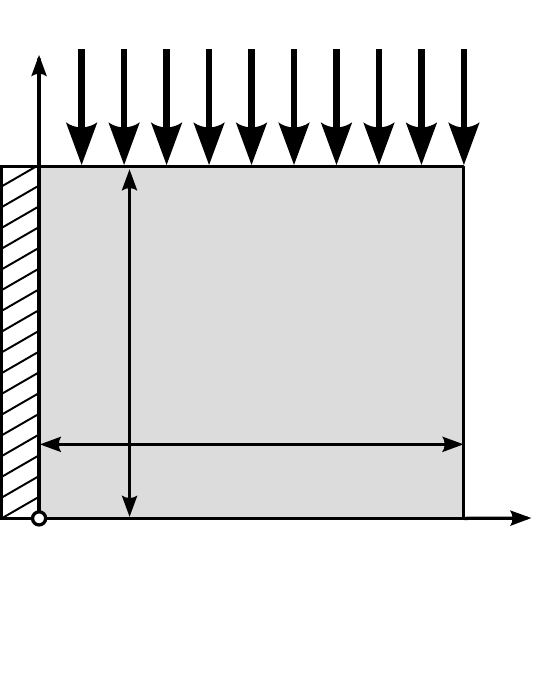}};
        \node[below] at (.075\linewidth,.275\linewidth) {$(0,0)$};
        \node[above] at (.925\linewidth,.185\linewidth) {$x$};
        \node[left] at (.06\linewidth,1.1\linewidth) {$y$};
        \node[above] at (.5\linewidth,1.15\linewidth) {$\overline{\tensorOne{t}} = - (F/a) \tensorOne{n}$};            
        \node[below] at (.5\linewidth,.425\linewidth) {$a$}; 
        \node[above] at (.175\linewidth,.55\linewidth) {$b$};       
      \end{tikzpicture}
    \caption{}
    \label{fig:cantilever_sketch}
  \end{subfigure}
  \hfill\null
  \caption{Domain sketch for Mandel's problem (a) and the cantilevered square block (b). For Mandel's problem, symmetry allows for modeling only the gray region shown in (a), i.e. a quarter of the domain.}
\end{figure}

We now consider three sets of numerical experiments.
In the first set, Mandel's problem (Fig.~\ref{fig:mandel_sketch}), a classical benchmark of linear poroelasticity for which an analytical solution is available \cite{Mandel1953}, is used to validate our MFD-VEM formulation.
In the second set, we numerically investigate the scheme's convergence properties based on a manufactured regular solution \cite{BofBotDiP16}, solving the IBVP \eqref{eq:IBVP} on the unit square domain $\Omega = [0,1]^2$.
In the third set, a cantilevered square block (Fig.~\ref{fig:cantilever_sketch}) is considered to demonstrate the robustness of the proposed stabilization with respect to pressure oscillations in the incompressible limit. In particular, we emphasize its beneficial effects on the iterative solver convergence.
Details of the material and parameter values used in the three test cases are summarized in Table \ref{tab:parameters}.

In all tests, we assume the zero vector as initial guess for (non-restarted) GMRES and terminate the iterations when the initial residual has been reduced by a factor of $10^6$.
As to the preconditioner \eqref{eq:blk_precond}, the operator $\tilde{\Mat{A}}\sub{uu}^{-1}$ is defined based on the so-called separate displacement component part of the stiffness matrix \cite{AxeGus78,Bla94}, i.e. a sparse approximation to $\Mat{A}\sub{uu}$ in which $x$- and $y$-displacement dofs are decoupled.
We consider either a sparse direct solver or algebraic multigrid (AMG) preconditioning for both $\tilde{\Mat{A}}\sub{uu}^{-1}$ and $\tilde{\Mat{C}}\sub{\pi\pi}^{-1}$.
Specifically, we use a classic AMG method \cite{RugStu87} as provided by the \texttt{HSL\_MI20} package \cite{BoyMihSco10} with default parameters (symmetric Gauss-Seidel smoother, single V-cycle, and a direct coarse solver) except for the coarsening failure criterion control (\texttt{c\_fail}), which is set to 2.
An $\ell_1$-Jacobi \cite{Baker11,VasYan14} smoother is always employed for $\widetilde{\Mat{B}}\sub{pp}^{-1}$.

\begin{table}[htbp]
	\caption{Parameter values used for the numerical examples.}
	\label{tab:parameters}
    \small
    \begin{tabular}{llllllll}
    \toprule
    Symbol & Parameter & Units &
 \begin{tabular}{@{}l@{}}Mandel's \\ problem \end{tabular} &
 \begin{tabular}{@{}l@{}}Manufactured \\ solution problem \end{tabular} &
 \begin{tabular}{@{}l@{}}Cantilever \\ problem \end{tabular}\\
    \midrule
    $\lambda$         & Lam\'e{}'s first parameter            & [Pa]    & $2.778 \times 10^{5}$ & $1.0$ & $1.429 \times 10^{5}$   \\
    $G$               & Shear modulus                         & [-]     & $4.167 \times 10^{5}$ & $1.0$ & $3.571 \times 10^{4}$   \\
    $b$               & Biot's coefficient                    & [-]     & 1.0 & 1.0 & 1.0 \\
    $S_{\varepsilon}$ & Constrained specific storage          & [Pa]    & 0.0 & 0.0 & 0.0 \\
    $\kappa$          & Isotropic permeability over viscosity & [m$^2$ $\cdot$ Pa$^{-1}$ $\cdot$ s$^{-1}$] & $1 \times 10^{-15}$ & $1.0$ & $1 \times 10^{-7}$ \\
    \midrule
    $a$               & Domain size in $x$-direction          & [m] & 1.0 & 1.0 & 1.0 \\
    $b$               & Domain size in $z$-direction          & [m] & 1.0 & 1.0 & 1.0 \\
    $F$               & Applied force magnitude               & [N $\cdot$ m$^{-1}$] & $2 \times 10^{2}$ & --- & $1.0$ \\
    \bottomrule
    \end{tabular}
\end{table}

% mainfile: main.tex
\subsection{Mandel's problem}
\label{sec:mandel}

%\begin{figure}[htbp]
%  \hfill
%  \begin{subfigure}{.3\linewidth}
%    \centerline{\includegraphics[width=
%      \linewidth]{./plotData/mandelTest/mandel_a}}
%    \caption{Physical domain.}
%    \label{fig:mandeldomain}
%  \end{subfigure}
%  \hfill
%  \begin{subfigure}{.3\linewidth}
%    \centerline{\includegraphics[width=
%      \linewidth]{./plotData/mandelTest/mandel_b}}
%    \caption{Computational domain.}
%    \label{fig:mandelcomputational}
%  \end{subfigure}
%  \hfill\null
%  \caption{Mandel's problem setup.}
%\end{figure}

%\begin{figure}[htbp]
%  \centering
%  \begin{subfigure}[b]{.49\linewidth}
%    \input{plotData/mandelTest/solution_fullycoupled_01000.tex}
%    \caption{$t = 0.1\cdot T_c$}
%  \end{subfigure}
%  \begin{subfigure}[b]{.49\linewidth}
%    \input{plotData/mandelTest/solution_fullycoupled_03000.tex}
%    \caption{$t = 0.3\cdot T_c$}
%  \end{subfigure}
%  \begin{subfigure}[b]{.49\linewidth}
%    \input{plotData/mandelTest/solution_fullycoupled_06000.tex}
%    \caption{$t = 0.6\cdot T_c$}
%  \end{subfigure}
%  \begin{subfigure}[b]{.49\linewidth}
%    \input{plotData/mandelTest/solution_fullycoupled_10000.tex}
%    \caption{$t = T_c$}
%  \end{subfigure}
%  \caption{Profiles of $\frac{p}{p^0}$ at $y = 0.5$}
%  \label{fig:mandelPressure}
%\end{figure}

\begin{figure}[htbp]
  \centering \small
  \begin{subfigure}{.475\linewidth}
    \begin{tikzpicture}
  \begin{axis}[%
    width=\linewidth,
    height=.85\linewidth,
    xmin=0,
    xmax=1,
    xlabel style={font=\color{white!15!black}},
    xlabel={$x / a$},
    ylabel={$p / p^0$},
    xlabel near ticks,
    xtick={0.0,0.2,...,1.0},
    ylabel near ticks,
    xtick={0.0,0.2,...,1.2},
    grid=major,
    ymin=0,
    ymax=1.2,
    axis background/.style={fill=white},
    legend style={font=\small},
    tick label style={font=\small},
    label style={font=\small},
    ]
    % Exact solution     
    \addplot [color=mycolor1]
    table[]{plotData/mandelTest/Cartesian/solution_fullycoupled_01000-2.tsv};
    \addlegendentry{$t / T_c = 0.1$}

    \addplot [color=mycolor2]
    table[]{plotData/mandelTest/Cartesian/solution_fullycoupled_03000-2.tsv};
    \addlegendentry{$t / T_c = 0.3$}
    
    \addplot [color=mycolor3]
    table[]{plotData/mandelTest/Cartesian/solution_fullycoupled_06000-2.tsv};
    \addlegendentry{$t / T_c = 0.6$}
    
    \addplot [color=mycolor4]
    table[]{plotData/mandelTest/Cartesian/solution_fullycoupled_10000-2.tsv};
    \addlegendentry{$t / T_c = 1.0$}
    
    % Numerical solution
    \addplot [only marks, mark=*, mark size =1.5pt, mark options={solid, fill=mycolor1}]
    table[]{plotData/mandelTest/Cartesian/solution_fullycoupled_01000-1.tsv};

    \addplot [only marks, mark=*, mark size =1.5pt, mark options={solid, fill=mycolor2}]
    table[]{plotData/mandelTest/Cartesian/solution_fullycoupled_03000-1.tsv};

    \addplot [only marks, mark=*, mark size =1.5pt, mark options={solid, fill=mycolor3}]
    table[]{plotData/mandelTest/Cartesian/solution_fullycoupled_06000-1.tsv};

    \addplot [only marks, mark=*, mark size =1.5pt, mark options={solid, fill=mycolor4}]
    table[]{plotData/mandelTest/Cartesian/solution_fullycoupled_10000-1.tsv};

  \end{axis}
\end{tikzpicture}
    \caption{\texttt{Cartesian}.}
    \label{fig:mandel_pressure_cartesian}
  \end{subfigure}
  \hfill
  \begin{subfigure}{.475\linewidth}
    \begin{tikzpicture}
  \begin{axis}[%
    width=\linewidth,
    height=.85\linewidth,
    xmin=0,
    xmax=1,
    xlabel style={font=\color{white!15!black}},
    xlabel={$x / a$},
    ylabel={$p / p^0$},
    xlabel near ticks,
    xtick={0.0,0.2,...,1.0},
    ylabel near ticks,
    xtick={0.0,0.2,...,1.2},
    grid=major,
    ymin=0,
    ymax=1.2,
    axis background/.style={fill=white},
    legend style={font=\small},
    tick label style={font=\small},
    label style={font=\small},
    ]
    % Exact solution     
    \addplot [color=mycolor1]
    table[]{plotData/mandelTest/Skewed/solution_fullycoupled_01000-2.tsv};
    \addlegendentry{$t / T_c = 0.1$}

    \addplot [color=mycolor2]
    table[]{plotData/mandelTest/Skewed/solution_fullycoupled_03000-2.tsv};
    \addlegendentry{$t / T_c = 0.3$}
    
    \addplot [color=mycolor3]
    table[]{plotData/mandelTest/Skewed/solution_fullycoupled_06000-2.tsv};
    \addlegendentry{$t / T_c = 0.6$}
    
    \addplot [color=mycolor4]
    table[]{plotData/mandelTest/Skewed/solution_fullycoupled_10000-2.tsv};
    \addlegendentry{$t / T_c = 1.0$}
    
    % Numerical solution
    \addplot [only marks, mark=*, mark size =1.5pt, mark options={solid, fill=mycolor1}]
    table[]{plotData/mandelTest/Skewed/solution_fullycoupled_01000-1.tsv};

    \addplot [only marks, mark=*, mark size =1.5pt, mark options={solid, fill=mycolor2}]
    table[]{plotData/mandelTest/Skewed/solution_fullycoupled_03000-1.tsv};

    \addplot [only marks, mark=*, mark size =1.5pt, mark options={solid, fill=mycolor3}]
    table[]{plotData/mandelTest/Skewed/solution_fullycoupled_06000-1.tsv};

    \addplot [only marks, mark=*, mark size =1.5pt, mark options={solid, fill=mycolor4}]
    table[]{plotData/mandelTest/Skewed/solution_fullycoupled_10000-1.tsv};

  \end{axis}
\end{tikzpicture}
    \caption{\texttt{Skewed}.}
    \label{fig:mandel_pressure_skewed}
  \end{subfigure}
  \hfill\null
  \begin{subfigure}{.475\linewidth}
    \begin{tikzpicture}
  \begin{axis}[%
    width=\linewidth,
    height=.85\linewidth,
    xmin=0,
    xmax=1,
    xlabel style={font=\color{white!15!black}},
    xlabel={$x / a$},
    ylabel={$p / p^0$},
    xlabel near ticks,
    xtick={0.0,0.2,...,1.0},
    ylabel near ticks,
    xtick={0.0,0.2,...,1.2},
    grid=major,
    ymin=0,
    ymax=1.2,
    axis background/.style={fill=white},
    legend style={font=\small},
    tick label style={font=\small},
    label style={font=\small},
    ]
    % Exact solution     
    \addplot [color=mycolor1]
    table[]{plotData/mandelTest/Hybrid/solution_fullycoupled_01000-2.tsv};
    \addlegendentry{$t / T_c = 0.1$}

    \addplot [color=mycolor2]
    table[]{plotData/mandelTest/Hybrid/solution_fullycoupled_03000-2.tsv};
    \addlegendentry{$t / T_c = 0.3$}
    
    \addplot [color=mycolor3]
    table[]{plotData/mandelTest/Hybrid/solution_fullycoupled_06000-2.tsv};
    \addlegendentry{$t / T_c = 0.6$}
    
    \addplot [color=mycolor4]
    table[]{plotData/mandelTest/Hybrid/solution_fullycoupled_10000-2.tsv};
    \addlegendentry{$t / T_c = 1.0$}
    
    % Numerical solution
    \addplot [only marks, mark=*, mark size =1.5pt, mark options={solid, fill=mycolor1}]
    table[]{plotData/mandelTest/Hybrid/solution_fullycoupled_01000-1.tsv};

    \addplot [only marks, mark=*, mark size =1.5pt, mark options={solid, fill=mycolor2}]
    table[]{plotData/mandelTest/Hybrid/solution_fullycoupled_03000-1.tsv};

    \addplot [only marks, mark=*, mark size =1.5pt, mark options={solid, fill=mycolor3}]
    table[]{plotData/mandelTest/Hybrid/solution_fullycoupled_06000-1.tsv};

    \addplot [only marks, mark=*, mark size =1.5pt, mark options={solid, fill=mycolor4}]
    table[]{plotData/mandelTest/Hybrid/solution_fullycoupled_10000-1.tsv};

  \end{axis}
\end{tikzpicture}
    \caption{\texttt{Hybrid}.}
    \label{fig:mandel_pressure_hybrid}
  \end{subfigure}
  \hfill
  \begin{subfigure}{.475\linewidth}
    \begin{tikzpicture}
  \begin{axis}[%
    width=\linewidth,
    height=.85\linewidth,
    xmin=0,
    xmax=1,
    xlabel style={font=\color{white!15!black}},
    xlabel={$x / a$},
    ylabel={$p / p^0$},
    xlabel near ticks,
    xtick={0.0,0.2,...,1.0},
    ylabel near ticks,
    xtick={0.0,0.2,...,1.2},
    grid=major,
    ymin=0,
    ymax=1.2,
    axis background/.style={fill=white},
    legend style={font=\small},
    tick label style={font=\small},
    label style={font=\small},
    ]
    % Exact solution     
    \addplot [color=mycolor1]
    table[]{plotData/mandelTest/Polymesher1/solution_fullycoupled_01000-2.tsv};
    \addlegendentry{$t / T_c = 0.1$}

    \addplot [color=mycolor2]
    table[]{plotData/mandelTest/Polymesher1/solution_fullycoupled_03000-2.tsv};
    \addlegendentry{$t / T_c = 0.3$}
    
    \addplot [color=mycolor3]
    table[]{plotData/mandelTest/Polymesher1/solution_fullycoupled_06000-2.tsv};
    \addlegendentry{$t / T_c = 0.6$}
    
    \addplot [color=mycolor4]
    table[]{plotData/mandelTest/Polymesher1/solution_fullycoupled_10000-2.tsv};
    \addlegendentry{$t / T_c = 1.0$}
    
    % Numerical solution
    \addplot [only marks, mark=*, mark size =1.5pt, mark options={solid, fill=mycolor1}]
    table[]{plotData/mandelTest/Polymesher1/solution_fullycoupled_01000-1.tsv};

    \addplot [only marks, mark=*, mark size =1.5pt, mark options={solid, fill=mycolor2}]
    table[]{plotData/mandelTest/Polymesher1/solution_fullycoupled_03000-1.tsv};

    \addplot [only marks, mark=*, mark size =1.5pt, mark options={solid, fill=mycolor3}]
    table[]{plotData/mandelTest/Polymesher1/solution_fullycoupled_06000-1.tsv};

    \addplot [only marks, mark=*, mark size =1.5pt, mark options={solid, fill=mycolor4}]
    table[]{plotData/mandelTest/Polymesher1/solution_fullycoupled_10000-1.tsv};

  \end{axis}
\end{tikzpicture}
    \caption{\texttt{Polymesher1}.}
    \label{fig:mandel_pressure_poly1}
  \end{subfigure}
  \hfill\null
  \begin{subfigure}{.475\linewidth}
    \begin{tikzpicture}
  \begin{axis}[%
    width=\linewidth,
    height=.85\linewidth,
    xmin=0,
    xmax=1,
    xlabel style={font=\color{white!15!black}},
    xlabel={$x / a$},
    ylabel={$p / p^0$},
    xlabel near ticks,
    xtick={0.0,0.2,...,1.0},
    ylabel near ticks,
    xtick={0.0,0.2,...,1.2},
    grid=major,
    ymin=0,
    ymax=1.2,
    axis background/.style={fill=white},
    legend style={font=\small},
    tick label style={font=\small},
    label style={font=\small},
    ]
    % Exact solution     
    \addplot [color=mycolor1]
    table[]{plotData/mandelTest/Polymesher20/solution_fullycoupled_01000-2.tsv};
    \addlegendentry{$t / T_c = 0.1$}

    \addplot [color=mycolor2]
    table[]{plotData/mandelTest/Polymesher20/solution_fullycoupled_03000-2.tsv};
    \addlegendentry{$t / T_c = 0.3$}
    
    \addplot [color=mycolor3]
    table[]{plotData/mandelTest/Polymesher20/solution_fullycoupled_06000-2.tsv};
    \addlegendentry{$t / T_c = 0.6$}
    
    \addplot [color=mycolor4]
    table[]{plotData/mandelTest/Polymesher20/solution_fullycoupled_10000-2.tsv};
    \addlegendentry{$t / T_c = 1.0$}
    
    % Numerical solution
    \addplot [only marks, mark=*, mark size =1.5pt, mark options={solid, fill=mycolor1}]
    table[]{plotData/mandelTest/Polymesher20/solution_fullycoupled_01000-1.tsv};

    \addplot [only marks, mark=*, mark size =1.5pt, mark options={solid, fill=mycolor2}]
    table[]{plotData/mandelTest/Polymesher20/solution_fullycoupled_03000-1.tsv};

    \addplot [only marks, mark=*, mark size =1.5pt, mark options={solid, fill=mycolor3}]
    table[]{plotData/mandelTest/Polymesher20/solution_fullycoupled_06000-1.tsv};

    \addplot [only marks, mark=*, mark size =1.5pt, mark options={solid, fill=mycolor4}]
    table[]{plotData/mandelTest/Polymesher20/solution_fullycoupled_10000-1.tsv};

  \end{axis}
\end{tikzpicture}
    \caption{\texttt{Polymesher20}.}
    \label{fig:mandel_pressure_poly20}
  \end{subfigure}

  \caption{Mandel's problem: normalized analytical (continuous
    profiles) and MFD-VEM (circle markers) pressure solutions along a
    cross section at $y = 0.5$ for a 400-cell mesh from the
    \texttt{Cartesian} (a), \texttt{Skewed} (b),
    \texttt{Hybrid} (c), \texttt{Polymesher1} (d) and \texttt{Polymesher20}
    (e) families, respectively. }
  \label{fig:mandelPressure}
\end{figure}
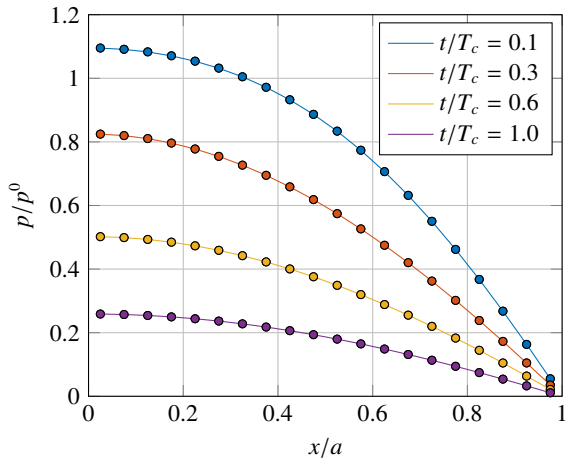
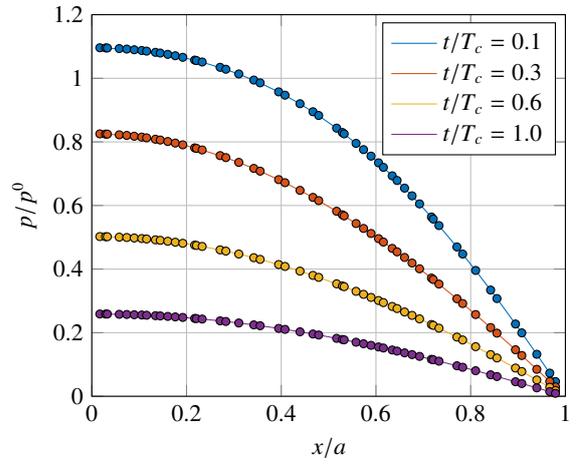
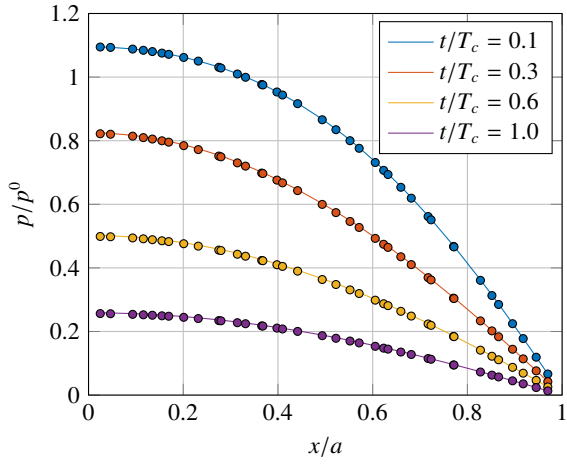
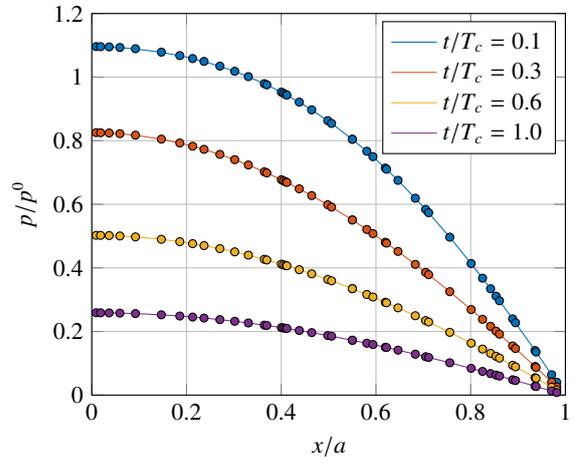
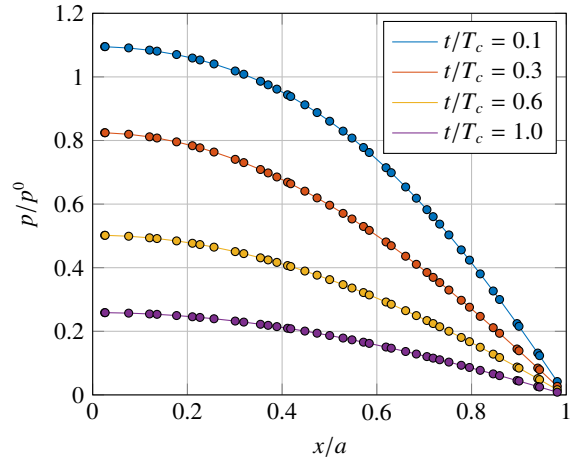

Mandel's problem \cite{Mandel1953} consists of a sample of saturated, isotropic, poroelastic material that is loaded under plane-strain conditions by a constant compressive force of magnitude $2F$ applied at time $t = 0$ on rigid, frictionless, impermeable plates.
The sample has dimensions $2a \times 2b$ (Fig.~\ref{fig:mandel_sketch}).
The left and right sides ($x = \pm a$) are stress free, drained and kept at constant ambient pressure.
The exact analytical solutions of this problem are known 
(see \cite[Appendix A.2]{castelletto2015}).

Given the symmetry of the problem, we solve it on a quarter of the domain, represented by the square $\Omega = (0,a)\times (0,b)$ (Fig.~\ref{fig:mandel_sketch}).
The simulation setup is the same as that used in \cite{castelletto2015}.
We use the unit square as the domain, i.e. $a = b = 1$ m.
The rigid plate constraint is accounted for by prescribing the vertical displacement at the loaded boundaries using the closed-form solution of the problem.
%
%We impose stress-free boundary conditions for the mechanics everywhere except for $u_1(0,y) = 0$, $u_2(x,0) = 0$, and setting the exact solution for $u_2(x,1)$.
%
%For the flow equation, we impose the exact pressure on the right boundary of
%the domain and no flow conditions on the other boundaries.

In Fig.~\ref{fig:mandelPressure} we compare analytical and MFD-VEM solutions for the pressure along a cross-section at $y = 0.5$ at different times.
The results are relative to 400-cell meshes from the mesh families in Fig.~\ref{fig:mesh_illustration}.
A constant time step size $\Delta t = 10^{-4}\cdot T_c$ is used in both cases, with $T_c = {a^2}/({\kappa(\lambda + 2G)})$ the characteristic time of the consolidation process.
We can see that the computed values match the expected profile very well, with relative error values of the order of $10^{-4}$.
Note that the parameters used are representative of the limit case of
incompressible solid ($b$ = 1.0) and fluid ($S_{\varepsilon}$ = 0) constituents, which maximizes the hydromechanical coupling, i.e. the so-called Mandel-Cryer effect.
These results confirm the robustness of the method with respect to badly shaped elements.
Because of the lack of regularity of the pressure field for Mandel's problem \cite{PhiWhe07,Yi14}, we do not perform a mesh refinement study to evaluate numerically the convergence rate of quantities of interest here. 
This is addressed in the next section.

%%% Local Variables:
%%% mode: latex
%%% TeX-master: "main"
%%% End:

% mainfile: main.tex
\subsection{Test with exact solution}
\label{sec:convTest}

\begin{figure}
  \small
   
  \hfill
  \begin{subfigure}{.32\linewidth} 
    \begin{tikzpicture}  
   \pgfplotsset{every axis legend/.append style={
     at={(0.5,1.1)},
     anchor=south}}

  \begin{loglogaxis}[
    width= \textwidth,
    height= \textwidth,
    xmin=1e-2,xmax=1e0,
    ymin=1e-6,ymax=1e-2,
    xtick={1e-2,1e-1,1e0},
    ytick={1e-6,1e-5,1e-4,1e-3,1e-2,1e-2},
    xlabel near ticks,    
    ylabel near ticks,
    x dir=reverse,
    xlabel={$h$},  
    ylabel={$e_p$},  
    grid=major,
    legend cell align={left},
    legend style={anchor=south},
    ]
    \logLogSlopeTriangle{0.6}{0.25}{0.2}{1}{black};
    
    \addplot [color=mycolor1,mark=*,mark size =1.5pt]
    table[]{./plotData/convergenceTest/Cartesian/totalErrors-1.tsv};
    \addlegendentry{{\footnotesize{\texttt{Cartesian}}} ($m$ = 1.05)}
    
    \addplot [color=mycolor2,mark=*,mark size =1.5pt]
    table[]{./plotData/convergenceTest/Hybrid/totalErrors-1.tsv};
    \addlegendentry{{\footnotesize{\texttt{Hybrid}}} ($m$ = 1.06)}

    \addplot [color=mycolor3,mark=*,mark size =1.5pt]
    table[]{./plotData/convergenceTest/Skewed/totalErrors-1.tsv};
    \addlegendentry{{\footnotesize{\texttt{Skewed}}} ($m$ = 1.07)}

    \addplot [color=mycolor4,mark=*,mark size =1.5pt]
    table[]{./plotData/convergenceTest/Polymesher1/totalErrors-1.tsv};
    \addlegendentry{{\footnotesize{\texttt{Polymesher1}}} ($m$ = 1.36)}
    
    \addplot [color=mycolor5,mark=*,mark size =1.5pt]
    table[]{./plotData/convergenceTest/Polymesher20/totalErrors-1.tsv};
    \addlegendentry{{\footnotesize{\texttt{PolyMesher20}}} ($m$ = 0.91)} 
  \end{loglogaxis}
\end{tikzpicture}
    \caption{}
  \end{subfigure}
  \hfill
  \begin{subfigure}{.32\linewidth}
	  \begin{tikzpicture}  
      \pgfplotsset{every axis legend/.append style={
       at={(0.5,1.1)},
       anchor=south}}

  \begin{loglogaxis}[
    width= \textwidth,
    height= \textwidth,
    xmin=1e-2,xmax=1e0,
    ymin=1e-7,ymax=1e-3,
    xtick={1e-2,1e-1,1e0},
    ytick={1e-7,1e-6,1e-5,1e-4,1e-3},
    xlabel near ticks,    
    ylabel near ticks,
    x dir=reverse,    
    xlabel={$h$},  
    ylabel={$e_{u}$},  
    grid=major,
    legend cell align={left},
    legend style={anchor=south}
    ]
     \logLogSlopeTriangle{0.6}{0.25}{0.2}{1}{black};
    
    \addplot [color=mycolor1,mark=*,mark size =1.5pt]
    table[]{./plotData/convergenceTest/Cartesian/totalErrors-3.tsv};
    \addlegendentry{{\footnotesize{\texttt{Cartesian}}} (m = 0.98)}
   
    \addplot [color=mycolor2,mark=*,mark size =1.5pt]
    table[]{./plotData/convergenceTest/Hybrid/totalErrors-3.tsv};
    \addlegendentry{{\footnotesize{\texttt{Hybrid}}} (m = 1.16)}

    \addplot [color=mycolor3,mark=*,mark size =1.5pt]
    table[]{./plotData/convergenceTest/Skewed/totalErrors-3.tsv};
    \addlegendentry{{\footnotesize{\texttt{Skewed}}} (m = 1.49)}

    \addplot [color=mycolor4,mark=*,mark size =1.5pt]
    table[]{./plotData/convergenceTest/Polymesher1/totalErrors-3.tsv};
    \addlegendentry{{\footnotesize{\texttt{Polymesher1}}} (m = 1.34)}
    
    \addplot [color=mycolor5,mark=*,mark size =1.5pt]
    table[]{./plotData/convergenceTest/Polymesher20/totalErrors-3.tsv};
    \addlegendentry{{\footnotesize{\texttt{PolyMesher20}}} (m = 1.01)} 
  \end{loglogaxis}
\end{tikzpicture}
    \caption{}
  \end{subfigure}
%  \hfill\null
%  
\hfill
  \begin{subfigure}{.32\linewidth}
      \begin{tikzpicture}  
      \pgfplotsset{every axis legend/.append style={
       at={(0.5,1.1)},
       anchor=south}}

  \begin{loglogaxis}[
    width= \textwidth,
    height= \textwidth,
    xmin=1e-2,xmax=1e0,
    ymin=1e-6,ymax=1e-2,
    xtick={1e-2,1e-1,1e0},
    ytick={1e-6,1e-5,1e-4,1e-3,1e-2,1e-2},
    xlabel near ticks,    
    ylabel near ticks,
    x dir=reverse,
%    legend style={font=\small},
%    tick label style={font=\small},
%    label style={font=\small},     
    xlabel={$h$},  
    ylabel={$e_{\sigma^{\prime}}$},  
    grid=major,
    %legend pos = north east,
    legend cell align={left},
    legend style={anchor=south},
    ]
    \logLogSlopeTriangle{0.6}{0.25}{0.2}{1}{black};
    
    \addplot [color=mycolor1,mark=*,mark size =1.5pt]%, mark options={solid, fill=mycolor1}]
    table[]{./plotData/convergenceTest/Cartesian/totalErrors-2.tsv};
    \addlegendentry{{\footnotesize{\texttt{Cartesian}}} (m = 0.99)}
    
    \addplot [color=mycolor2,mark=*,mark size =1.5pt]%, mark options={solid, fill=mycolor2}]
    table[]{./plotData/convergenceTest/Hybrid/totalErrors-2.tsv};
    \addlegendentry{{\footnotesize{\texttt{Hybrid}}} (m = 1.11)}

    \addplot [color=mycolor3,mark=*,mark size =1.5pt]%, mark options={solid, fill=mycolor2}]
    table[]{./plotData/convergenceTest/Skewed/totalErrors-2.tsv};
    \addlegendentry{{\footnotesize{\texttt{Skewed}}} (m = 1.89)}

    \addplot [color=mycolor4,mark=*,mark size =1.5pt]%, mark options={solid, fill=mycolor3}]
    table[]{./plotData/convergenceTest/Polymesher1/totalErrors-2.tsv};
    \addlegendentry{{\footnotesize{\texttt{Polymesher1}}} (m = 1.15)}
    
    \addplot [color=mycolor5,mark=*,mark size =1.5pt]%, mark options={solid, fill=mycolor4}]
    table[]{./plotData/convergenceTest/Polymesher20/totalErrors-2.tsv};
    \addlegendentry{{\footnotesize{\texttt{PolyMesher20}}} (m = 0.93)} 
  \end{loglogaxis}
\end{tikzpicture}
    \caption{}
  \end{subfigure}
  \hfill\null  
  
  \caption{Test with manufactured exact solution: total error behavior for the pressure (a), displacement (b), and effective stress (c). For each profile the convergence rates are provided in parenthesis in the legend.}
  \label{fig:totalErrors}
\end{figure}
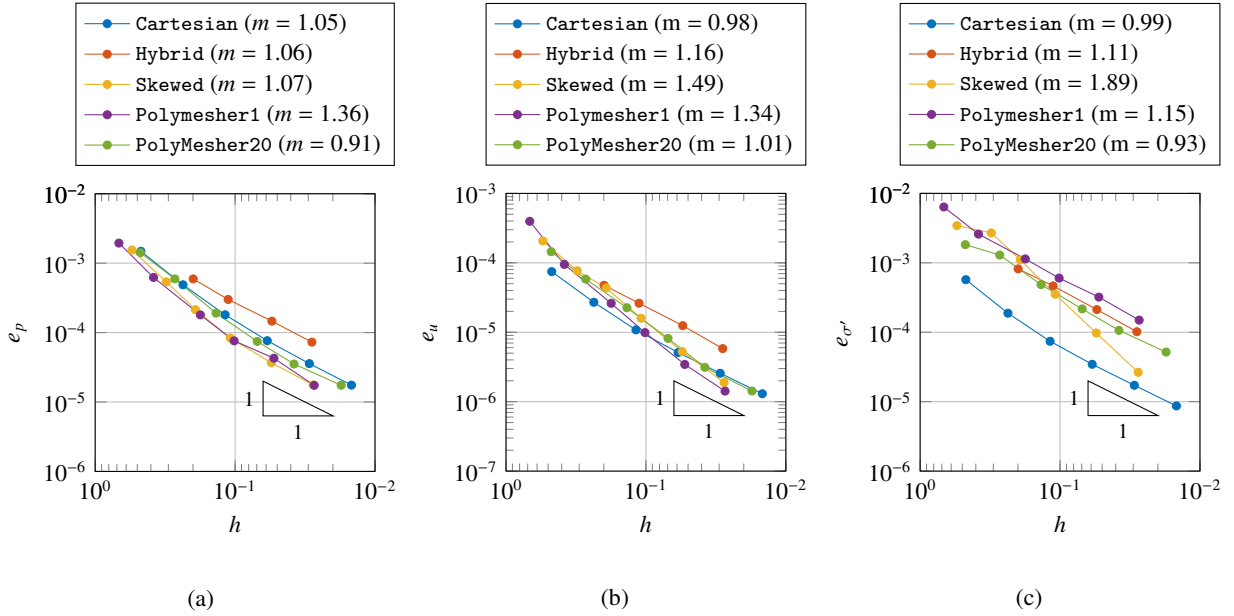

\begin{figure}
  \centering
  \begin{subfigure}{.45\linewidth}
    \centering
    \begin{tikzpicture}
  \begin{loglogaxis}[%
    width= 0.7\textwidth,
    height= 0.75\textwidth,
    xmin=1e-3,xmax=2e-1,
    ymin=1e-5,ymax=5e-3,
    xtick={1e-4,1e-3,1e-2,1e-1,1e0},
    ytick={1e-5,1e-4,1e-3,1e-2,1e-1},
    x dir=reverse,
    legend style={font=\small},
    tick label style={font=\small},
    label style={font=\small},
    xlabel={$\tau$},
    grid=major,
    legend pos = outer north east,
    legend style = {draw},
    every axis plot/.append style={mark=*},
    legend columns=1,
    ]
    \logLogSlopeTriangle{0.3}{0.2}{0.15}{1}{black};

    \addplot [color=mycolor1]
    table[]{plotData/convergenceTest/Cartesian/pressureDifferentDiam-1.tsv};
    \addlegendentry{h = 0.4714% , m = 0.01
    }

    \addplot [color=mycolor2]
    table[]{plotData/convergenceTest/Cartesian/pressureDifferentDiam-2.tsv};
    \addlegendentry{h = 0.2357% , m = 0.01
    }

    \addplot [color=mycolor3]
    table[]{plotData/convergenceTest/Cartesian/pressureDifferentDiam-3.tsv};
    \addlegendentry{h = 0.11785% , m = 0.05
    }

    \addplot [color=mycolor4]
    table[]{plotData/convergenceTest/Cartesian/pressureDifferentDiam-4.tsv};
    \addlegendentry{h = 0.058926% , m = 0.38
    }

    \addplot [color=mycolor5]
    table[]{plotData/convergenceTest/Cartesian/pressureDifferentDiam-5.tsv};
    \addlegendentry{h = 0.029463% , m = 0.88
    }

    \addplot [color=mycolor6]
    table[]{plotData/convergenceTest/Cartesian/pressureDifferentDiam-6.tsv};
    \addlegendentry{h = 0.014731% , m = 0.99
    }

    % \addplot [color=mycolor7, dashed]
    % table[]{plotData/convergenceTest/Cartesian/pressureDifferentDiam-7.tsv};
    % \addlegendentry{ref m=1}

  \end{loglogaxis}
\end{tikzpicture}%
    \caption{\texttt{Cartesian}.}
  \end{subfigure}
  \begin{subfigure}{.45\linewidth}
    \centering
    \begin{tikzpicture}

  \begin{loglogaxis}[%
    width= 0.7\textwidth,
    height= 0.75\textwidth,
    xmin=1e-3,xmax=2e-1,
    ymin=1e-5,ymax=5e-3,
    xtick={1e-4,1e-3,1e-2,1e-1,1e0},
    ytick={1e-5,1e-4,1e-3,1e-2,1e-1},
    x dir=reverse,
    legend style={font=\small},
    tick label style={font=\small},
    label style={font=\small},
    xlabel={$\tau$},
    grid=major,
    legend pos = outer north east,
    legend style = {draw},
    every axis plot/.append style={mark=*},
    legend columns=1,
    ]

    \logLogSlopeTriangle{0.3}{0.2}{0.15}{1}{black};

    \addplot [color=mycolor1]
    table[]{plotData/convergenceTest/Skewed/pressureDifferentDiam-1.tsv};
    \addlegendentry{h = 0.54565% , m = 0.01
    }

    \addplot [color=mycolor2]
    table[]{plotData/convergenceTest/Skewed/pressureDifferentDiam-2.tsv};
    \addlegendentry{h = 0.30995% , m = 0.01
    }

    \addplot [color=mycolor3]
    table[]{plotData/convergenceTest/Skewed/pressureDifferentDiam-3.tsv};
    \addlegendentry{h = 0.1921% , m = 0.03
    }

    \addplot [color=mycolor4]
    table[]{plotData/convergenceTest/Skewed/pressureDifferentDiam-4.tsv};
    \addlegendentry{h = 0.10842% , m = 0.23
    }

    \addplot [color=mycolor5]
    table[]{plotData/convergenceTest/Skewed/pressureDifferentDiam-5.tsv};
    \addlegendentry{h = 0.055085% , m = 0.79
    }

    \addplot [color=mycolor6]
    table[]{plotData/convergenceTest/Skewed/pressureDifferentDiam-6.tsv};
    \addlegendentry{h = 0.027653% , m = 0.98
    }

    % \addplot [color=mycolor7, dashed]
    % table[]{plotData/convergenceTest/Skewed/pressureDifferentDiam-7.tsv};
    % \addlegendentry{ref m=1}

  \end{loglogaxis}
\end{tikzpicture}%    
    \caption{\texttt{Skewed}.}
  \end{subfigure}
  \begin{subfigure}{.45\linewidth}
    \centering
     \begin{tikzpicture}
  \begin{loglogaxis}[%
    width= 0.7\textwidth,
    height= 0.75\textwidth,
    xmin=1e-2,xmax=2e-1,
    ymin=5e-5,ymax=1e-3,
    xtick={1e-4,1e-3,1e-2,1e-1,1e0},
    ytick={1e-5,1e-4,1e-3,1e-2,1e-1},
    x dir=reverse,
    legend style={font=\small},
    tick label style={font=\small},
    label style={font=\small},
    xlabel={$\tau$},
    grid=major,
    legend pos = outer north east,
    legend style = {draw},
    every axis plot/.append style={mark=*},
    legend columns=1,
    ]

    \logLogSlopeTriangle{0.5}{0.2}{0.15}{1}{black};

    \addplot [color=mycolor1]
    table[]{plotData/convergenceTest/Hybrid/pressureDifferentDiam-1.tsv};
    \addlegendentry{h = 0.19939% , m = 0.16
    }

    \addplot [color=mycolor2]
    table[]{plotData/convergenceTest/Hybrid/pressureDifferentDiam-2.tsv};
    \addlegendentry{h = 0.11207% , m = 0.65
    }

    \addplot [color=mycolor3]
    table[]{plotData/convergenceTest/Hybrid/pressureDifferentDiam-3.tsv};
    \addlegendentry{h = 0.054572% , m = 0.96
    }

    \addplot [color=mycolor4]
    table[]{plotData/convergenceTest/Hybrid/pressureDifferentDiam-4.tsv};
    \addlegendentry{h = 0.028265% , m = 1.00
    }

    % \addplot [color=mycolor5, dashed]
    % table[]{plotData/convergenceTest/Hybrid/pressureDifferentDiam-5.tsv};
    % \addlegendentry{ref m=1}

  \end{loglogaxis}
\end{tikzpicture}%
    \caption{\texttt{Hybrid}.}
  \end{subfigure}
  \begin{subfigure}{.45\linewidth}
    \centering
    \begin{tikzpicture}
  \begin{loglogaxis}[%
    width= 0.7\textwidth,
    height= 0.75\textwidth,
    xmin=1e-3,xmax=2e-1,
    ymin=1e-5,ymax=5e-3,
    xtick={1e-4,1e-3,1e-2,1e-1,1e0},
    ytick={1e-5,1e-4,1e-3,1e-2,1e-1},
    x dir=reverse,
    legend style={font=\small},
    tick label style={font=\small},
    label style={font=\small},
    xlabel={$\tau$},
    grid=major,
    legend pos = outer north east,
    legend style = {draw},
    every axis plot/.append style={mark=*},
    legend columns=1,
    ]
    \logLogSlopeTriangle{0.3}{0.2}{0.15}{1}{black};
    
    \addplot [color=mycolor1]
    table[]{plotData/convergenceTest/Polymesher1/pressureDifferentDiam-1.tsv};
    \addlegendentry{h = 0.67702% , m = 0.01
    }

    \addplot [color=mycolor2]
    table[]{plotData/convergenceTest/Polymesher1/pressureDifferentDiam-2.tsv};
    \addlegendentry{h = 0.38297% , m = 0.01
    }

    \addplot [color=mycolor3]
    table[]{plotData/convergenceTest/Polymesher1/pressureDifferentDiam-3.tsv};
    \addlegendentry{h = 0.17723% , m = 0.05
    }

    \addplot [color=mycolor4]
    table[]{plotData/convergenceTest/Polymesher1/pressureDifferentDiam-4.tsv};
    \addlegendentry{h = 0.10159% , m = 0.40
    }

    \addplot [color=mycolor5]
    table[]{plotData/convergenceTest/Polymesher1/pressureDifferentDiam-5.tsv};
    \addlegendentry{h = 0.052792% , m = 0.50
    }

    \addplot [color=mycolor6]
    table[]{plotData/convergenceTest/Polymesher1/pressureDifferentDiam-6.tsv};
    \addlegendentry{h = 0.027202% , m = 1.00
    }

    % \addplot [color=mycolor7, dashed]
    % table[]{plotData/convergenceTest/Polymesher1/pressureDifferentDiam-7.tsv};
    % \addlegendentry{ref m=1}

  \end{loglogaxis}
\end{tikzpicture}%
    \caption{\texttt{Polymesher1}.}
  \end{subfigure}
  \begin{subfigure}{.45\linewidth}
    \centering
    \begin{tikzpicture}
  \begin{loglogaxis}[%
    width= 0.7\textwidth,
    height= 0.75\textwidth,
    xmin=1e-3,xmax=2e-1,
    ymin=1e-5,ymax=5e-3,
    xtick={1e-4,1e-3,1e-2,1e-1,1e0},
    ytick={1e-5,1e-4,1e-3,1e-2,1e-1},
    x dir=reverse,
    legend style={font=\small},
    tick label style={font=\small},
    label style={font=\small},
    xlabel={$\tau$},
    grid=major,
    legend pos = outer north east,
    legend style = {draw},
    every axis plot/.append style={mark=*},
    legend columns=1,
    ]

    \logLogSlopeTriangle{0.3}{0.2}{0.15}{1}{black};

    \addplot [color=mycolor1]
    table[]{plotData/convergenceTest/Polymesher20/pressureDifferentDiam-1.tsv};
    \addlegendentry{h = 0.47563% , m = 0.01
    }

    \addplot [color=mycolor2]
    table[]{plotData/convergenceTest/Polymesher20/pressureDifferentDiam-2.tsv};
    \addlegendentry{h = 0.26968% , m = 0.01
    }

    \addplot [color=mycolor3]
    table[]{plotData/convergenceTest/Polymesher20/pressureDifferentDiam-3.tsv};
    \addlegendentry{h = 0.13683% , m = 0.04
    }

    \addplot [color=mycolor4]
    table[]{plotData/convergenceTest/Polymesher20/pressureDifferentDiam-4.tsv};
    \addlegendentry{h = 0.069522% , m = 0.45
    }

    \addplot [color=mycolor5]
    table[]{plotData/convergenceTest/Polymesher20/pressureDifferentDiam-5.tsv};
    \addlegendentry{h = 0.037938% , m = 0.93
    }

    \addplot [color=mycolor6]
    table[]{plotData/convergenceTest/Polymesher20/pressureDifferentDiam-6.tsv};
    \addlegendentry{h = 0.017468% , m = 0.99
    }

    % \addplot [color=mycolor7, dashed]
    % table[]{plotData/convergenceTest/Polymesher20/pressureDifferentDiam-7.tsv};
    % \addlegendentry{ref m=1}

  \end{loglogaxis}
\end{tikzpicture}%    
    \caption{\texttt{Polymesher20}.}
  \end{subfigure}
  \caption{Test with manufactured exact solution: pressure error behavior through space refinements.}
  \label{fig:pressureDifferentDiams}
\end{figure}
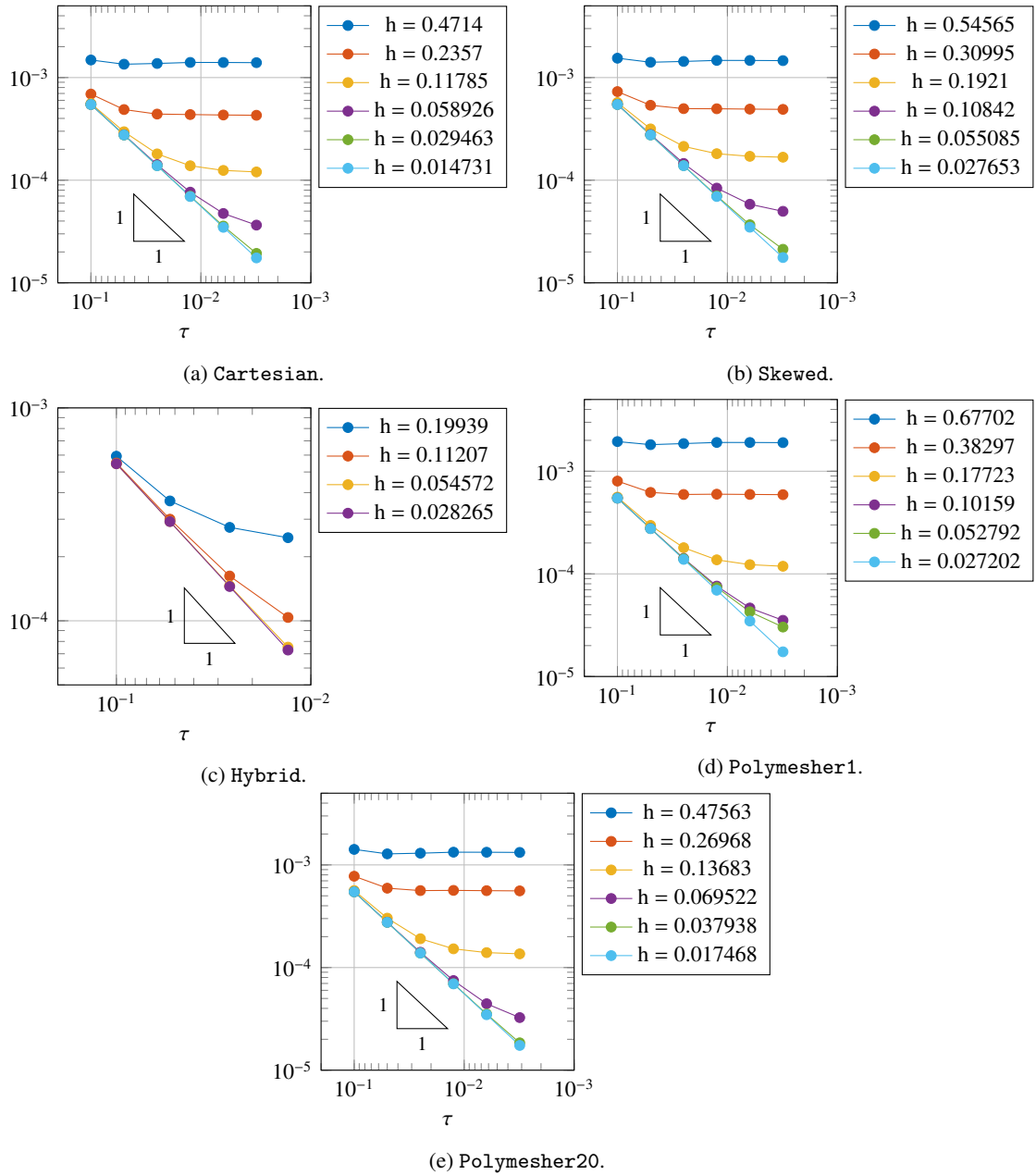

\begin{figure}
  \centering
  \begin{subfigure}{.49\linewidth}
    \centering
    \begin{tikzpicture}

  \begin{loglogaxis}[%
    width= 0.7\textwidth,
    height= 0.75\textwidth,
    xmin=1e-3,xmax=2e-1,
    ymin=1e-6,ymax=1e-4,
    xtick={1e-4,1e-3,1e-2,1e-1,1e0},
    ytick={1e-6,1e-5,1e-4,1e-3,1e-2,1e-1},
    x dir=reverse,
    legend style={font=\small},
    tick label style={font=\small},
    label style={font=\small},
    xlabel={$\tau$},
    grid=major,
    legend pos = outer north east,
    legend style = {draw},
    every axis plot/.append style={mark=*},
    legend columns=1,
    ]
    \logLogSlopeTriangle{0.3}{0.2}{0.15}{1}{black};

    \addplot [color=mycolor1]
    table[]{plotData/convergenceTest/Cartesian/displDifferentDiam-1.tsv};
    \addlegendentry{h = 0.4714% , m = 0.00
    }

    \addplot [color=mycolor2]
    table[]{plotData/convergenceTest/Cartesian/displDifferentDiam-2.tsv};
    \addlegendentry{h = 0.2357% , m = -0.01
    }

    \addplot [color=mycolor3]
    table[]{plotData/convergenceTest/Cartesian/displDifferentDiam-3.tsv};
    \addlegendentry{h = 0.11785% , m = 0.01
    }

    \addplot [color=mycolor4]
    table[]{plotData/convergenceTest/Cartesian/displDifferentDiam-4.tsv};
    \addlegendentry{h = 0.058926% , m = 0.52
    }

    \addplot [color=mycolor5]
    table[]{plotData/convergenceTest/Cartesian/displDifferentDiam-5.tsv};
    \addlegendentry{h = 0.029463% , m = 1.00
    }

    \addplot [color=mycolor6]
    table[]{plotData/convergenceTest/Cartesian/displDifferentDiam-6.tsv};
    \addlegendentry{h = 0.014731% , m = 1.01
    }

    % \addplot [color=mycolor7, dashed]
    % table[]{plotData/convergenceTest/Cartesian/displDifferentDiam-7.tsv};
    % \addlegendentry{ref m=1}

  \end{loglogaxis}
\end{tikzpicture}%    
    \caption{\texttt{Cartesian}.}
  \end{subfigure}
  \begin{subfigure}{.49\linewidth}
    \centering
    \begin{tikzpicture}

  \begin{loglogaxis}[%
    width= 0.7\textwidth,
    height= 0.75\textwidth,
    xmin=1e-3,xmax=2e-1,
    ymin=1e-6,ymax=1e-4,
    xtick={1e-4,1e-3,1e-2,1e-1,1e0},
    ytick={1e-6,1e-5,1e-4,1e-3,1e-2,1e-1},
    x dir=reverse,
    legend style={font=\small},
    tick label style={font=\small},
    label style={font=\small},
    xlabel={$\tau$},
    grid=major,
    legend pos = outer north east,
    legend style = {draw},
    every axis plot/.append style={mark=*},
    legend columns=1,
    ]
    \logLogSlopeTriangle{0.3}{0.2}{0.15}{1}{black};

    \addplot [color=mycolor1]
    table[]{plotData/convergenceTest/Skewed/displDifferentDiam-1.tsv};
    \addlegendentry{h = 0.54565% , m = -0.00
    }

    \addplot [color=mycolor2]
    table[]{plotData/convergenceTest/Skewed/displDifferentDiam-2.tsv};
    \addlegendentry{h = 0.30995% , m = 0.01
    }

    \addplot [color=mycolor3]
    table[]{plotData/convergenceTest/Skewed/displDifferentDiam-3.tsv};
    \addlegendentry{h = 0.1921% , m = 0.02
    }

    \addplot [color=mycolor4]
    table[]{plotData/convergenceTest/Skewed/displDifferentDiam-4.tsv};
    \addlegendentry{h = 0.10842% , m = 0.06
    }

    \addplot [color=mycolor5]
    table[]{plotData/convergenceTest/Skewed/displDifferentDiam-5.tsv};
    \addlegendentry{h = 0.055085% , m = 0.26
    }

    \addplot [color=mycolor6]
    table[]{plotData/convergenceTest/Skewed/displDifferentDiam-6.tsv};
    \addlegendentry{h = 0.027653% , m = 0.72
    }

    % \addplot [color=mycolor7, dashed]
    % table[]{plotData/convergenceTest/Skewed/displDifferentDiam-7.tsv};
    % \addlegendentry{ref m=1}

  \end{loglogaxis}
\end{tikzpicture}%
    \caption{\texttt{Skewed}.}
  \end{subfigure}
  \begin{subfigure}{.49\linewidth}
    \centering
    \begin{tikzpicture}

  \begin{loglogaxis}[%
    width= 0.7\textwidth,
    height= 0.75\textwidth,
    xmin=1e-2,xmax=2e-1,
    ymin=5e-6,ymax=1e-4,
    xtick={1e-4,1e-3,1e-2,1e-1,1e0},
    ytick={1e-6,1e-5,1e-4,1e-3,1e-2,1e-1},
    x dir=reverse,
    legend style={font=\small},
    tick label style={font=\small},
    label style={font=\small},
    xlabel={$\tau$},
    grid=major,
    legend pos = outer north east,
    legend style = {draw},
    every axis plot/.append style={mark=*},
    legend columns=1,
    ]
    \logLogSlopeTriangle{0.3}{0.2}{0.15}{1}{black};

    \addplot [color=mycolor1]
    table[]{plotData/convergenceTest/Hybrid/displDifferentDiam-1.tsv};
    \addlegendentry{h = 0.19939% , m = 0.08
    }

    \addplot [color=mycolor2]
    table[]{plotData/convergenceTest/Hybrid/displDifferentDiam-2.tsv};
    \addlegendentry{h = 0.11207% , m = 0.42
    }

    \addplot [color=mycolor3]
    table[]{plotData/convergenceTest/Hybrid/displDifferentDiam-3.tsv};
    \addlegendentry{h = 0.054572% , m = 0.81
    }

    \addplot [color=mycolor4]
    table[]{plotData/convergenceTest/Hybrid/displDifferentDiam-4.tsv};
    \addlegendentry{h = 0.028265% , m = 0.96
    }

    % \addplot [color=mycolor5, dashed]
    % table[]{plotData/convergenceTest/Hybrid/displDifferentDiam-5.tsv};
    % \addlegendentry{ref m=1}

  \end{loglogaxis}
\end{tikzpicture}%
    \caption{\texttt{Hybrid}.}
  \end{subfigure}
  \begin{subfigure}{.49\linewidth}
    \centering
    \begin{tikzpicture}

  \begin{loglogaxis}[%
    width= 0.7\textwidth,
    height= 0.75\textwidth,
    xmin=1e-3,xmax=2e-1,
    ymin=1e-6,ymax=7e-4,
    xtick={1e-4,1e-3,1e-2,1e-1,1e0},
    ytick={1e-6,1e-5,1e-4,1e-3,1e-2,1e-1},
    x dir=reverse,
    legend style={font=\small},
    tick label style={font=\small},
    label style={font=\small},
    xlabel={$\tau$},
    grid=major,
    legend pos = outer north east,
    legend style = {draw},
    every axis plot/.append style={mark=*},
    legend columns=1,
    ]
    \logLogSlopeTriangle{0.3}{0.2}{0.15}{1}{black};

    \addplot [color=mycolor1]
    table[]{plotData/convergenceTest/Polymesher1/displDifferentDiam-1.tsv};
    \addlegendentry{h = 0.67702% , m = -0.00
    }

    \addplot [color=mycolor2]
    table[]{plotData/convergenceTest/Polymesher1/displDifferentDiam-2.tsv};
    \addlegendentry{h = 0.38297% , m = -0.01
    }

    \addplot [color=mycolor3]
    table[]{plotData/convergenceTest/Polymesher1/displDifferentDiam-3.tsv};
    \addlegendentry{h = 0.17723% , m = 0.04
    }

    \addplot [color=mycolor4]
    table[]{plotData/convergenceTest/Polymesher1/displDifferentDiam-4.tsv};
    \addlegendentry{h = 0.10159% , m = 0.09
    }

    \addplot [color=mycolor5]
    table[]{plotData/convergenceTest/Polymesher1/displDifferentDiam-5.tsv};
    \addlegendentry{h = 0.052792% , m = 0.37
    }

    \addplot [color=mycolor6]
    table[]{plotData/convergenceTest/Polymesher1/displDifferentDiam-6.tsv};
    \addlegendentry{h = 0.027202% , m = 0.93
    }

    % \addplot [color=mycolor7, dashed]
    % table[]{plotData/convergenceTest/Polymesher1/displDifferentDiam-7.tsv};
    % \addlegendentry{ref m=1}

  \end{loglogaxis}
\end{tikzpicture}%
    \caption{\texttt{Polymesher1}.}
  \end{subfigure}
  \begin{subfigure}{.49\linewidth}
    \centering
    \begin{tikzpicture}

  \begin{loglogaxis}[%
    width= 0.7\textwidth,
    height= 0.75\textwidth,
    xmin=1e-3,xmax=2e-1,
    ymin=1e-6,ymax=1e-4,
    xtick={1e-4,1e-3,1e-2,1e-1,1e0},
    ytick={1e-6,1e-5,1e-4,1e-3,1e-2,1e-1},
    x dir=reverse,
    legend style={font=\small},
    tick label style={font=\small},
    label style={font=\small},
    xlabel={$\tau$},
    grid=major,
    legend pos = outer north east,
    legend style = {draw},
    every axis plot/.append style={mark=*},
    legend columns=1,
    ]
    \logLogSlopeTriangle{0.3}{0.2}{0.15}{1}{black};

    \addplot [color=mycolor1]
    table[]{plotData/convergenceTest/Polymesher20/displDifferentDiam-1.tsv};
    \addlegendentry{h = 0.47563% , m = 0.01
    }

    \addplot [color=mycolor2]
    table[]{plotData/convergenceTest/Polymesher20/displDifferentDiam-2.tsv};
    \addlegendentry{h = 0.26968% , m = 0.01
    }

    \addplot [color=mycolor3]
    table[]{plotData/convergenceTest/Polymesher20/displDifferentDiam-3.tsv};
    \addlegendentry{h = 0.13683% , m = 0.08
    }

    \addplot [color=mycolor4]
    table[]{plotData/convergenceTest/Polymesher20/displDifferentDiam-4.tsv};
    \addlegendentry{h = 0.069522% , m = 0.22
    }

    \addplot [color=mycolor5]
    table[]{plotData/convergenceTest/Polymesher20/displDifferentDiam-5.tsv};
    \addlegendentry{h = 0.037938% , m = 0.66
    }

    \addplot [color=mycolor6]
    table[]{plotData/convergenceTest/Polymesher20/displDifferentDiam-6.tsv};
    \addlegendentry{h = 0.017468% , m = 0.94
    }

    % \addplot [color=mycolor7, dashed]
    % table[]{plotData/convergenceTest/Polymesher20/displDifferentDiam-7.tsv};
    % \addlegendentry{ref m=1}

  \end{loglogaxis}
\end{tikzpicture}%
    \caption{\texttt{Polymesher20}.}
  \end{subfigure}
  \caption{Test with manufactured exact solution: displacement error behavior through space refinements}
  \label{fig:displacementDifferentDiams}
\end{figure}
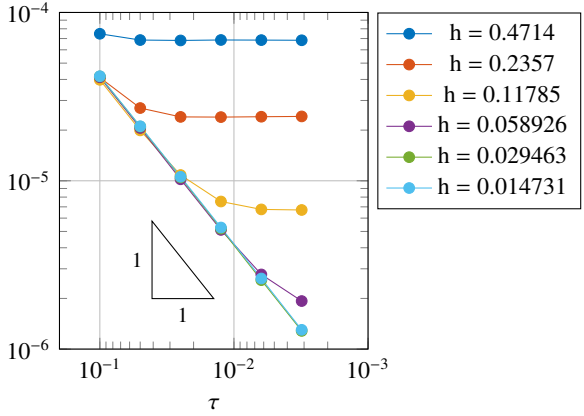
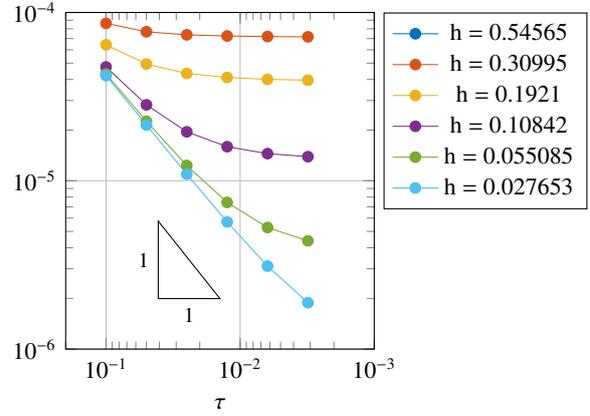
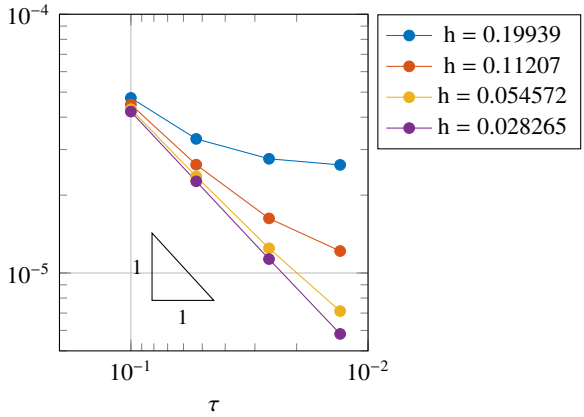
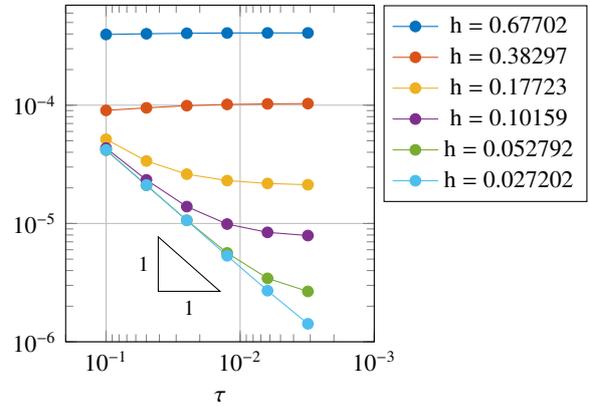
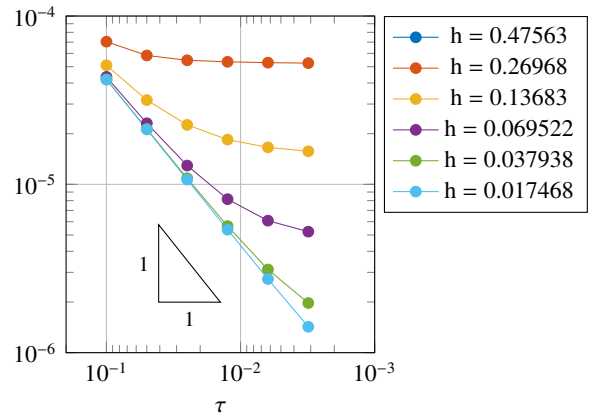

\begin{figure}
  \centering
  \begin{subfigure}{.32\linewidth}
    \centering \begin{tikzpicture}

  \begin{axis}[%
    width=\linewidth,
    height=\linewidth,
    xmode=log,
    xmin=0.01,
    xmax=1,
    xminorticks=true,
    tick align=outside,
    xlabel style={font=\color{white!15!black}},
    xlabel={$h$},
    ymode=log,
    ymin=0.003125,
    ymax=0.1,
    yminorticks=true,
    ylabel style={font=\color{white!15!black}},
    ylabel={$\tau$},
    zmode=log,
    zmin=1e-05,
    zmax=0.00148356799423852,
    zminorticks=true,
    view={-37.5}{30},
    axis background/.style={fill=white},
    axis x line*=bottom,
    axis y line*=left,
    axis z line*=left,
    legend style={at={(1.03,1)}, anchor=north west, legend cell align=left, align=left, draw=white!15!black}
    ]

\addplot3[%
surf,
shader=flat corner, draw=black, z buffer=sort, colormap={mymap}{[1pt] rgb(0pt)=(0.2422,0.1504,0.6603); rgb(1pt)=(0.25039,0.164995,0.707614); rgb(2pt)=(0.257771,0.181781,0.751138); rgb(3pt)=(0.264729,0.197757,0.795214); rgb(4pt)=(0.270648,0.214676,0.836371); rgb(5pt)=(0.275114,0.234238,0.870986); rgb(6pt)=(0.2783,0.255871,0.899071); rgb(7pt)=(0.280333,0.278233,0.9221); rgb(8pt)=(0.281338,0.300595,0.941376); rgb(9pt)=(0.281014,0.322757,0.957886); rgb(10pt)=(0.279467,0.344671,0.971676); rgb(11pt)=(0.275971,0.366681,0.982905); rgb(12pt)=(0.269914,0.3892,0.9906); rgb(13pt)=(0.260243,0.412329,0.995157); rgb(14pt)=(0.244033,0.435833,0.998833); rgb(15pt)=(0.220643,0.460257,0.997286); rgb(16pt)=(0.196333,0.484719,0.989152); rgb(17pt)=(0.183405,0.507371,0.979795); rgb(18pt)=(0.178643,0.528857,0.968157); rgb(19pt)=(0.176438,0.549905,0.952019); rgb(20pt)=(0.168743,0.570262,0.935871); rgb(21pt)=(0.154,0.5902,0.9218); rgb(22pt)=(0.146029,0.609119,0.907857); rgb(23pt)=(0.138024,0.627629,0.89729); rgb(24pt)=(0.124814,0.645929,0.888343); rgb(25pt)=(0.111252,0.6635,0.876314); rgb(26pt)=(0.0952095,0.679829,0.859781); rgb(27pt)=(0.0688714,0.694771,0.839357); rgb(28pt)=(0.0296667,0.708167,0.816333); rgb(29pt)=(0.00357143,0.720267,0.7917); rgb(30pt)=(0.00665714,0.731214,0.766014); rgb(31pt)=(0.0433286,0.741095,0.73941); rgb(32pt)=(0.0963952,0.75,0.712038); rgb(33pt)=(0.140771,0.7584,0.684157); rgb(34pt)=(0.1717,0.766962,0.655443); rgb(35pt)=(0.193767,0.775767,0.6251); rgb(36pt)=(0.216086,0.7843,0.5923); rgb(37pt)=(0.246957,0.791795,0.556743); rgb(38pt)=(0.290614,0.79729,0.518829); rgb(39pt)=(0.340643,0.8008,0.478857); rgb(40pt)=(0.3909,0.802871,0.435448); rgb(41pt)=(0.445629,0.802419,0.390919); rgb(42pt)=(0.5044,0.7993,0.348); rgb(43pt)=(0.561562,0.794233,0.304481); rgb(44pt)=(0.617395,0.787619,0.261238); rgb(45pt)=(0.671986,0.779271,0.2227); rgb(46pt)=(0.7242,0.769843,0.191029); rgb(47pt)=(0.773833,0.759805,0.16461); rgb(48pt)=(0.820314,0.749814,0.153529); rgb(49pt)=(0.863433,0.7406,0.159633); rgb(50pt)=(0.903543,0.733029,0.177414); rgb(51pt)=(0.939257,0.728786,0.209957); rgb(52pt)=(0.972757,0.729771,0.239443); rgb(53pt)=(0.995648,0.743371,0.237148); rgb(54pt)=(0.996986,0.765857,0.219943); rgb(55pt)=(0.995205,0.789252,0.202762); rgb(56pt)=(0.9892,0.813567,0.188533); rgb(57pt)=(0.978629,0.838629,0.176557); rgb(58pt)=(0.967648,0.8639,0.16429); rgb(59pt)=(0.96101,0.889019,0.153676); rgb(60pt)=(0.959671,0.913457,0.142257); rgb(61pt)=(0.962795,0.937338,0.12651); rgb(62pt)=(0.969114,0.960629,0.106362); rgb(63pt)=(0.9769,0.9839,0.0805)}, mesh/rows=6]
table[point meta=\thisrow{c}] {%
plotData/convergenceTest/Cartesian/pressureSurf-1.tsv};

  \end{axis}
\end{tikzpicture}%
    \caption{$e_{p}$}
  \end{subfigure}
  \begin{subfigure}{.32\linewidth}
    \centering % This file was created by matlab2tikz.
%
\begin{tikzpicture}

\begin{axis}[%
width=\linewidth,
height=\linewidth,
xmode=log,
xmin=0.01,
xmax=1,
xminorticks=true,
tick align=outside,
xlabel style={font=\color{white!15!black}},
xlabel={$h$},
ymode=log,
ymin=0.003125,
ymax=0.1,
yminorticks=true,
ylabel style={font=\color{white!15!black}},
ylabel={$\tau$},
zmode=log,
zmin=8.69464144916329e-06,
zmax=0.001,
zminorticks=true,
view={-37.5}{30},
axis background/.style={fill=white},
axis x line*=bottom,
axis y line*=left,
axis z line*=left,
legend style={at={(1.03,1)}, anchor=north west, legend cell align=left, align=left, draw=white!15!black}
]

\addplot3[%
surf,
shader=flat corner, draw=black, z buffer=sort, colormap={mymap}{[1pt] rgb(0pt)=(0.2422,0.1504,0.6603); rgb(1pt)=(0.25039,0.164995,0.707614); rgb(2pt)=(0.257771,0.181781,0.751138); rgb(3pt)=(0.264729,0.197757,0.795214); rgb(4pt)=(0.270648,0.214676,0.836371); rgb(5pt)=(0.275114,0.234238,0.870986); rgb(6pt)=(0.2783,0.255871,0.899071); rgb(7pt)=(0.280333,0.278233,0.9221); rgb(8pt)=(0.281338,0.300595,0.941376); rgb(9pt)=(0.281014,0.322757,0.957886); rgb(10pt)=(0.279467,0.344671,0.971676); rgb(11pt)=(0.275971,0.366681,0.982905); rgb(12pt)=(0.269914,0.3892,0.9906); rgb(13pt)=(0.260243,0.412329,0.995157); rgb(14pt)=(0.244033,0.435833,0.998833); rgb(15pt)=(0.220643,0.460257,0.997286); rgb(16pt)=(0.196333,0.484719,0.989152); rgb(17pt)=(0.183405,0.507371,0.979795); rgb(18pt)=(0.178643,0.528857,0.968157); rgb(19pt)=(0.176438,0.549905,0.952019); rgb(20pt)=(0.168743,0.570262,0.935871); rgb(21pt)=(0.154,0.5902,0.9218); rgb(22pt)=(0.146029,0.609119,0.907857); rgb(23pt)=(0.138024,0.627629,0.89729); rgb(24pt)=(0.124814,0.645929,0.888343); rgb(25pt)=(0.111252,0.6635,0.876314); rgb(26pt)=(0.0952095,0.679829,0.859781); rgb(27pt)=(0.0688714,0.694771,0.839357); rgb(28pt)=(0.0296667,0.708167,0.816333); rgb(29pt)=(0.00357143,0.720267,0.7917); rgb(30pt)=(0.00665714,0.731214,0.766014); rgb(31pt)=(0.0433286,0.741095,0.73941); rgb(32pt)=(0.0963952,0.75,0.712038); rgb(33pt)=(0.140771,0.7584,0.684157); rgb(34pt)=(0.1717,0.766962,0.655443); rgb(35pt)=(0.193767,0.775767,0.6251); rgb(36pt)=(0.216086,0.7843,0.5923); rgb(37pt)=(0.246957,0.791795,0.556743); rgb(38pt)=(0.290614,0.79729,0.518829); rgb(39pt)=(0.340643,0.8008,0.478857); rgb(40pt)=(0.3909,0.802871,0.435448); rgb(41pt)=(0.445629,0.802419,0.390919); rgb(42pt)=(0.5044,0.7993,0.348); rgb(43pt)=(0.561562,0.794233,0.304481); rgb(44pt)=(0.617395,0.787619,0.261238); rgb(45pt)=(0.671986,0.779271,0.2227); rgb(46pt)=(0.7242,0.769843,0.191029); rgb(47pt)=(0.773833,0.759805,0.16461); rgb(48pt)=(0.820314,0.749814,0.153529); rgb(49pt)=(0.863433,0.7406,0.159633); rgb(50pt)=(0.903543,0.733029,0.177414); rgb(51pt)=(0.939257,0.728786,0.209957); rgb(52pt)=(0.972757,0.729771,0.239443); rgb(53pt)=(0.995648,0.743371,0.237148); rgb(54pt)=(0.996986,0.765857,0.219943); rgb(55pt)=(0.995205,0.789252,0.202762); rgb(56pt)=(0.9892,0.813567,0.188533); rgb(57pt)=(0.978629,0.838629,0.176557); rgb(58pt)=(0.967648,0.8639,0.16429); rgb(59pt)=(0.96101,0.889019,0.153676); rgb(60pt)=(0.959671,0.913457,0.142257); rgb(61pt)=(0.962795,0.937338,0.12651); rgb(62pt)=(0.969114,0.960629,0.106362); rgb(63pt)=(0.9769,0.9839,0.0805)}, mesh/rows=6]
table[point meta=\thisrow{c}] {%
plotData/convergenceTest/Cartesian/energySurf-1.tsv};

\end{axis}
\end{tikzpicture}%
    \caption{$e_{\stress}$}
  \end{subfigure}
  \begin{subfigure}{.32\linewidth}
    \centering \begin{tikzpicture}

\begin{axis}[%
width=\linewidth,
height=\linewidth,
xmode=log,
xmin=0.01,
xmax=1,
xminorticks=true,
tick align=outside,
xlabel style={font=\color{white!15!black}},
xlabel={$h$},
ymode=log,
ymin=0.003125,
ymax=0.1,
yminorticks=true,
ylabel style={font=\color{white!15!black}},
ylabel={$\tau$},
zmode=log,
zmin=1e-06,
zmax=0.0001,
zminorticks=true,
view={-37.5}{30},
axis background/.style={fill=white},
axis x line*=bottom,
axis y line*=left,
axis z line*=left,
legend style={at={(1.03,1)}, anchor=north west, legend cell align=left, align=left, draw=white!15!black}
]

\addplot3[%
surf,
shader=flat corner, draw=black, z buffer=sort, colormap={mymap}{[1pt] rgb(0pt)=(0.2422,0.1504,0.6603); rgb(1pt)=(0.25039,0.164995,0.707614); rgb(2pt)=(0.257771,0.181781,0.751138); rgb(3pt)=(0.264729,0.197757,0.795214); rgb(4pt)=(0.270648,0.214676,0.836371); rgb(5pt)=(0.275114,0.234238,0.870986); rgb(6pt)=(0.2783,0.255871,0.899071); rgb(7pt)=(0.280333,0.278233,0.9221); rgb(8pt)=(0.281338,0.300595,0.941376); rgb(9pt)=(0.281014,0.322757,0.957886); rgb(10pt)=(0.279467,0.344671,0.971676); rgb(11pt)=(0.275971,0.366681,0.982905); rgb(12pt)=(0.269914,0.3892,0.9906); rgb(13pt)=(0.260243,0.412329,0.995157); rgb(14pt)=(0.244033,0.435833,0.998833); rgb(15pt)=(0.220643,0.460257,0.997286); rgb(16pt)=(0.196333,0.484719,0.989152); rgb(17pt)=(0.183405,0.507371,0.979795); rgb(18pt)=(0.178643,0.528857,0.968157); rgb(19pt)=(0.176438,0.549905,0.952019); rgb(20pt)=(0.168743,0.570262,0.935871); rgb(21pt)=(0.154,0.5902,0.9218); rgb(22pt)=(0.146029,0.609119,0.907857); rgb(23pt)=(0.138024,0.627629,0.89729); rgb(24pt)=(0.124814,0.645929,0.888343); rgb(25pt)=(0.111252,0.6635,0.876314); rgb(26pt)=(0.0952095,0.679829,0.859781); rgb(27pt)=(0.0688714,0.694771,0.839357); rgb(28pt)=(0.0296667,0.708167,0.816333); rgb(29pt)=(0.00357143,0.720267,0.7917); rgb(30pt)=(0.00665714,0.731214,0.766014); rgb(31pt)=(0.0433286,0.741095,0.73941); rgb(32pt)=(0.0963952,0.75,0.712038); rgb(33pt)=(0.140771,0.7584,0.684157); rgb(34pt)=(0.1717,0.766962,0.655443); rgb(35pt)=(0.193767,0.775767,0.6251); rgb(36pt)=(0.216086,0.7843,0.5923); rgb(37pt)=(0.246957,0.791795,0.556743); rgb(38pt)=(0.290614,0.79729,0.518829); rgb(39pt)=(0.340643,0.8008,0.478857); rgb(40pt)=(0.3909,0.802871,0.435448); rgb(41pt)=(0.445629,0.802419,0.390919); rgb(42pt)=(0.5044,0.7993,0.348); rgb(43pt)=(0.561562,0.794233,0.304481); rgb(44pt)=(0.617395,0.787619,0.261238); rgb(45pt)=(0.671986,0.779271,0.2227); rgb(46pt)=(0.7242,0.769843,0.191029); rgb(47pt)=(0.773833,0.759805,0.16461); rgb(48pt)=(0.820314,0.749814,0.153529); rgb(49pt)=(0.863433,0.7406,0.159633); rgb(50pt)=(0.903543,0.733029,0.177414); rgb(51pt)=(0.939257,0.728786,0.209957); rgb(52pt)=(0.972757,0.729771,0.239443); rgb(53pt)=(0.995648,0.743371,0.237148); rgb(54pt)=(0.996986,0.765857,0.219943); rgb(55pt)=(0.995205,0.789252,0.202762); rgb(56pt)=(0.9892,0.813567,0.188533); rgb(57pt)=(0.978629,0.838629,0.176557); rgb(58pt)=(0.967648,0.8639,0.16429); rgb(59pt)=(0.96101,0.889019,0.153676); rgb(60pt)=(0.959671,0.913457,0.142257); rgb(61pt)=(0.962795,0.937338,0.12651); rgb(62pt)=(0.969114,0.960629,0.106362); rgb(63pt)=(0.9769,0.9839,0.0805)}, mesh/rows=6]
table[point meta=\thisrow{c}] {%
plotData/convergenceTest/Cartesian/displSurf-1.tsv};

\end{axis}
\end{tikzpicture}%
    \caption{$e_{\displ}$}
  \end{subfigure}
  \caption{Test with manufactured exact solution: error through space and time refinements, for the \texttt{Cartesian} mesh family.}
  \label{fig:cartsurf}
\end{figure}
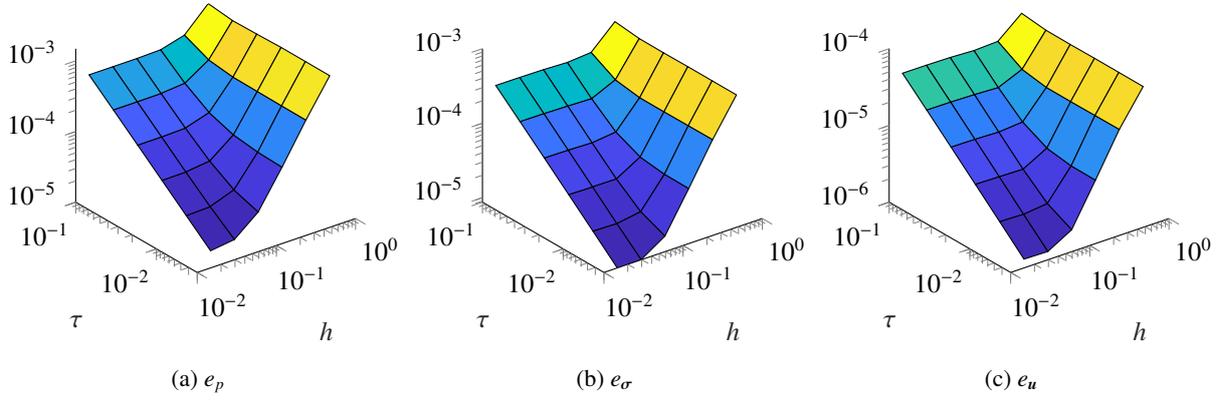

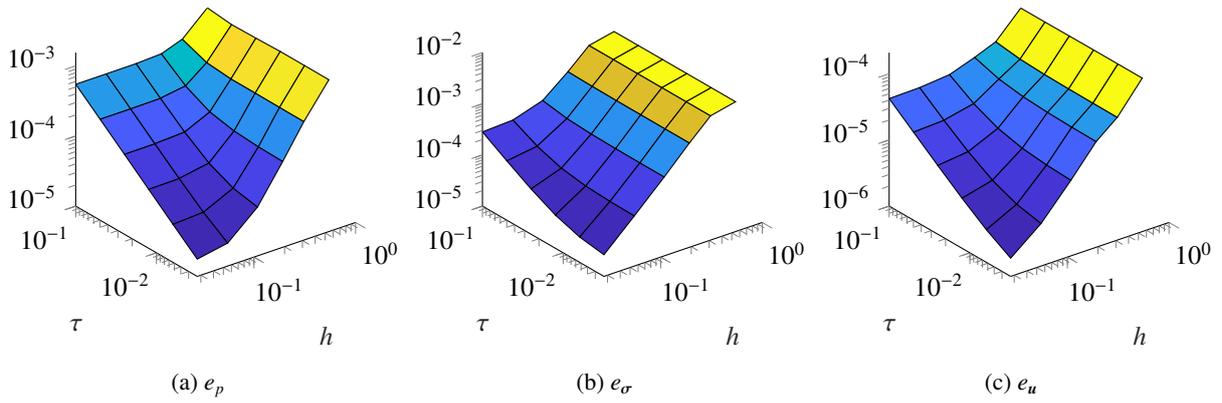
\begin{figure}
  \centering
  \begin{subfigure}{.32\linewidth}
    \centering \begin{tikzpicture}

  \begin{axis}[%
    width=\linewidth,
    height=\linewidth,
    xmode=log,
    xmin=0.0276528251152793,
    xmax=1,
    xminorticks=true,
    tick align=outside,
    xlabel style={font=\color{white!15!black}},
    xlabel={$h$},
    ymode=log,
    ymin=0.003125,
    ymax=0.1,
    yminorticks=true,
    ylabel style={font=\color{white!15!black}},
    ylabel={$\tau$},
    zmode=log,
    zmin=1e-05,
    zmax=0.00154467896794605,
    zminorticks=true,
    view={-37.5}{30},
    axis background/.style={fill=white},
    axis x line*=bottom,
    axis y line*=left,
    axis z line*=left,
    legend style={at={(1.03,1)}, anchor=north west, legend cell align=left, align=left, draw=white!15!black}
    ]

\addplot3[%
surf,
shader=flat corner, draw=black, z buffer=sort, colormap={mymap}{[1pt] rgb(0pt)=(0.2422,0.1504,0.6603); rgb(1pt)=(0.25039,0.164995,0.707614); rgb(2pt)=(0.257771,0.181781,0.751138); rgb(3pt)=(0.264729,0.197757,0.795214); rgb(4pt)=(0.270648,0.214676,0.836371); rgb(5pt)=(0.275114,0.234238,0.870986); rgb(6pt)=(0.2783,0.255871,0.899071); rgb(7pt)=(0.280333,0.278233,0.9221); rgb(8pt)=(0.281338,0.300595,0.941376); rgb(9pt)=(0.281014,0.322757,0.957886); rgb(10pt)=(0.279467,0.344671,0.971676); rgb(11pt)=(0.275971,0.366681,0.982905); rgb(12pt)=(0.269914,0.3892,0.9906); rgb(13pt)=(0.260243,0.412329,0.995157); rgb(14pt)=(0.244033,0.435833,0.998833); rgb(15pt)=(0.220643,0.460257,0.997286); rgb(16pt)=(0.196333,0.484719,0.989152); rgb(17pt)=(0.183405,0.507371,0.979795); rgb(18pt)=(0.178643,0.528857,0.968157); rgb(19pt)=(0.176438,0.549905,0.952019); rgb(20pt)=(0.168743,0.570262,0.935871); rgb(21pt)=(0.154,0.5902,0.9218); rgb(22pt)=(0.146029,0.609119,0.907857); rgb(23pt)=(0.138024,0.627629,0.89729); rgb(24pt)=(0.124814,0.645929,0.888343); rgb(25pt)=(0.111252,0.6635,0.876314); rgb(26pt)=(0.0952095,0.679829,0.859781); rgb(27pt)=(0.0688714,0.694771,0.839357); rgb(28pt)=(0.0296667,0.708167,0.816333); rgb(29pt)=(0.00357143,0.720267,0.7917); rgb(30pt)=(0.00665714,0.731214,0.766014); rgb(31pt)=(0.0433286,0.741095,0.73941); rgb(32pt)=(0.0963952,0.75,0.712038); rgb(33pt)=(0.140771,0.7584,0.684157); rgb(34pt)=(0.1717,0.766962,0.655443); rgb(35pt)=(0.193767,0.775767,0.6251); rgb(36pt)=(0.216086,0.7843,0.5923); rgb(37pt)=(0.246957,0.791795,0.556743); rgb(38pt)=(0.290614,0.79729,0.518829); rgb(39pt)=(0.340643,0.8008,0.478857); rgb(40pt)=(0.3909,0.802871,0.435448); rgb(41pt)=(0.445629,0.802419,0.390919); rgb(42pt)=(0.5044,0.7993,0.348); rgb(43pt)=(0.561562,0.794233,0.304481); rgb(44pt)=(0.617395,0.787619,0.261238); rgb(45pt)=(0.671986,0.779271,0.2227); rgb(46pt)=(0.7242,0.769843,0.191029); rgb(47pt)=(0.773833,0.759805,0.16461); rgb(48pt)=(0.820314,0.749814,0.153529); rgb(49pt)=(0.863433,0.7406,0.159633); rgb(50pt)=(0.903543,0.733029,0.177414); rgb(51pt)=(0.939257,0.728786,0.209957); rgb(52pt)=(0.972757,0.729771,0.239443); rgb(53pt)=(0.995648,0.743371,0.237148); rgb(54pt)=(0.996986,0.765857,0.219943); rgb(55pt)=(0.995205,0.789252,0.202762); rgb(56pt)=(0.9892,0.813567,0.188533); rgb(57pt)=(0.978629,0.838629,0.176557); rgb(58pt)=(0.967648,0.8639,0.16429); rgb(59pt)=(0.96101,0.889019,0.153676); rgb(60pt)=(0.959671,0.913457,0.142257); rgb(61pt)=(0.962795,0.937338,0.12651); rgb(62pt)=(0.969114,0.960629,0.106362); rgb(63pt)=(0.9769,0.9839,0.0805)}, mesh/rows=6]
table[point meta=\thisrow{c}] {%
plotData/convergenceTest/Skewed/pressureSurf-1.tsv};

  \end{axis}
\end{tikzpicture}%
    \caption{$e_{p}$}
  \end{subfigure}
  \begin{subfigure}{.32\linewidth}
    \centering \begin{tikzpicture}

  \begin{axis}[%
    width=\linewidth,
    height=\linewidth,
    xmode=log,
    xmin=0.0276528251152793,
    xmax=1,
    xminorticks=true,
    tick align=outside,
    xlabel style={font=\color{white!15!black}},
    xlabel={$h$},
    ymode=log,
    ymin=0.003125,
    ymax=0.1,
    yminorticks=true,
    ylabel style={font=\color{white!15!black}},
    ylabel={$\tau$},
    zmode=log,
    zmin=1e-05,
    zmax=0.01,
    zminorticks=true,
    view={-37.5}{30},
    axis background/.style={fill=white},
    axis x line*=bottom,
    axis y line*=left,
    axis z line*=left,
    legend style={at={(1.03,1)}, anchor=north west, legend cell align=left, align=left, draw=white!15!black}
    ]

\addplot3[%
surf,
shader=flat corner, draw=black, z buffer=sort, colormap={mymap}{[1pt] rgb(0pt)=(0.2422,0.1504,0.6603); rgb(1pt)=(0.25039,0.164995,0.707614); rgb(2pt)=(0.257771,0.181781,0.751138); rgb(3pt)=(0.264729,0.197757,0.795214); rgb(4pt)=(0.270648,0.214676,0.836371); rgb(5pt)=(0.275114,0.234238,0.870986); rgb(6pt)=(0.2783,0.255871,0.899071); rgb(7pt)=(0.280333,0.278233,0.9221); rgb(8pt)=(0.281338,0.300595,0.941376); rgb(9pt)=(0.281014,0.322757,0.957886); rgb(10pt)=(0.279467,0.344671,0.971676); rgb(11pt)=(0.275971,0.366681,0.982905); rgb(12pt)=(0.269914,0.3892,0.9906); rgb(13pt)=(0.260243,0.412329,0.995157); rgb(14pt)=(0.244033,0.435833,0.998833); rgb(15pt)=(0.220643,0.460257,0.997286); rgb(16pt)=(0.196333,0.484719,0.989152); rgb(17pt)=(0.183405,0.507371,0.979795); rgb(18pt)=(0.178643,0.528857,0.968157); rgb(19pt)=(0.176438,0.549905,0.952019); rgb(20pt)=(0.168743,0.570262,0.935871); rgb(21pt)=(0.154,0.5902,0.9218); rgb(22pt)=(0.146029,0.609119,0.907857); rgb(23pt)=(0.138024,0.627629,0.89729); rgb(24pt)=(0.124814,0.645929,0.888343); rgb(25pt)=(0.111252,0.6635,0.876314); rgb(26pt)=(0.0952095,0.679829,0.859781); rgb(27pt)=(0.0688714,0.694771,0.839357); rgb(28pt)=(0.0296667,0.708167,0.816333); rgb(29pt)=(0.00357143,0.720267,0.7917); rgb(30pt)=(0.00665714,0.731214,0.766014); rgb(31pt)=(0.0433286,0.741095,0.73941); rgb(32pt)=(0.0963952,0.75,0.712038); rgb(33pt)=(0.140771,0.7584,0.684157); rgb(34pt)=(0.1717,0.766962,0.655443); rgb(35pt)=(0.193767,0.775767,0.6251); rgb(36pt)=(0.216086,0.7843,0.5923); rgb(37pt)=(0.246957,0.791795,0.556743); rgb(38pt)=(0.290614,0.79729,0.518829); rgb(39pt)=(0.340643,0.8008,0.478857); rgb(40pt)=(0.3909,0.802871,0.435448); rgb(41pt)=(0.445629,0.802419,0.390919); rgb(42pt)=(0.5044,0.7993,0.348); rgb(43pt)=(0.561562,0.794233,0.304481); rgb(44pt)=(0.617395,0.787619,0.261238); rgb(45pt)=(0.671986,0.779271,0.2227); rgb(46pt)=(0.7242,0.769843,0.191029); rgb(47pt)=(0.773833,0.759805,0.16461); rgb(48pt)=(0.820314,0.749814,0.153529); rgb(49pt)=(0.863433,0.7406,0.159633); rgb(50pt)=(0.903543,0.733029,0.177414); rgb(51pt)=(0.939257,0.728786,0.209957); rgb(52pt)=(0.972757,0.729771,0.239443); rgb(53pt)=(0.995648,0.743371,0.237148); rgb(54pt)=(0.996986,0.765857,0.219943); rgb(55pt)=(0.995205,0.789252,0.202762); rgb(56pt)=(0.9892,0.813567,0.188533); rgb(57pt)=(0.978629,0.838629,0.176557); rgb(58pt)=(0.967648,0.8639,0.16429); rgb(59pt)=(0.96101,0.889019,0.153676); rgb(60pt)=(0.959671,0.913457,0.142257); rgb(61pt)=(0.962795,0.937338,0.12651); rgb(62pt)=(0.969114,0.960629,0.106362); rgb(63pt)=(0.9769,0.9839,0.0805)}, mesh/rows=6]
table[point meta=\thisrow{c}] {%
plotData/convergenceTest/Skewed/energySurf-1.tsv};

  \end{axis}
\end{tikzpicture}%
    \caption{$e_{\stress}$}
  \end{subfigure}
  \begin{subfigure}{.32\linewidth}
    \centering \begin{tikzpicture}

  \begin{axis}[%
    width=\linewidth,
    height=\linewidth,
    xmode=log,
    xmin=0.0276528251152793,
    xmax=1,
    xminorticks=true,
    tick align=outside,
    xlabel style={font=\color{white!15!black}},
    xlabel={$h$},
    ymode=log,
    ymin=0.003125,
    ymax=0.1,
    yminorticks=true,
    ylabel style={font=\color{white!15!black}},
    ylabel={$\tau$},
    zmode=log,
    zmin=1e-06,
    zmax=0.000208726292749776,
    zminorticks=true,
    view={-37.5}{30},
    axis background/.style={fill=white},
    axis x line*=bottom,
    axis y line*=left,
    axis z line*=left,
    legend style={at={(1.03,1)}, anchor=north west, legend cell align=left, align=left, draw=white!15!black}
    ]

\addplot3[%
surf,
shader=flat corner, draw=black, z buffer=sort, colormap={mymap}{[1pt] rgb(0pt)=(0.2422,0.1504,0.6603); rgb(1pt)=(0.25039,0.164995,0.707614); rgb(2pt)=(0.257771,0.181781,0.751138); rgb(3pt)=(0.264729,0.197757,0.795214); rgb(4pt)=(0.270648,0.214676,0.836371); rgb(5pt)=(0.275114,0.234238,0.870986); rgb(6pt)=(0.2783,0.255871,0.899071); rgb(7pt)=(0.280333,0.278233,0.9221); rgb(8pt)=(0.281338,0.300595,0.941376); rgb(9pt)=(0.281014,0.322757,0.957886); rgb(10pt)=(0.279467,0.344671,0.971676); rgb(11pt)=(0.275971,0.366681,0.982905); rgb(12pt)=(0.269914,0.3892,0.9906); rgb(13pt)=(0.260243,0.412329,0.995157); rgb(14pt)=(0.244033,0.435833,0.998833); rgb(15pt)=(0.220643,0.460257,0.997286); rgb(16pt)=(0.196333,0.484719,0.989152); rgb(17pt)=(0.183405,0.507371,0.979795); rgb(18pt)=(0.178643,0.528857,0.968157); rgb(19pt)=(0.176438,0.549905,0.952019); rgb(20pt)=(0.168743,0.570262,0.935871); rgb(21pt)=(0.154,0.5902,0.9218); rgb(22pt)=(0.146029,0.609119,0.907857); rgb(23pt)=(0.138024,0.627629,0.89729); rgb(24pt)=(0.124814,0.645929,0.888343); rgb(25pt)=(0.111252,0.6635,0.876314); rgb(26pt)=(0.0952095,0.679829,0.859781); rgb(27pt)=(0.0688714,0.694771,0.839357); rgb(28pt)=(0.0296667,0.708167,0.816333); rgb(29pt)=(0.00357143,0.720267,0.7917); rgb(30pt)=(0.00665714,0.731214,0.766014); rgb(31pt)=(0.0433286,0.741095,0.73941); rgb(32pt)=(0.0963952,0.75,0.712038); rgb(33pt)=(0.140771,0.7584,0.684157); rgb(34pt)=(0.1717,0.766962,0.655443); rgb(35pt)=(0.193767,0.775767,0.6251); rgb(36pt)=(0.216086,0.7843,0.5923); rgb(37pt)=(0.246957,0.791795,0.556743); rgb(38pt)=(0.290614,0.79729,0.518829); rgb(39pt)=(0.340643,0.8008,0.478857); rgb(40pt)=(0.3909,0.802871,0.435448); rgb(41pt)=(0.445629,0.802419,0.390919); rgb(42pt)=(0.5044,0.7993,0.348); rgb(43pt)=(0.561562,0.794233,0.304481); rgb(44pt)=(0.617395,0.787619,0.261238); rgb(45pt)=(0.671986,0.779271,0.2227); rgb(46pt)=(0.7242,0.769843,0.191029); rgb(47pt)=(0.773833,0.759805,0.16461); rgb(48pt)=(0.820314,0.749814,0.153529); rgb(49pt)=(0.863433,0.7406,0.159633); rgb(50pt)=(0.903543,0.733029,0.177414); rgb(51pt)=(0.939257,0.728786,0.209957); rgb(52pt)=(0.972757,0.729771,0.239443); rgb(53pt)=(0.995648,0.743371,0.237148); rgb(54pt)=(0.996986,0.765857,0.219943); rgb(55pt)=(0.995205,0.789252,0.202762); rgb(56pt)=(0.9892,0.813567,0.188533); rgb(57pt)=(0.978629,0.838629,0.176557); rgb(58pt)=(0.967648,0.8639,0.16429); rgb(59pt)=(0.96101,0.889019,0.153676); rgb(60pt)=(0.959671,0.913457,0.142257); rgb(61pt)=(0.962795,0.937338,0.12651); rgb(62pt)=(0.969114,0.960629,0.106362); rgb(63pt)=(0.9769,0.9839,0.0805)}, mesh/rows=6]
table[point meta=\thisrow{c}] {%
plotData/convergenceTest/Skewed/displSurf-1.tsv};

  \end{axis}
\end{tikzpicture}%
    \caption{$e_{\displ}$}
  \end{subfigure}
  \caption{Test with manufactured exact solution: error through space and time refinements, for the \texttt{Skewed} mesh family.}
  \label{fig:skewedsurf}
\end{figure}

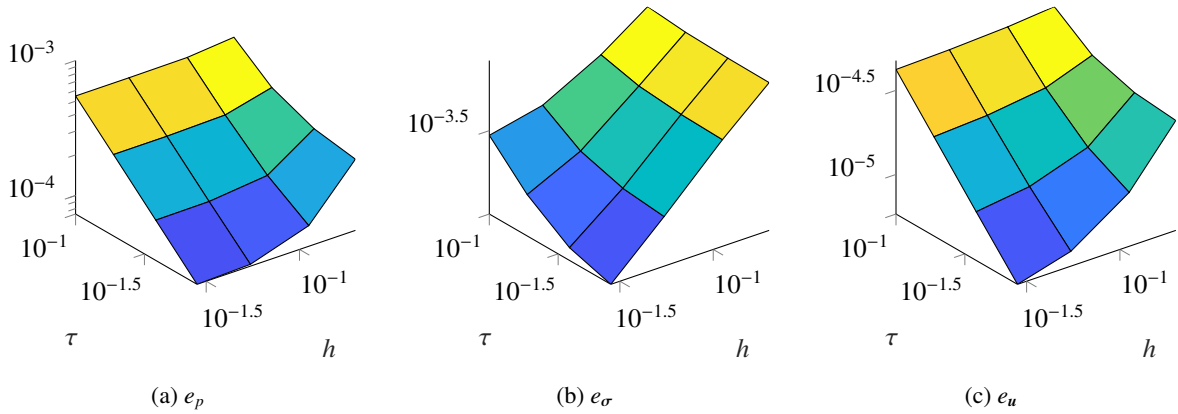
\begin{figure}
  \centering
  \begin{subfigure}{.32\linewidth}
    \centering \begin{tikzpicture}

  \begin{axis}[%
    width=\linewidth,
    height=\linewidth,
    xmode=log,
    xmin=0.0282647121473865,
    xmax=0.199394115718455,
    xminorticks=true,
    tick align=outside,
    xlabel style={font=\color{white!15!black}},
    xlabel={$h$},
    ymode=log,
    ymin=0.01311113494956,
    ymax=0.1,
    yminorticks=true,
    ylabel style={font=\color{white!15!black}},
    ylabel={$\tau$},
    zmode=log,
    zmin=7.28708163869506e-05,
    zmax=0.001,
    zminorticks=true,
    view={-37.5}{30},
    axis background/.style={fill=white},
    axis x line*=bottom,
    axis y line*=left,
    axis z line*=left,
    legend style={at={(1.03,1)}, anchor=north west, legend cell align=left, align=left, draw=white!15!black}
    ]

    \addplot3[%
    surf,
    shader=flat corner, draw=black, z buffer=sort, colormap={mymap}{[1pt] rgb(0pt)=(0.2422,0.1504,0.6603); rgb(1pt)=(0.25039,0.164995,0.707614); rgb(2pt)=(0.257771,0.181781,0.751138); rgb(3pt)=(0.264729,0.197757,0.795214); rgb(4pt)=(0.270648,0.214676,0.836371); rgb(5pt)=(0.275114,0.234238,0.870986); rgb(6pt)=(0.2783,0.255871,0.899071); rgb(7pt)=(0.280333,0.278233,0.9221); rgb(8pt)=(0.281338,0.300595,0.941376); rgb(9pt)=(0.281014,0.322757,0.957886); rgb(10pt)=(0.279467,0.344671,0.971676); rgb(11pt)=(0.275971,0.366681,0.982905); rgb(12pt)=(0.269914,0.3892,0.9906); rgb(13pt)=(0.260243,0.412329,0.995157); rgb(14pt)=(0.244033,0.435833,0.998833); rgb(15pt)=(0.220643,0.460257,0.997286); rgb(16pt)=(0.196333,0.484719,0.989152); rgb(17pt)=(0.183405,0.507371,0.979795); rgb(18pt)=(0.178643,0.528857,0.968157); rgb(19pt)=(0.176438,0.549905,0.952019); rgb(20pt)=(0.168743,0.570262,0.935871); rgb(21pt)=(0.154,0.5902,0.9218); rgb(22pt)=(0.146029,0.609119,0.907857); rgb(23pt)=(0.138024,0.627629,0.89729); rgb(24pt)=(0.124814,0.645929,0.888343); rgb(25pt)=(0.111252,0.6635,0.876314); rgb(26pt)=(0.0952095,0.679829,0.859781); rgb(27pt)=(0.0688714,0.694771,0.839357); rgb(28pt)=(0.0296667,0.708167,0.816333); rgb(29pt)=(0.00357143,0.720267,0.7917); rgb(30pt)=(0.00665714,0.731214,0.766014); rgb(31pt)=(0.0433286,0.741095,0.73941); rgb(32pt)=(0.0963952,0.75,0.712038); rgb(33pt)=(0.140771,0.7584,0.684157); rgb(34pt)=(0.1717,0.766962,0.655443); rgb(35pt)=(0.193767,0.775767,0.6251); rgb(36pt)=(0.216086,0.7843,0.5923); rgb(37pt)=(0.246957,0.791795,0.556743); rgb(38pt)=(0.290614,0.79729,0.518829); rgb(39pt)=(0.340643,0.8008,0.478857); rgb(40pt)=(0.3909,0.802871,0.435448); rgb(41pt)=(0.445629,0.802419,0.390919); rgb(42pt)=(0.5044,0.7993,0.348); rgb(43pt)=(0.561562,0.794233,0.304481); rgb(44pt)=(0.617395,0.787619,0.261238); rgb(45pt)=(0.671986,0.779271,0.2227); rgb(46pt)=(0.7242,0.769843,0.191029); rgb(47pt)=(0.773833,0.759805,0.16461); rgb(48pt)=(0.820314,0.749814,0.153529); rgb(49pt)=(0.863433,0.7406,0.159633); rgb(50pt)=(0.903543,0.733029,0.177414); rgb(51pt)=(0.939257,0.728786,0.209957); rgb(52pt)=(0.972757,0.729771,0.239443); rgb(53pt)=(0.995648,0.743371,0.237148); rgb(54pt)=(0.996986,0.765857,0.219943); rgb(55pt)=(0.995205,0.789252,0.202762); rgb(56pt)=(0.9892,0.813567,0.188533); rgb(57pt)=(0.978629,0.838629,0.176557); rgb(58pt)=(0.967648,0.8639,0.16429); rgb(59pt)=(0.96101,0.889019,0.153676); rgb(60pt)=(0.959671,0.913457,0.142257); rgb(61pt)=(0.962795,0.937338,0.12651); rgb(62pt)=(0.969114,0.960629,0.106362); rgb(63pt)=(0.9769,0.9839,0.0805)}, mesh/rows=4]
    table[point meta=\thisrow{c}] {%
      plotData/convergenceTest/Hybrid/pressureSurf-1.tsv};

  \end{axis}
\end{tikzpicture}%
    \caption{$e_{p}$}
  \end{subfigure}
  \begin{subfigure}{.32\linewidth}
    \centering \begin{tikzpicture}

  \begin{axis}[%
    width=\linewidth,
    height=\linewidth,
    xmode=log,
    xmin=0.0282647121473865,
    xmax=0.199394115718455,
    xminorticks=true,
    tick align=outside,
    xlabel style={font=\color{white!15!black}},
    xlabel={$h$},
    ymode=log,
    ymin=0.01311113494956,
    ymax=0.1,
    yminorticks=true,
    ylabel style={font=\color{white!15!black}},
    ylabel={$\tau$},
    zmode=log,
    zmin=0.000101864803190054,
    zmax=0.000819073139428861,
    zminorticks=true,
    view={-37.5}{30},
    axis background/.style={fill=white},
    axis x line*=bottom,
    axis y line*=left,
    axis z line*=left,
    legend style={at={(1.03,1)}, anchor=north west, legend cell align=left, align=left, draw=white!15!black}
    ]

    \addplot3[%
    surf,
    shader=flat corner, draw=black, z buffer=sort, colormap={mymap}{[1pt] rgb(0pt)=(0.2422,0.1504,0.6603); rgb(1pt)=(0.25039,0.164995,0.707614); rgb(2pt)=(0.257771,0.181781,0.751138); rgb(3pt)=(0.264729,0.197757,0.795214); rgb(4pt)=(0.270648,0.214676,0.836371); rgb(5pt)=(0.275114,0.234238,0.870986); rgb(6pt)=(0.2783,0.255871,0.899071); rgb(7pt)=(0.280333,0.278233,0.9221); rgb(8pt)=(0.281338,0.300595,0.941376); rgb(9pt)=(0.281014,0.322757,0.957886); rgb(10pt)=(0.279467,0.344671,0.971676); rgb(11pt)=(0.275971,0.366681,0.982905); rgb(12pt)=(0.269914,0.3892,0.9906); rgb(13pt)=(0.260243,0.412329,0.995157); rgb(14pt)=(0.244033,0.435833,0.998833); rgb(15pt)=(0.220643,0.460257,0.997286); rgb(16pt)=(0.196333,0.484719,0.989152); rgb(17pt)=(0.183405,0.507371,0.979795); rgb(18pt)=(0.178643,0.528857,0.968157); rgb(19pt)=(0.176438,0.549905,0.952019); rgb(20pt)=(0.168743,0.570262,0.935871); rgb(21pt)=(0.154,0.5902,0.9218); rgb(22pt)=(0.146029,0.609119,0.907857); rgb(23pt)=(0.138024,0.627629,0.89729); rgb(24pt)=(0.124814,0.645929,0.888343); rgb(25pt)=(0.111252,0.6635,0.876314); rgb(26pt)=(0.0952095,0.679829,0.859781); rgb(27pt)=(0.0688714,0.694771,0.839357); rgb(28pt)=(0.0296667,0.708167,0.816333); rgb(29pt)=(0.00357143,0.720267,0.7917); rgb(30pt)=(0.00665714,0.731214,0.766014); rgb(31pt)=(0.0433286,0.741095,0.73941); rgb(32pt)=(0.0963952,0.75,0.712038); rgb(33pt)=(0.140771,0.7584,0.684157); rgb(34pt)=(0.1717,0.766962,0.655443); rgb(35pt)=(0.193767,0.775767,0.6251); rgb(36pt)=(0.216086,0.7843,0.5923); rgb(37pt)=(0.246957,0.791795,0.556743); rgb(38pt)=(0.290614,0.79729,0.518829); rgb(39pt)=(0.340643,0.8008,0.478857); rgb(40pt)=(0.3909,0.802871,0.435448); rgb(41pt)=(0.445629,0.802419,0.390919); rgb(42pt)=(0.5044,0.7993,0.348); rgb(43pt)=(0.561562,0.794233,0.304481); rgb(44pt)=(0.617395,0.787619,0.261238); rgb(45pt)=(0.671986,0.779271,0.2227); rgb(46pt)=(0.7242,0.769843,0.191029); rgb(47pt)=(0.773833,0.759805,0.16461); rgb(48pt)=(0.820314,0.749814,0.153529); rgb(49pt)=(0.863433,0.7406,0.159633); rgb(50pt)=(0.903543,0.733029,0.177414); rgb(51pt)=(0.939257,0.728786,0.209957); rgb(52pt)=(0.972757,0.729771,0.239443); rgb(53pt)=(0.995648,0.743371,0.237148); rgb(54pt)=(0.996986,0.765857,0.219943); rgb(55pt)=(0.995205,0.789252,0.202762); rgb(56pt)=(0.9892,0.813567,0.188533); rgb(57pt)=(0.978629,0.838629,0.176557); rgb(58pt)=(0.967648,0.8639,0.16429); rgb(59pt)=(0.96101,0.889019,0.153676); rgb(60pt)=(0.959671,0.913457,0.142257); rgb(61pt)=(0.962795,0.937338,0.12651); rgb(62pt)=(0.969114,0.960629,0.106362); rgb(63pt)=(0.9769,0.9839,0.0805)}, mesh/rows=4]
    table[point meta=\thisrow{c}] {%
      plotData/convergenceTest/Hybrid/energySurf-1.tsv};

  \end{axis}
\end{tikzpicture}%
    \caption{$e_{\stress}$}
  \end{subfigure}
  \begin{subfigure}{.32\linewidth}
    \centering \begin{tikzpicture}

  \begin{axis}[%
    width=\linewidth,
    height=\linewidth,
    xmode=log,
    xmin=0.0282647121473865,
    xmax=0.199394115718455,
    xminorticks=true,
    tick align=outside,
    xlabel style={font=\color{white!15!black}},
    xlabel={$h$},
    ymode=log,
    ymin=0.01311113494956,
    ymax=0.1,
    yminorticks=true,
    ylabel style={font=\color{white!15!black}},
    ylabel={$\tau$},
    zmode=log,
    zmin=5.81512980289804e-06,
    zmax=4.75072143422008e-05,
    zminorticks=true,
    view={-37.5}{30},
    axis background/.style={fill=white},
    axis x line*=bottom,
    axis y line*=left,
    axis z line*=left,
    legend style={at={(1.03,1)}, anchor=north west, legend cell align=left, align=left, draw=white!15!black}
    ]

    \addplot3[%
    surf,
    shader=flat corner, draw=black, z buffer=sort, colormap={mymap}{[1pt] rgb(0pt)=(0.2422,0.1504,0.6603); rgb(1pt)=(0.25039,0.164995,0.707614); rgb(2pt)=(0.257771,0.181781,0.751138); rgb(3pt)=(0.264729,0.197757,0.795214); rgb(4pt)=(0.270648,0.214676,0.836371); rgb(5pt)=(0.275114,0.234238,0.870986); rgb(6pt)=(0.2783,0.255871,0.899071); rgb(7pt)=(0.280333,0.278233,0.9221); rgb(8pt)=(0.281338,0.300595,0.941376); rgb(9pt)=(0.281014,0.322757,0.957886); rgb(10pt)=(0.279467,0.344671,0.971676); rgb(11pt)=(0.275971,0.366681,0.982905); rgb(12pt)=(0.269914,0.3892,0.9906); rgb(13pt)=(0.260243,0.412329,0.995157); rgb(14pt)=(0.244033,0.435833,0.998833); rgb(15pt)=(0.220643,0.460257,0.997286); rgb(16pt)=(0.196333,0.484719,0.989152); rgb(17pt)=(0.183405,0.507371,0.979795); rgb(18pt)=(0.178643,0.528857,0.968157); rgb(19pt)=(0.176438,0.549905,0.952019); rgb(20pt)=(0.168743,0.570262,0.935871); rgb(21pt)=(0.154,0.5902,0.9218); rgb(22pt)=(0.146029,0.609119,0.907857); rgb(23pt)=(0.138024,0.627629,0.89729); rgb(24pt)=(0.124814,0.645929,0.888343); rgb(25pt)=(0.111252,0.6635,0.876314); rgb(26pt)=(0.0952095,0.679829,0.859781); rgb(27pt)=(0.0688714,0.694771,0.839357); rgb(28pt)=(0.0296667,0.708167,0.816333); rgb(29pt)=(0.00357143,0.720267,0.7917); rgb(30pt)=(0.00665714,0.731214,0.766014); rgb(31pt)=(0.0433286,0.741095,0.73941); rgb(32pt)=(0.0963952,0.75,0.712038); rgb(33pt)=(0.140771,0.7584,0.684157); rgb(34pt)=(0.1717,0.766962,0.655443); rgb(35pt)=(0.193767,0.775767,0.6251); rgb(36pt)=(0.216086,0.7843,0.5923); rgb(37pt)=(0.246957,0.791795,0.556743); rgb(38pt)=(0.290614,0.79729,0.518829); rgb(39pt)=(0.340643,0.8008,0.478857); rgb(40pt)=(0.3909,0.802871,0.435448); rgb(41pt)=(0.445629,0.802419,0.390919); rgb(42pt)=(0.5044,0.7993,0.348); rgb(43pt)=(0.561562,0.794233,0.304481); rgb(44pt)=(0.617395,0.787619,0.261238); rgb(45pt)=(0.671986,0.779271,0.2227); rgb(46pt)=(0.7242,0.769843,0.191029); rgb(47pt)=(0.773833,0.759805,0.16461); rgb(48pt)=(0.820314,0.749814,0.153529); rgb(49pt)=(0.863433,0.7406,0.159633); rgb(50pt)=(0.903543,0.733029,0.177414); rgb(51pt)=(0.939257,0.728786,0.209957); rgb(52pt)=(0.972757,0.729771,0.239443); rgb(53pt)=(0.995648,0.743371,0.237148); rgb(54pt)=(0.996986,0.765857,0.219943); rgb(55pt)=(0.995205,0.789252,0.202762); rgb(56pt)=(0.9892,0.813567,0.188533); rgb(57pt)=(0.978629,0.838629,0.176557); rgb(58pt)=(0.967648,0.8639,0.16429); rgb(59pt)=(0.96101,0.889019,0.153676); rgb(60pt)=(0.959671,0.913457,0.142257); rgb(61pt)=(0.962795,0.937338,0.12651); rgb(62pt)=(0.969114,0.960629,0.106362); rgb(63pt)=(0.9769,0.9839,0.0805)}, mesh/rows=4]
    table[point meta=\thisrow{c}] {%
      plotData/convergenceTest/Hybrid/displSurf-1.tsv};
  \end{axis}
\end{tikzpicture}%
    \caption{$e_{\displ}$}
  \end{subfigure}
  \caption{Test with manufactured exact solution: error through space and time refinements, for the \texttt{Hybrid} mesh family.}
  \label{fig:hybrid}
\end{figure}

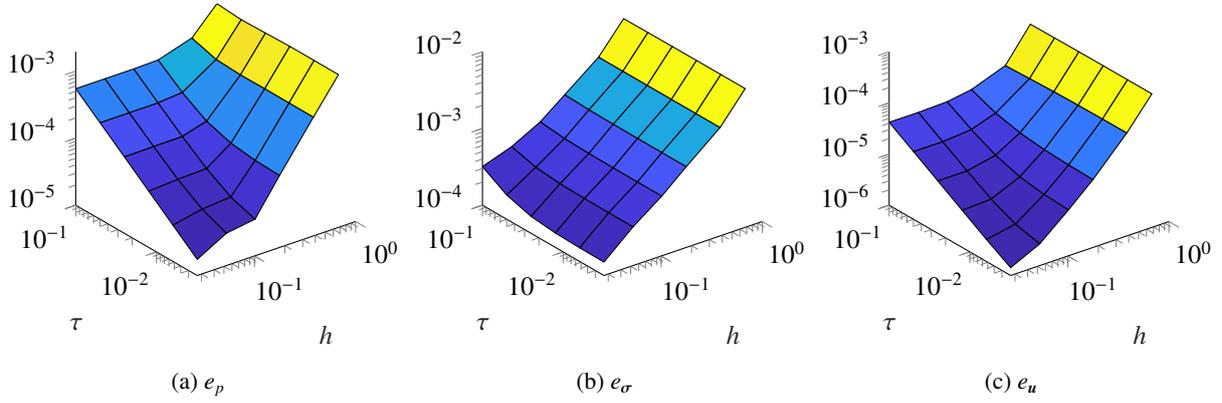
\begin{figure}
  \centering
  \begin{subfigure}{.32\linewidth}
    \centering
\begin{tikzpicture}

\begin{axis}[%
width=\linewidth,
height=\linewidth,
xmode=log,
xmin=0.0272022512213662,
xmax=1,
xminorticks=true,
tick align=outside,
xlabel style={font=\color{white!15!black}},
xlabel={$h$},
ymode=log,
ymin=0.003125,
ymax=0.1,
yminorticks=true,
ylabel style={font=\color{white!15!black}},
ylabel={$\tau$},
zmode=log,
zmin=1e-05,
zmax=0.00194597131782464,
zminorticks=true,
view={-37.5}{30},
axis background/.style={fill=white},
axis x line*=bottom,
axis y line*=left,
axis z line*=left,
legend style={at={(1.03,1)}, anchor=north west, legend cell align=left, align=left, draw=white!15!black}
]

\addplot3[%
surf,
shader=flat corner, draw=black, z buffer=sort, colormap={mymap}{[1pt] rgb(0pt)=(0.2422,0.1504,0.6603); rgb(1pt)=(0.25039,0.164995,0.707614); rgb(2pt)=(0.257771,0.181781,0.751138); rgb(3pt)=(0.264729,0.197757,0.795214); rgb(4pt)=(0.270648,0.214676,0.836371); rgb(5pt)=(0.275114,0.234238,0.870986); rgb(6pt)=(0.2783,0.255871,0.899071); rgb(7pt)=(0.280333,0.278233,0.9221); rgb(8pt)=(0.281338,0.300595,0.941376); rgb(9pt)=(0.281014,0.322757,0.957886); rgb(10pt)=(0.279467,0.344671,0.971676); rgb(11pt)=(0.275971,0.366681,0.982905); rgb(12pt)=(0.269914,0.3892,0.9906); rgb(13pt)=(0.260243,0.412329,0.995157); rgb(14pt)=(0.244033,0.435833,0.998833); rgb(15pt)=(0.220643,0.460257,0.997286); rgb(16pt)=(0.196333,0.484719,0.989152); rgb(17pt)=(0.183405,0.507371,0.979795); rgb(18pt)=(0.178643,0.528857,0.968157); rgb(19pt)=(0.176438,0.549905,0.952019); rgb(20pt)=(0.168743,0.570262,0.935871); rgb(21pt)=(0.154,0.5902,0.9218); rgb(22pt)=(0.146029,0.609119,0.907857); rgb(23pt)=(0.138024,0.627629,0.89729); rgb(24pt)=(0.124814,0.645929,0.888343); rgb(25pt)=(0.111252,0.6635,0.876314); rgb(26pt)=(0.0952095,0.679829,0.859781); rgb(27pt)=(0.0688714,0.694771,0.839357); rgb(28pt)=(0.0296667,0.708167,0.816333); rgb(29pt)=(0.00357143,0.720267,0.7917); rgb(30pt)=(0.00665714,0.731214,0.766014); rgb(31pt)=(0.0433286,0.741095,0.73941); rgb(32pt)=(0.0963952,0.75,0.712038); rgb(33pt)=(0.140771,0.7584,0.684157); rgb(34pt)=(0.1717,0.766962,0.655443); rgb(35pt)=(0.193767,0.775767,0.6251); rgb(36pt)=(0.216086,0.7843,0.5923); rgb(37pt)=(0.246957,0.791795,0.556743); rgb(38pt)=(0.290614,0.79729,0.518829); rgb(39pt)=(0.340643,0.8008,0.478857); rgb(40pt)=(0.3909,0.802871,0.435448); rgb(41pt)=(0.445629,0.802419,0.390919); rgb(42pt)=(0.5044,0.7993,0.348); rgb(43pt)=(0.561562,0.794233,0.304481); rgb(44pt)=(0.617395,0.787619,0.261238); rgb(45pt)=(0.671986,0.779271,0.2227); rgb(46pt)=(0.7242,0.769843,0.191029); rgb(47pt)=(0.773833,0.759805,0.16461); rgb(48pt)=(0.820314,0.749814,0.153529); rgb(49pt)=(0.863433,0.7406,0.159633); rgb(50pt)=(0.903543,0.733029,0.177414); rgb(51pt)=(0.939257,0.728786,0.209957); rgb(52pt)=(0.972757,0.729771,0.239443); rgb(53pt)=(0.995648,0.743371,0.237148); rgb(54pt)=(0.996986,0.765857,0.219943); rgb(55pt)=(0.995205,0.789252,0.202762); rgb(56pt)=(0.9892,0.813567,0.188533); rgb(57pt)=(0.978629,0.838629,0.176557); rgb(58pt)=(0.967648,0.8639,0.16429); rgb(59pt)=(0.96101,0.889019,0.153676); rgb(60pt)=(0.959671,0.913457,0.142257); rgb(61pt)=(0.962795,0.937338,0.12651); rgb(62pt)=(0.969114,0.960629,0.106362); rgb(63pt)=(0.9769,0.9839,0.0805)}, mesh/rows=6]
table[point meta=\thisrow{c}] {%
plotData/convergenceTest/Polymesher1/pressureSurf-1.tsv};

\end{axis}
\end{tikzpicture}%
    \caption{$e_{p}$}
  \end{subfigure}
  \begin{subfigure}{.32\linewidth}
    \centering % This file was created by matlab2tikz.
%
\begin{tikzpicture}

\begin{axis}[%
width=\linewidth,
height=\linewidth,
xmode=log,
xmin=0.0272022512213662,
xmax=1,
xminorticks=true,
tick align=outside,
xlabel style={font=\color{white!15!black}},
xlabel={$h$},
ymode=log,
ymin=0.003125,
ymax=0.1,
yminorticks=true,
ylabel style={font=\color{white!15!black}},
ylabel={$\tau$},
zmode=log,
zmin=0.0001,
zmax=0.01,
zminorticks=true,
view={-37.5}{30},
axis background/.style={fill=white},
axis x line*=bottom,
axis y line*=left,
axis z line*=left,
legend style={at={(1.03,1)}, anchor=north west, legend cell align=left, align=left, draw=white!15!black}
]

\addplot3[%
surf,
shader=flat corner, draw=black, z buffer=sort, colormap={mymap}{[1pt] rgb(0pt)=(0.2422,0.1504,0.6603); rgb(1pt)=(0.25039,0.164995,0.707614); rgb(2pt)=(0.257771,0.181781,0.751138); rgb(3pt)=(0.264729,0.197757,0.795214); rgb(4pt)=(0.270648,0.214676,0.836371); rgb(5pt)=(0.275114,0.234238,0.870986); rgb(6pt)=(0.2783,0.255871,0.899071); rgb(7pt)=(0.280333,0.278233,0.9221); rgb(8pt)=(0.281338,0.300595,0.941376); rgb(9pt)=(0.281014,0.322757,0.957886); rgb(10pt)=(0.279467,0.344671,0.971676); rgb(11pt)=(0.275971,0.366681,0.982905); rgb(12pt)=(0.269914,0.3892,0.9906); rgb(13pt)=(0.260243,0.412329,0.995157); rgb(14pt)=(0.244033,0.435833,0.998833); rgb(15pt)=(0.220643,0.460257,0.997286); rgb(16pt)=(0.196333,0.484719,0.989152); rgb(17pt)=(0.183405,0.507371,0.979795); rgb(18pt)=(0.178643,0.528857,0.968157); rgb(19pt)=(0.176438,0.549905,0.952019); rgb(20pt)=(0.168743,0.570262,0.935871); rgb(21pt)=(0.154,0.5902,0.9218); rgb(22pt)=(0.146029,0.609119,0.907857); rgb(23pt)=(0.138024,0.627629,0.89729); rgb(24pt)=(0.124814,0.645929,0.888343); rgb(25pt)=(0.111252,0.6635,0.876314); rgb(26pt)=(0.0952095,0.679829,0.859781); rgb(27pt)=(0.0688714,0.694771,0.839357); rgb(28pt)=(0.0296667,0.708167,0.816333); rgb(29pt)=(0.00357143,0.720267,0.7917); rgb(30pt)=(0.00665714,0.731214,0.766014); rgb(31pt)=(0.0433286,0.741095,0.73941); rgb(32pt)=(0.0963952,0.75,0.712038); rgb(33pt)=(0.140771,0.7584,0.684157); rgb(34pt)=(0.1717,0.766962,0.655443); rgb(35pt)=(0.193767,0.775767,0.6251); rgb(36pt)=(0.216086,0.7843,0.5923); rgb(37pt)=(0.246957,0.791795,0.556743); rgb(38pt)=(0.290614,0.79729,0.518829); rgb(39pt)=(0.340643,0.8008,0.478857); rgb(40pt)=(0.3909,0.802871,0.435448); rgb(41pt)=(0.445629,0.802419,0.390919); rgb(42pt)=(0.5044,0.7993,0.348); rgb(43pt)=(0.561562,0.794233,0.304481); rgb(44pt)=(0.617395,0.787619,0.261238); rgb(45pt)=(0.671986,0.779271,0.2227); rgb(46pt)=(0.7242,0.769843,0.191029); rgb(47pt)=(0.773833,0.759805,0.16461); rgb(48pt)=(0.820314,0.749814,0.153529); rgb(49pt)=(0.863433,0.7406,0.159633); rgb(50pt)=(0.903543,0.733029,0.177414); rgb(51pt)=(0.939257,0.728786,0.209957); rgb(52pt)=(0.972757,0.729771,0.239443); rgb(53pt)=(0.995648,0.743371,0.237148); rgb(54pt)=(0.996986,0.765857,0.219943); rgb(55pt)=(0.995205,0.789252,0.202762); rgb(56pt)=(0.9892,0.813567,0.188533); rgb(57pt)=(0.978629,0.838629,0.176557); rgb(58pt)=(0.967648,0.8639,0.16429); rgb(59pt)=(0.96101,0.889019,0.153676); rgb(60pt)=(0.959671,0.913457,0.142257); rgb(61pt)=(0.962795,0.937338,0.12651); rgb(62pt)=(0.969114,0.960629,0.106362); rgb(63pt)=(0.9769,0.9839,0.0805)}, mesh/rows=6]
table[point meta=\thisrow{c}] {%
plotData/convergenceTest/Polymesher1/energySurf-1.tsv};

\end{axis}
\end{tikzpicture}%
    \caption{$e_{\stress}$}
  \end{subfigure}
  \begin{subfigure}{.32\linewidth}
    \centering % This file was created by matlab2tikz.
%
\begin{tikzpicture}

\begin{axis}[%
width=\linewidth,
height=\linewidth,
xmode=log,
xmin=0.0272022512213662,
xmax=1,
xminorticks=true,
tick align=outside,
xlabel style={font=\color{white!15!black}},
xlabel={$h$},
ymode=log,
ymin=0.003125,
ymax=0.1,
yminorticks=true,
ylabel style={font=\color{white!15!black}},
ylabel={$\tau$},
zmode=log,
zmin=1e-06,
zmax=0.001,
zminorticks=true,
view={-37.5}{30},
axis background/.style={fill=white},
axis x line*=bottom,
axis y line*=left,
axis z line*=left,
legend style={at={(1.03,1)}, anchor=north west, legend cell align=left, align=left, draw=white!15!black}
]

\addplot3[%
surf,
shader=flat corner, draw=black, z buffer=sort, colormap={mymap}{[1pt] rgb(0pt)=(0.2422,0.1504,0.6603); rgb(1pt)=(0.25039,0.164995,0.707614); rgb(2pt)=(0.257771,0.181781,0.751138); rgb(3pt)=(0.264729,0.197757,0.795214); rgb(4pt)=(0.270648,0.214676,0.836371); rgb(5pt)=(0.275114,0.234238,0.870986); rgb(6pt)=(0.2783,0.255871,0.899071); rgb(7pt)=(0.280333,0.278233,0.9221); rgb(8pt)=(0.281338,0.300595,0.941376); rgb(9pt)=(0.281014,0.322757,0.957886); rgb(10pt)=(0.279467,0.344671,0.971676); rgb(11pt)=(0.275971,0.366681,0.982905); rgb(12pt)=(0.269914,0.3892,0.9906); rgb(13pt)=(0.260243,0.412329,0.995157); rgb(14pt)=(0.244033,0.435833,0.998833); rgb(15pt)=(0.220643,0.460257,0.997286); rgb(16pt)=(0.196333,0.484719,0.989152); rgb(17pt)=(0.183405,0.507371,0.979795); rgb(18pt)=(0.178643,0.528857,0.968157); rgb(19pt)=(0.176438,0.549905,0.952019); rgb(20pt)=(0.168743,0.570262,0.935871); rgb(21pt)=(0.154,0.5902,0.9218); rgb(22pt)=(0.146029,0.609119,0.907857); rgb(23pt)=(0.138024,0.627629,0.89729); rgb(24pt)=(0.124814,0.645929,0.888343); rgb(25pt)=(0.111252,0.6635,0.876314); rgb(26pt)=(0.0952095,0.679829,0.859781); rgb(27pt)=(0.0688714,0.694771,0.839357); rgb(28pt)=(0.0296667,0.708167,0.816333); rgb(29pt)=(0.00357143,0.720267,0.7917); rgb(30pt)=(0.00665714,0.731214,0.766014); rgb(31pt)=(0.0433286,0.741095,0.73941); rgb(32pt)=(0.0963952,0.75,0.712038); rgb(33pt)=(0.140771,0.7584,0.684157); rgb(34pt)=(0.1717,0.766962,0.655443); rgb(35pt)=(0.193767,0.775767,0.6251); rgb(36pt)=(0.216086,0.7843,0.5923); rgb(37pt)=(0.246957,0.791795,0.556743); rgb(38pt)=(0.290614,0.79729,0.518829); rgb(39pt)=(0.340643,0.8008,0.478857); rgb(40pt)=(0.3909,0.802871,0.435448); rgb(41pt)=(0.445629,0.802419,0.390919); rgb(42pt)=(0.5044,0.7993,0.348); rgb(43pt)=(0.561562,0.794233,0.304481); rgb(44pt)=(0.617395,0.787619,0.261238); rgb(45pt)=(0.671986,0.779271,0.2227); rgb(46pt)=(0.7242,0.769843,0.191029); rgb(47pt)=(0.773833,0.759805,0.16461); rgb(48pt)=(0.820314,0.749814,0.153529); rgb(49pt)=(0.863433,0.7406,0.159633); rgb(50pt)=(0.903543,0.733029,0.177414); rgb(51pt)=(0.939257,0.728786,0.209957); rgb(52pt)=(0.972757,0.729771,0.239443); rgb(53pt)=(0.995648,0.743371,0.237148); rgb(54pt)=(0.996986,0.765857,0.219943); rgb(55pt)=(0.995205,0.789252,0.202762); rgb(56pt)=(0.9892,0.813567,0.188533); rgb(57pt)=(0.978629,0.838629,0.176557); rgb(58pt)=(0.967648,0.8639,0.16429); rgb(59pt)=(0.96101,0.889019,0.153676); rgb(60pt)=(0.959671,0.913457,0.142257); rgb(61pt)=(0.962795,0.937338,0.12651); rgb(62pt)=(0.969114,0.960629,0.106362); rgb(63pt)=(0.9769,0.9839,0.0805)}, mesh/rows=6]
table[point meta=\thisrow{c}] {%
plotData/convergenceTest/Polymesher1/displSurf-1.tsv};

\end{axis}
\end{tikzpicture}%
    \caption{$e_{\displ}$}
  \end{subfigure}
  \caption{Test with manufactured exact solution: error through space and time refinements, for the \texttt{Polymesher1} mesh family.}
  \label{fig:poly1surf}
\end{figure}

\begin{figure}
  \centering
  \begin{subfigure}{.32\linewidth}
    \centering
    \begin{tikzpicture}

  \begin{axis}[%
    width=\linewidth,
    height=\linewidth,
    xmode=log,
    xmin=0.01,
    xmax=1,
    xminorticks=true,
    tick align=outside,
    xlabel style={font=\color{white!15!black}},
    xlabel={$h$},
    ymode=log,
    ymin=0.003125,
    ymax=0.1,
    yminorticks=true,
    ylabel style={font=\color{white!15!black}},
    ylabel={$\tau$},
    zmode=log,
    zmin=1e-05,
    zmax=0.00141706211202513,
    zminorticks=true,
    view={-37.5}{30},
    axis background/.style={fill=white},
    axis x line*=bottom,
    axis y line*=left,
    axis z line*=left,
    legend style={at={(1.03,1)}, anchor=north west, legend cell align=left, align=left, draw=white!15!black}
    ]

\addplot3[%
surf,
shader=flat corner, draw=black, z buffer=sort, colormap={mymap}{[1pt] rgb(0pt)=(0.2422,0.1504,0.6603); rgb(1pt)=(0.25039,0.164995,0.707614); rgb(2pt)=(0.257771,0.181781,0.751138); rgb(3pt)=(0.264729,0.197757,0.795214); rgb(4pt)=(0.270648,0.214676,0.836371); rgb(5pt)=(0.275114,0.234238,0.870986); rgb(6pt)=(0.2783,0.255871,0.899071); rgb(7pt)=(0.280333,0.278233,0.9221); rgb(8pt)=(0.281338,0.300595,0.941376); rgb(9pt)=(0.281014,0.322757,0.957886); rgb(10pt)=(0.279467,0.344671,0.971676); rgb(11pt)=(0.275971,0.366681,0.982905); rgb(12pt)=(0.269914,0.3892,0.9906); rgb(13pt)=(0.260243,0.412329,0.995157); rgb(14pt)=(0.244033,0.435833,0.998833); rgb(15pt)=(0.220643,0.460257,0.997286); rgb(16pt)=(0.196333,0.484719,0.989152); rgb(17pt)=(0.183405,0.507371,0.979795); rgb(18pt)=(0.178643,0.528857,0.968157); rgb(19pt)=(0.176438,0.549905,0.952019); rgb(20pt)=(0.168743,0.570262,0.935871); rgb(21pt)=(0.154,0.5902,0.9218); rgb(22pt)=(0.146029,0.609119,0.907857); rgb(23pt)=(0.138024,0.627629,0.89729); rgb(24pt)=(0.124814,0.645929,0.888343); rgb(25pt)=(0.111252,0.6635,0.876314); rgb(26pt)=(0.0952095,0.679829,0.859781); rgb(27pt)=(0.0688714,0.694771,0.839357); rgb(28pt)=(0.0296667,0.708167,0.816333); rgb(29pt)=(0.00357143,0.720267,0.7917); rgb(30pt)=(0.00665714,0.731214,0.766014); rgb(31pt)=(0.0433286,0.741095,0.73941); rgb(32pt)=(0.0963952,0.75,0.712038); rgb(33pt)=(0.140771,0.7584,0.684157); rgb(34pt)=(0.1717,0.766962,0.655443); rgb(35pt)=(0.193767,0.775767,0.6251); rgb(36pt)=(0.216086,0.7843,0.5923); rgb(37pt)=(0.246957,0.791795,0.556743); rgb(38pt)=(0.290614,0.79729,0.518829); rgb(39pt)=(0.340643,0.8008,0.478857); rgb(40pt)=(0.3909,0.802871,0.435448); rgb(41pt)=(0.445629,0.802419,0.390919); rgb(42pt)=(0.5044,0.7993,0.348); rgb(43pt)=(0.561562,0.794233,0.304481); rgb(44pt)=(0.617395,0.787619,0.261238); rgb(45pt)=(0.671986,0.779271,0.2227); rgb(46pt)=(0.7242,0.769843,0.191029); rgb(47pt)=(0.773833,0.759805,0.16461); rgb(48pt)=(0.820314,0.749814,0.153529); rgb(49pt)=(0.863433,0.7406,0.159633); rgb(50pt)=(0.903543,0.733029,0.177414); rgb(51pt)=(0.939257,0.728786,0.209957); rgb(52pt)=(0.972757,0.729771,0.239443); rgb(53pt)=(0.995648,0.743371,0.237148); rgb(54pt)=(0.996986,0.765857,0.219943); rgb(55pt)=(0.995205,0.789252,0.202762); rgb(56pt)=(0.9892,0.813567,0.188533); rgb(57pt)=(0.978629,0.838629,0.176557); rgb(58pt)=(0.967648,0.8639,0.16429); rgb(59pt)=(0.96101,0.889019,0.153676); rgb(60pt)=(0.959671,0.913457,0.142257); rgb(61pt)=(0.962795,0.937338,0.12651); rgb(62pt)=(0.969114,0.960629,0.106362); rgb(63pt)=(0.9769,0.9839,0.0805)}, mesh/rows=6]
table[point meta=\thisrow{c}] {%
plotData/convergenceTest/Polymesher20/pressureSurf-1.tsv};

  \end{axis}
\end{tikzpicture}%
    \caption{$e_{p}$}
  \end{subfigure}
  \begin{subfigure}{.32\linewidth}
    \centering
    \begin{tikzpicture}

  \begin{axis}[%
    width=\linewidth,
    height=\linewidth,
    xmode=log,
    xmin=0.01,
    xmax=1,
    xminorticks=true,
    tick align=outside,
    xlabel style={font=\color{white!15!black}},
    xlabel={$h$},
    ymode=log,
    ymin=0.003125,
    ymax=0.1,
    yminorticks=true,
    ylabel style={font=\color{white!15!black}},
    ylabel={$\tau$},
    zmode=log,
    zmin=5.18744293214344e-05,
    zmax=0.0018288103006384,
    zminorticks=true,
    view={-37.5}{30},
    axis background/.style={fill=white},
    axis x line*=bottom,
    axis y line*=left,
    axis z line*=left,
    legend style={at={(1.03,1)}, anchor=north west, legend cell align=left, align=left, draw=white!15!black}
    ]

\addplot3[%
surf,
shader=flat corner, draw=black, z buffer=sort, colormap={mymap}{[1pt] rgb(0pt)=(0.2422,0.1504,0.6603); rgb(1pt)=(0.25039,0.164995,0.707614); rgb(2pt)=(0.257771,0.181781,0.751138); rgb(3pt)=(0.264729,0.197757,0.795214); rgb(4pt)=(0.270648,0.214676,0.836371); rgb(5pt)=(0.275114,0.234238,0.870986); rgb(6pt)=(0.2783,0.255871,0.899071); rgb(7pt)=(0.280333,0.278233,0.9221); rgb(8pt)=(0.281338,0.300595,0.941376); rgb(9pt)=(0.281014,0.322757,0.957886); rgb(10pt)=(0.279467,0.344671,0.971676); rgb(11pt)=(0.275971,0.366681,0.982905); rgb(12pt)=(0.269914,0.3892,0.9906); rgb(13pt)=(0.260243,0.412329,0.995157); rgb(14pt)=(0.244033,0.435833,0.998833); rgb(15pt)=(0.220643,0.460257,0.997286); rgb(16pt)=(0.196333,0.484719,0.989152); rgb(17pt)=(0.183405,0.507371,0.979795); rgb(18pt)=(0.178643,0.528857,0.968157); rgb(19pt)=(0.176438,0.549905,0.952019); rgb(20pt)=(0.168743,0.570262,0.935871); rgb(21pt)=(0.154,0.5902,0.9218); rgb(22pt)=(0.146029,0.609119,0.907857); rgb(23pt)=(0.138024,0.627629,0.89729); rgb(24pt)=(0.124814,0.645929,0.888343); rgb(25pt)=(0.111252,0.6635,0.876314); rgb(26pt)=(0.0952095,0.679829,0.859781); rgb(27pt)=(0.0688714,0.694771,0.839357); rgb(28pt)=(0.0296667,0.708167,0.816333); rgb(29pt)=(0.00357143,0.720267,0.7917); rgb(30pt)=(0.00665714,0.731214,0.766014); rgb(31pt)=(0.0433286,0.741095,0.73941); rgb(32pt)=(0.0963952,0.75,0.712038); rgb(33pt)=(0.140771,0.7584,0.684157); rgb(34pt)=(0.1717,0.766962,0.655443); rgb(35pt)=(0.193767,0.775767,0.6251); rgb(36pt)=(0.216086,0.7843,0.5923); rgb(37pt)=(0.246957,0.791795,0.556743); rgb(38pt)=(0.290614,0.79729,0.518829); rgb(39pt)=(0.340643,0.8008,0.478857); rgb(40pt)=(0.3909,0.802871,0.435448); rgb(41pt)=(0.445629,0.802419,0.390919); rgb(42pt)=(0.5044,0.7993,0.348); rgb(43pt)=(0.561562,0.794233,0.304481); rgb(44pt)=(0.617395,0.787619,0.261238); rgb(45pt)=(0.671986,0.779271,0.2227); rgb(46pt)=(0.7242,0.769843,0.191029); rgb(47pt)=(0.773833,0.759805,0.16461); rgb(48pt)=(0.820314,0.749814,0.153529); rgb(49pt)=(0.863433,0.7406,0.159633); rgb(50pt)=(0.903543,0.733029,0.177414); rgb(51pt)=(0.939257,0.728786,0.209957); rgb(52pt)=(0.972757,0.729771,0.239443); rgb(53pt)=(0.995648,0.743371,0.237148); rgb(54pt)=(0.996986,0.765857,0.219943); rgb(55pt)=(0.995205,0.789252,0.202762); rgb(56pt)=(0.9892,0.813567,0.188533); rgb(57pt)=(0.978629,0.838629,0.176557); rgb(58pt)=(0.967648,0.8639,0.16429); rgb(59pt)=(0.96101,0.889019,0.153676); rgb(60pt)=(0.959671,0.913457,0.142257); rgb(61pt)=(0.962795,0.937338,0.12651); rgb(62pt)=(0.969114,0.960629,0.106362); rgb(63pt)=(0.9769,0.9839,0.0805)}, mesh/rows=6]
table[point meta=\thisrow{c}] {%
plotData/convergenceTest/Polymesher20/energySurf-1.tsv};

  \end{axis}
\end{tikzpicture}%
    \caption{$e_{\stress}$}
  \end{subfigure}
  \begin{subfigure}{.32\linewidth}
    \centering \begin{tikzpicture}

  \begin{axis}[%
    width=\linewidth,
    height=\linewidth,
    xmode=log,
    xmin=0.01,
    xmax=1,
    xminorticks=true,
    tick align=outside,
    xlabel style={font=\color{white!15!black}},
    xlabel={$h$},
    ymode=log,
    ymin=0.003125,
    ymax=0.1,
    yminorticks=true,
    ylabel style={font=\color{white!15!black}},
    ylabel={$\tau$},
    zmode=log,
    zmin=1e-06,
    zmax=0.000144686321693441,
    zminorticks=true,
    view={-37.5}{30},
    axis background/.style={fill=white},
    axis x line*=bottom,
    axis y line*=left,
    axis z line*=left,
    legend style={at={(1.03,1)}, anchor=north west, legend cell align=left, align=left, draw=white!15!black}
    ]

\addplot3[%
surf,
shader=flat corner, draw=black, z buffer=sort, colormap={mymap}{[1pt] rgb(0pt)=(0.2422,0.1504,0.6603); rgb(1pt)=(0.25039,0.164995,0.707614); rgb(2pt)=(0.257771,0.181781,0.751138); rgb(3pt)=(0.264729,0.197757,0.795214); rgb(4pt)=(0.270648,0.214676,0.836371); rgb(5pt)=(0.275114,0.234238,0.870986); rgb(6pt)=(0.2783,0.255871,0.899071); rgb(7pt)=(0.280333,0.278233,0.9221); rgb(8pt)=(0.281338,0.300595,0.941376); rgb(9pt)=(0.281014,0.322757,0.957886); rgb(10pt)=(0.279467,0.344671,0.971676); rgb(11pt)=(0.275971,0.366681,0.982905); rgb(12pt)=(0.269914,0.3892,0.9906); rgb(13pt)=(0.260243,0.412329,0.995157); rgb(14pt)=(0.244033,0.435833,0.998833); rgb(15pt)=(0.220643,0.460257,0.997286); rgb(16pt)=(0.196333,0.484719,0.989152); rgb(17pt)=(0.183405,0.507371,0.979795); rgb(18pt)=(0.178643,0.528857,0.968157); rgb(19pt)=(0.176438,0.549905,0.952019); rgb(20pt)=(0.168743,0.570262,0.935871); rgb(21pt)=(0.154,0.5902,0.9218); rgb(22pt)=(0.146029,0.609119,0.907857); rgb(23pt)=(0.138024,0.627629,0.89729); rgb(24pt)=(0.124814,0.645929,0.888343); rgb(25pt)=(0.111252,0.6635,0.876314); rgb(26pt)=(0.0952095,0.679829,0.859781); rgb(27pt)=(0.0688714,0.694771,0.839357); rgb(28pt)=(0.0296667,0.708167,0.816333); rgb(29pt)=(0.00357143,0.720267,0.7917); rgb(30pt)=(0.00665714,0.731214,0.766014); rgb(31pt)=(0.0433286,0.741095,0.73941); rgb(32pt)=(0.0963952,0.75,0.712038); rgb(33pt)=(0.140771,0.7584,0.684157); rgb(34pt)=(0.1717,0.766962,0.655443); rgb(35pt)=(0.193767,0.775767,0.6251); rgb(36pt)=(0.216086,0.7843,0.5923); rgb(37pt)=(0.246957,0.791795,0.556743); rgb(38pt)=(0.290614,0.79729,0.518829); rgb(39pt)=(0.340643,0.8008,0.478857); rgb(40pt)=(0.3909,0.802871,0.435448); rgb(41pt)=(0.445629,0.802419,0.390919); rgb(42pt)=(0.5044,0.7993,0.348); rgb(43pt)=(0.561562,0.794233,0.304481); rgb(44pt)=(0.617395,0.787619,0.261238); rgb(45pt)=(0.671986,0.779271,0.2227); rgb(46pt)=(0.7242,0.769843,0.191029); rgb(47pt)=(0.773833,0.759805,0.16461); rgb(48pt)=(0.820314,0.749814,0.153529); rgb(49pt)=(0.863433,0.7406,0.159633); rgb(50pt)=(0.903543,0.733029,0.177414); rgb(51pt)=(0.939257,0.728786,0.209957); rgb(52pt)=(0.972757,0.729771,0.239443); rgb(53pt)=(0.995648,0.743371,0.237148); rgb(54pt)=(0.996986,0.765857,0.219943); rgb(55pt)=(0.995205,0.789252,0.202762); rgb(56pt)=(0.9892,0.813567,0.188533); rgb(57pt)=(0.978629,0.838629,0.176557); rgb(58pt)=(0.967648,0.8639,0.16429); rgb(59pt)=(0.96101,0.889019,0.153676); rgb(60pt)=(0.959671,0.913457,0.142257); rgb(61pt)=(0.962795,0.937338,0.12651); rgb(62pt)=(0.969114,0.960629,0.106362); rgb(63pt)=(0.9769,0.9839,0.0805)}, mesh/rows=6]
table[point meta=\thisrow{c}] {%
plotData/convergenceTest/Polymesher20/displSurf-1.tsv};

  \end{axis}
\end{tikzpicture}%
    \caption{$e_{\displ}$}
  \end{subfigure}
  \caption{Test with manufactured exact solution: error through space and time refinements, for the \texttt{Polymesher20} mesh family.}
  \label{fig:poly20surf}
\end{figure}
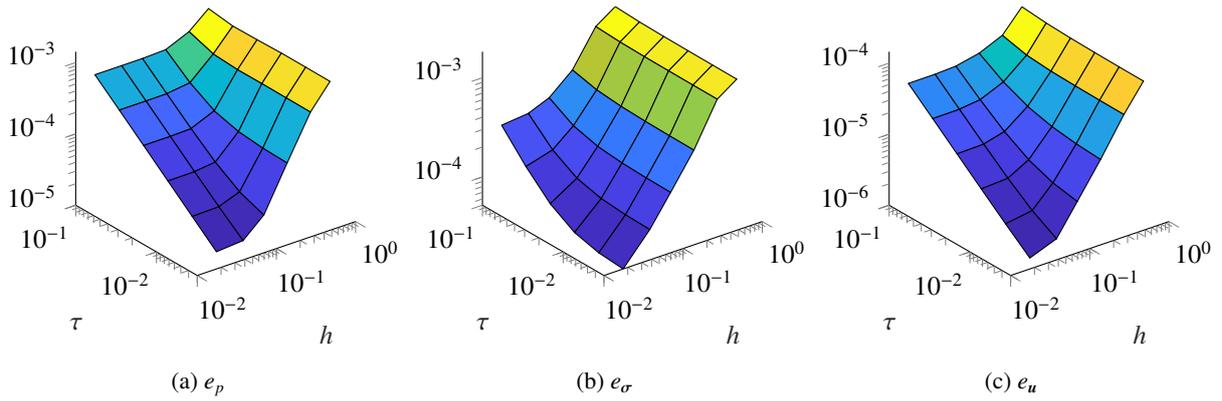

We consider the manufactured regular exact solution proposed in \cite{BofBotDiP16}. 
%
%Again, incompressible solid ($b = 0$) and fluid constituents ($S_{\varepsilon} = 0$) are assumed, with $G = \lambda = \alpha = \mu = 1$ and $\tensorTwo{\kappa} = \begin{pmatrix}1&0\\0&1 \end{pmatrix}$.
%
A body force $\tensorOne{b}$ is introduced in \eqref{eq:momentumBalance} such that the exact distributions of pressure and displacement read 
\begin{align}
  p(x,y,t) &= -\cos(\pi t)\sin(\pi x) \sin(\pi y) \,,
  &
  \boldsymbol{u}(x,y,t) &= \sin(\pi t)
                          \begin{pmatrix}
                            -\cos(\pi x)\cos(\pi y) \\ sin(\pi
                            x)\sin(\pi y)
                          \end{pmatrix}\,.
\end{align}
The following error measures are considered:
\begin{align}
  e_p &= \left(\int_0^T\norm[\lebl{\Omega}]
        {p(t)- p_h(t)}^2 \right)^{\frac12} \,,
  &
    e_{\stress^{\prime}}
  &= \left( \int_0^T \sum_{K\in\Th}
    \norm[K]{
    \tilde{\stress}^\prime_K\left(\II{\displ(t)}\right)
    - \tilde{\stress}^\prime_K\left(\displ_h(t)\right)}^2
    \right)^{\frac12} \,,
  \\
  e_{\displ} &= \left( \int_0^T \sum_{K\in\Th}
               \norm[K]
               {\tilde{\displ}_K(t) -
               \tilde{\displ}_{h,K}(t)}^2 \right)^{\frac12}\,.
\end{align}
The above quantities are computed on the families of meshes illustrated
in Fig.~\ref{fig:mesh_illustration}.
These meshes are refined together
with the time discretization parameter $\Delta t$, in such a way that,
if $(N^{\mathrm{cells}}_{k}, \Delta t_k)$ is the pair of
discretization parameters at the $k$-th refinement (number of cells
and time discretization parameter, respectively), then
$(N^{\mathrm{cells}}_{k+1},\Delta t_{k+1}) = (4N^{\mathrm{cells}}_k,
\Delta t_k/2)$. In Fig.~\ref{fig:totalErrors} we display the
behavior of the above quantities with respect to the maximum diameter
of the discretization $h$, at each refinement step, reporting the
experimental convergence rates in the legend, computed using the last
two points of each line. We can see that the rates of convergence are
in very good accordance with the theory on the more regular meshes
(\texttt{Polymesher20} and \texttt{Cartesian}), while on less regular
meshes we still see a preasymptotic behavior, due to the fact that
$h$ does not scale exactly as
$\frac{1}{\sqrt{N^{\mathrm{cells}}}}$. To check the correct
convergence rates in time, in Figs.~\ref{fig:pressureDifferentDiams}
and \ref{fig:displacementDifferentDiams} we display the behavior of
$e_p$ and $e_{\displ}$ for fixed $h$, as functions of the time
discretization parameter $\Delta t$. We can see that, once the time
component of the error becomes dominant, we obtain the expected rates
of convergence. Finally, Figs.~\ref{fig:poly20surf} to
\ref{fig:skewedsurf} display the errors computed for each choice of
$\Delta t$ and $h$ and confirm the behaviors expected from
theoretical results. Indeed, it can be checked that all quantities are
asymptotically decreasing linearly with respect to $h$; $e_p$
and $e_{\displ}$ decrease linearly with respect to $\Delta t$, and
$e_{\stress}$ is asymptotically constant with respect to $\Delta t$.

%%% Local Variables:
%%% mode: latex
%%% TeX-master: "main"
%%% End:

\subsection{Cantilevered square block problem}
\label{sec:cantilever}

\begin{table}
\small
\hfill
\begin{tabular}{crrrr}
\toprule
level &
\multicolumn{1}{c}{$|\mathcal{V}|$} & %\multicolumn{1}{r}{\begin{tabular}{@{}c@{} } number of  \\ vertices \end{tabular}}&
\multicolumn{1}{c}{$|\mathcal{T}|$} & %\multicolumn{1}{r}{\begin{tabular}{@{}c@{} } number of  \\ cells \end{tabular}}&
\multicolumn{1}{c}{$|\mathcal{F}|$} & %\multicolumn{1}{r}{\begin{tabular}{@{}c@{} } number of  \\ faces \end{tabular}}&
\multicolumn{1}{r}{\begin{tabular}{@{}c@{} } number of  \\ unknowns \end{tabular}} \\ 
\midrule
 0 &     121 &     100 &     220 &     562 \\ 
 1 &     441 &     400 &     840 &   2,122 \\ 
 2 &   1,681 &   1,600 &   3,280 &   8,242 \\ 
 3 &   6,561 &   6,400 &  12,960 &  32,482 \\ 
 4 &  25,921 &  25,600 &  51,520 & 128,962 \\ 
 5 & 103,041 & 102,400 & 205,440 & 513,922 \\ 
\bottomrule
\\
\multicolumn{5}{c}{(a) \texttt{Cartesian} and \texttt{Skewed}.}
\end{tabular} 
\hfill
\begin{tabular}{crrrr}
\toprule
level &
\multicolumn{1}{c}{$|\mathcal{V}|$} & %\multicolumn{1}{r}{\begin{tabular}{@{}c@{} } number of  \\ vertices \end{tabular}}&
\multicolumn{1}{c}{$|\mathcal{T}|$} & %\multicolumn{1}{r}{\begin{tabular}{@{}c@{} } number of  \\ cells \end{tabular}}&
\multicolumn{1}{c}{$|\mathcal{F}|$} & %\multicolumn{1}{r}{\begin{tabular}{@{}c@{} } number of  \\ faces \end{tabular}}&
\multicolumn{1}{r}{\begin{tabular}{@{}c@{} } number of  \\ unknowns \end{tabular}} \\ 
\midrule
 0 &     118 &     110 &     227 &     573 \\ 
 1 &     382 &     391 &     772 &   1,927 \\ 
 2 &   1,500 &   1,608 &   3,107 &   7,715 \\ 
 3 &   5,847 &   6,399 &  12,245 &  30,338 \\ 
 4 &  22,983 &  25,543 &  48,525 & 120,034 \\ 
 5 &  91,883 & 102,396 & 194,278 & 480,440 \\ 
\bottomrule
\\
\multicolumn{5}{c}{(b) \texttt{Hybrid}.}
\end{tabular} 
\hfill\null

\bigskip

\hfill
\begin{tabular}{crrrr}
\toprule
level &
\multicolumn{1}{c}{$|\mathcal{V}|$} & %\multicolumn{1}{r}{\begin{tabular}{@{}c@{} } number of  \\ vertices \end{tabular}}&
\multicolumn{1}{c}{$|\mathcal{T}|$} & %\multicolumn{1}{r}{\begin{tabular}{@{}c@{} } number of  \\ cells \end{tabular}}&
\multicolumn{1}{c}{$|\mathcal{F}|$} & %\multicolumn{1}{r}{\begin{tabular}{@{}c@{} } number of  \\ faces \end{tabular}}&
\multicolumn{1}{r}{\begin{tabular}{@{}c@{} } number of  \\ unknowns \end{tabular}} \\ 
\midrule
 0 &     202 &     100 &     301 &     805 \\ 
 1 &     802 &     400 &   1,201 &   3,205 \\ 
 2 &   3,201 &   1,600 &   4,800 &  12,802 \\ 
 3 &  12,802 &   6,400 &  19,201 &  51,205 \\ 
 4 &  51,194 &  25,600 &  76,793 & 204,781 \\ 
 5 & 204,732 & 102,400 & 307,131 & 818,995 \\ 
\bottomrule
\\
\multicolumn{5}{c}{(c) \texttt{Polymesher1}.}
\end{tabular} 
\hfill
\begin{tabular}{crrrr}
\toprule
level &
\multicolumn{1}{c}{$|\mathcal{V}|$} & %\multicolumn{1}{r}{\begin{tabular}{@{}c@{} } number of  \\ vertices \end{tabular}}&
\multicolumn{1}{c}{$|\mathcal{T}|$} & %\multicolumn{1}{r}{\begin{tabular}{@{}c@{} } number of  \\ cells \end{tabular}}&
\multicolumn{1}{c}{$|\mathcal{F}|$} & %\multicolumn{1}{r}{\begin{tabular}{@{}c@{} } number of  \\ faces \end{tabular}}&
\multicolumn{1}{r}{\begin{tabular}{@{}c@{} } number of  \\ unknowns \end{tabular}} \\ 
\midrule
 0 &     202 &     100 &     301 &     805 \\ 
 1 &     802 &     400 &   1,201 &   3,205 \\ 
 2 &   3,202 &   1,600 &   4,801 &  12,805 \\ 
 3 &  12,802 &   6,400 &  19,201 &  51,205 \\ 
 4 &  51,200 &  25,600 &  76,799 & 204,799 \\ 
 5 & 204,783 & 102,400 & 307,182 & 819,148 \\ 
\bottomrule
\\
\multicolumn{5}{c}{(d) \texttt{Polymesher20}.}
\end{tabular} 
\hfill\null
\caption{Cantilevered square block: grid refinement and problem size.}
\label{tab:grids}
\end{table}

We consider the cantilever problem discussed in \cite{phillips2009overcoming} and later used in \cite{frigo2020efficient} to demonstrate the robustness and computational efficiency of coupled schemes relying on the unstructured macro-element stabilization.
A sketch of the setup of the problem is shown in Fig.~\ref{fig:cantilever_sketch}.
We assume the problem domain $\Omega$ to be the unit square.
The mechanical boundary conditions fix the displacement ($\tensorOne{u} = \tensorOne{0}$) on the left, and impose a unit downward traction at the top and a zero traction on the right and at the bottom.
We impose no-flow boundary conditions on the four sides.
The problem parameters are identical to those used in \cite{frigo2020efficient} and are provided in Table \ref{tab:parameters}.
Sequences of refined meshes for each family are employed for a total of 6 levels (Table \ref{tab:grids}), with a problem size ranging three orders of magnitude in terms of global number of unknowns.

\begin{figure}  
    \small
  
    \hfill
    \begin{subfigure}{.18\linewidth}
      \centerline{\includegraphics[width=\linewidth]{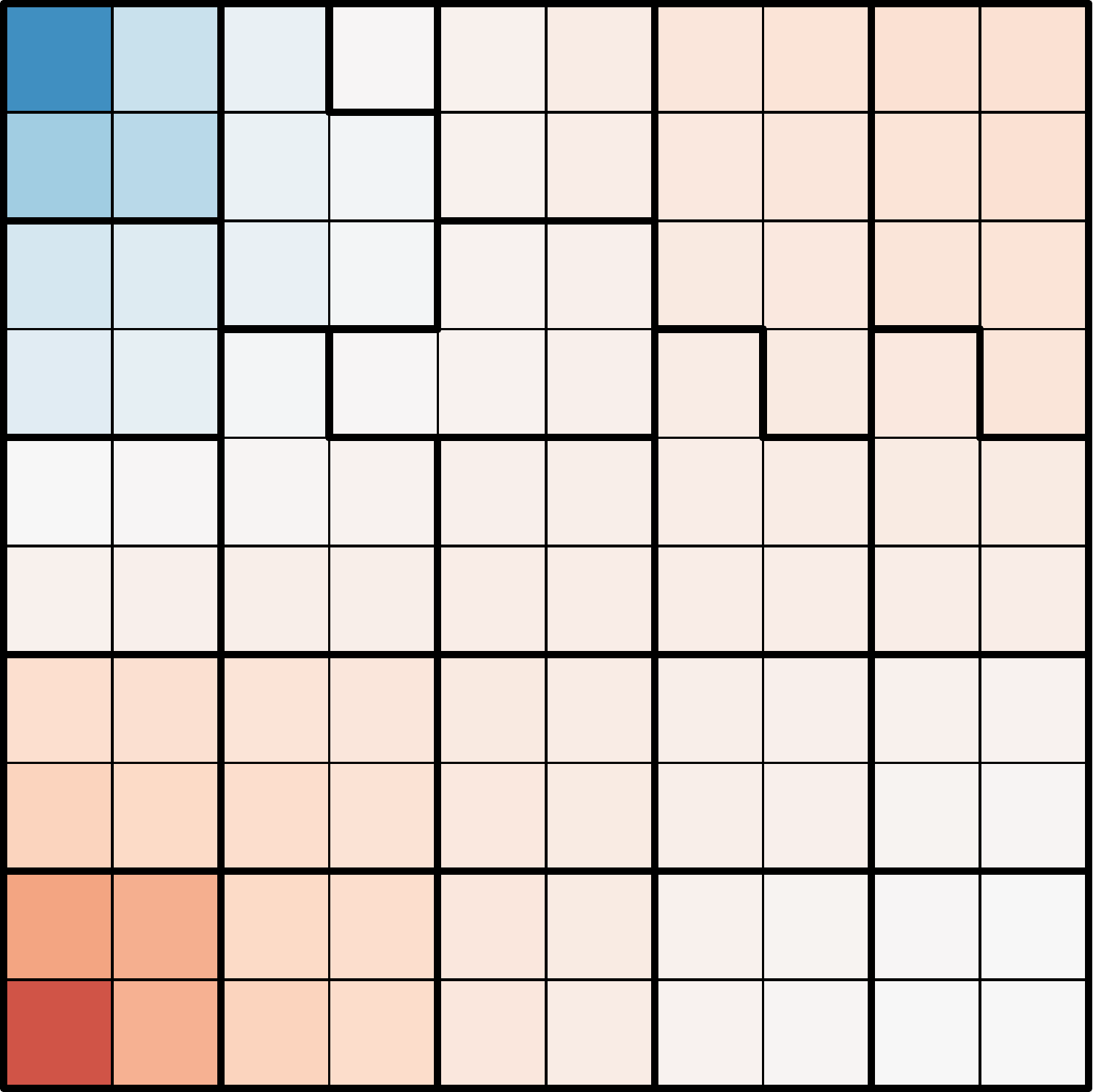}}
      \caption{{\footnotesize{\texttt{Cartesian}}}\textsuperscript{($\star$)}}
    \end{subfigure}
    \hfill    
    \begin{subfigure}{.18\linewidth}
      \centerline{\includegraphics[width=\linewidth]{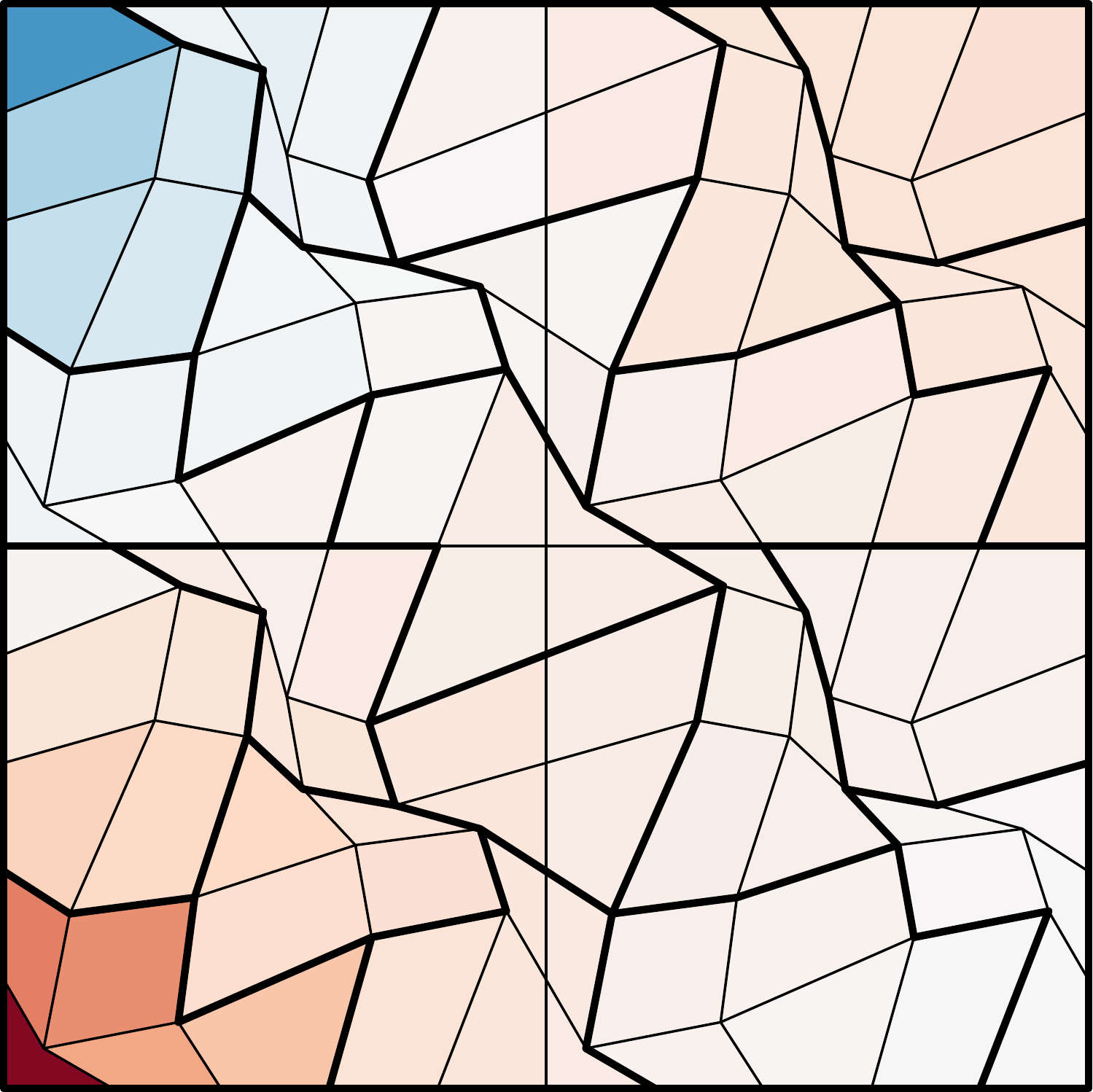}}
      \caption{{\footnotesize{\texttt{Skewed}}}\textsuperscript{($\star$)}}
    \end{subfigure}  
    \hfill    
    \begin{subfigure}{.18\linewidth}
      \centerline{\includegraphics[width=\linewidth]{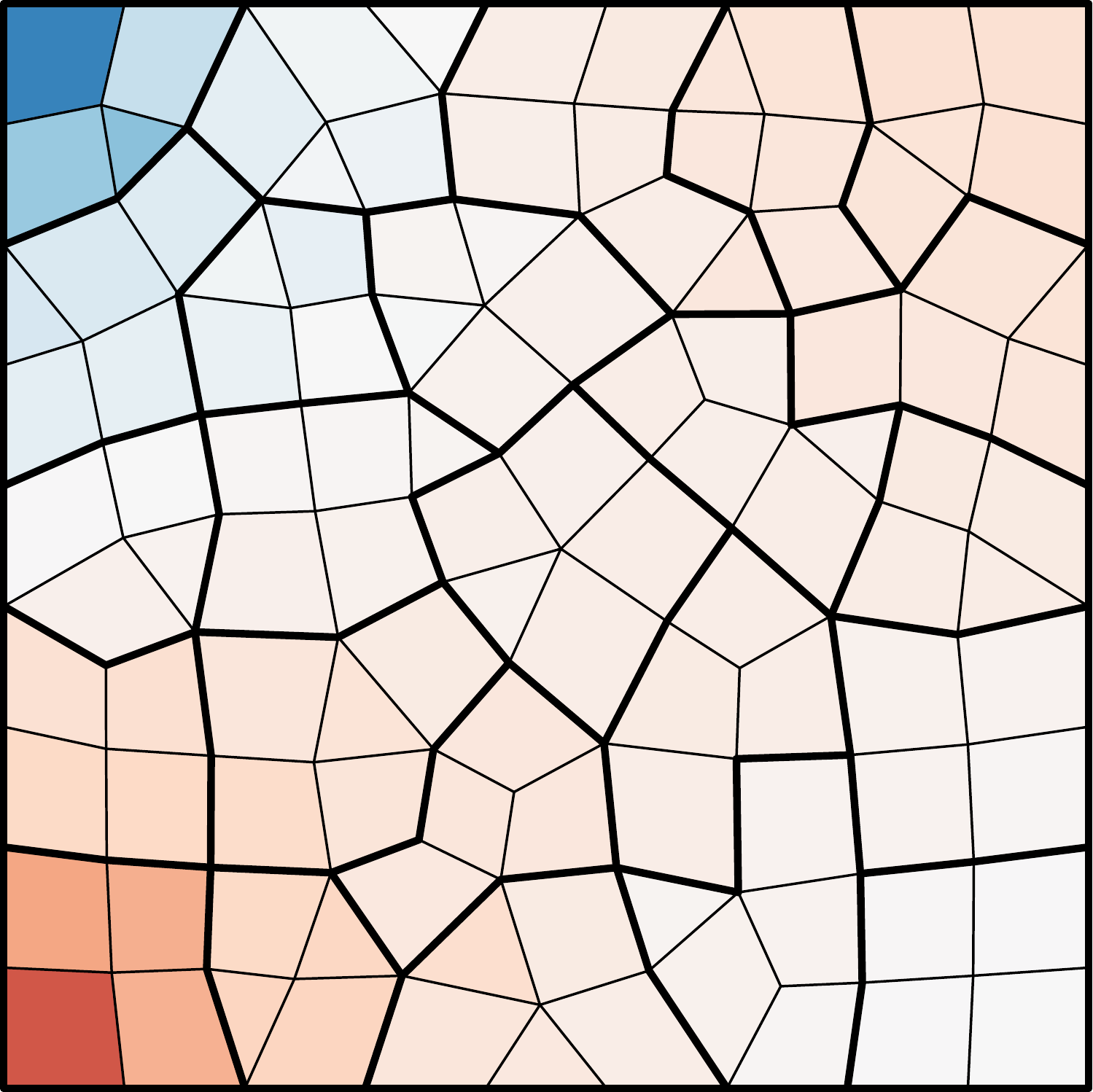}}
      \caption{{\footnotesize{\texttt{Hybrid}}}\textsuperscript{($\star$)}}
    \end{subfigure}
    \hfill    
    \begin{subfigure}{.18\linewidth}
      \centerline{\includegraphics[width=\linewidth]{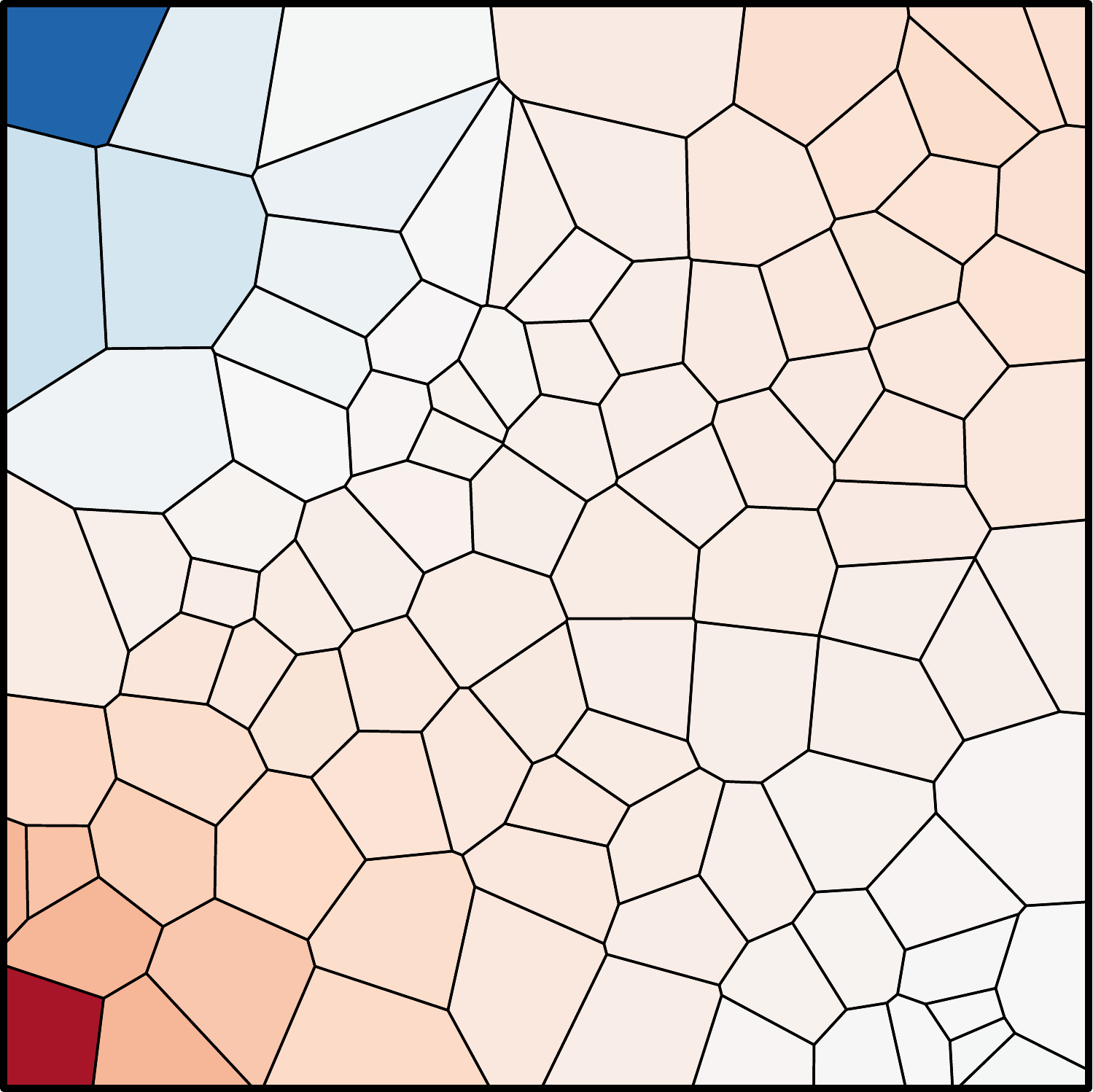}}
      \caption{{\footnotesize{\texttt{Polymesher1}}}}
    \end{subfigure}
    \hfill    
    \begin{subfigure}{.18\linewidth}
      \centerline{\includegraphics[width=\linewidth]{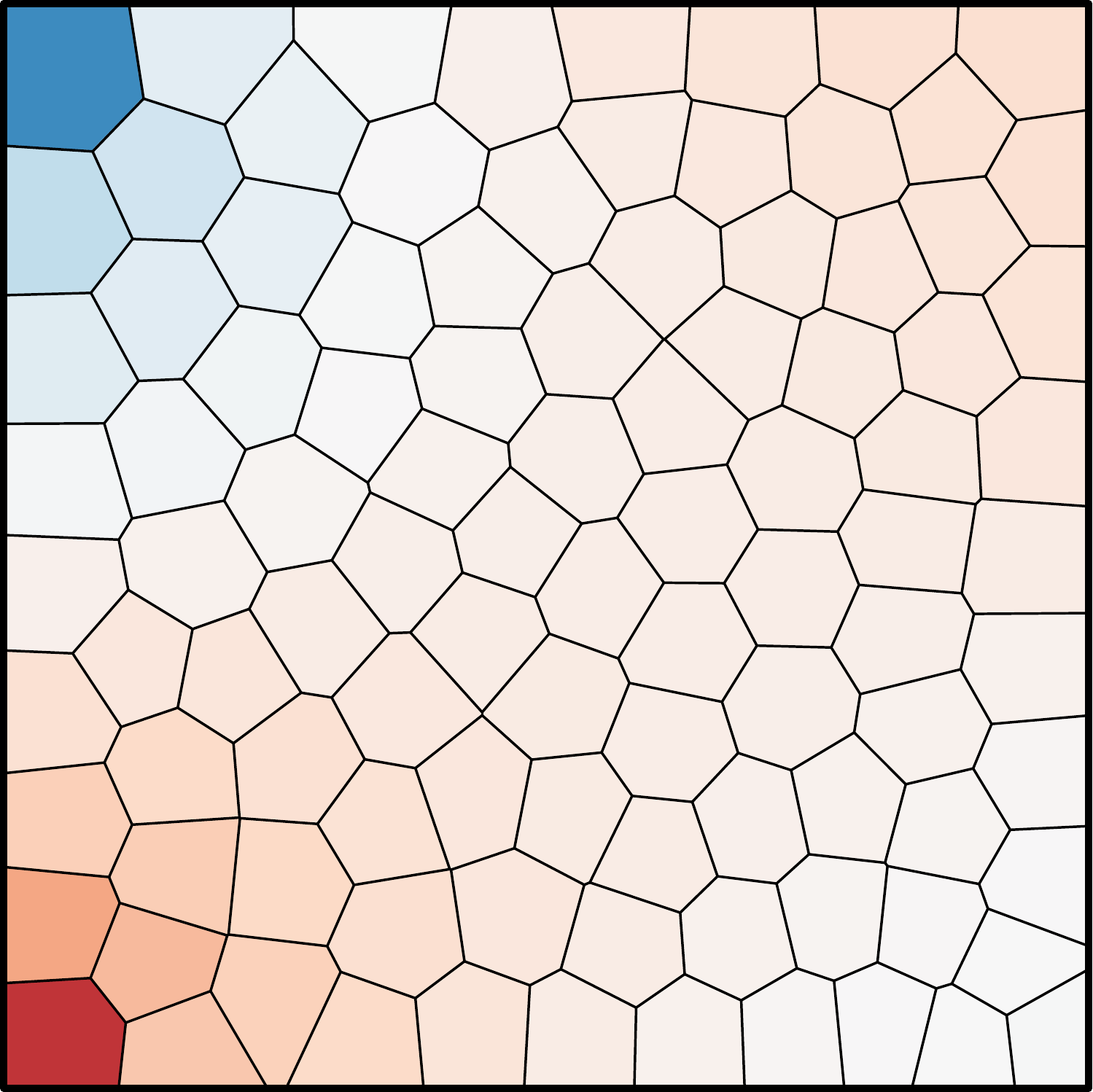}}
      \caption{{\footnotesize{\texttt{Polymesher20}}}}
    \end{subfigure}
    \hfill\null
  
    \bigskip
     
    \hfill
    \begin{subfigure}[t]{.18\linewidth}
      \centerline{\includegraphics[width=\linewidth]{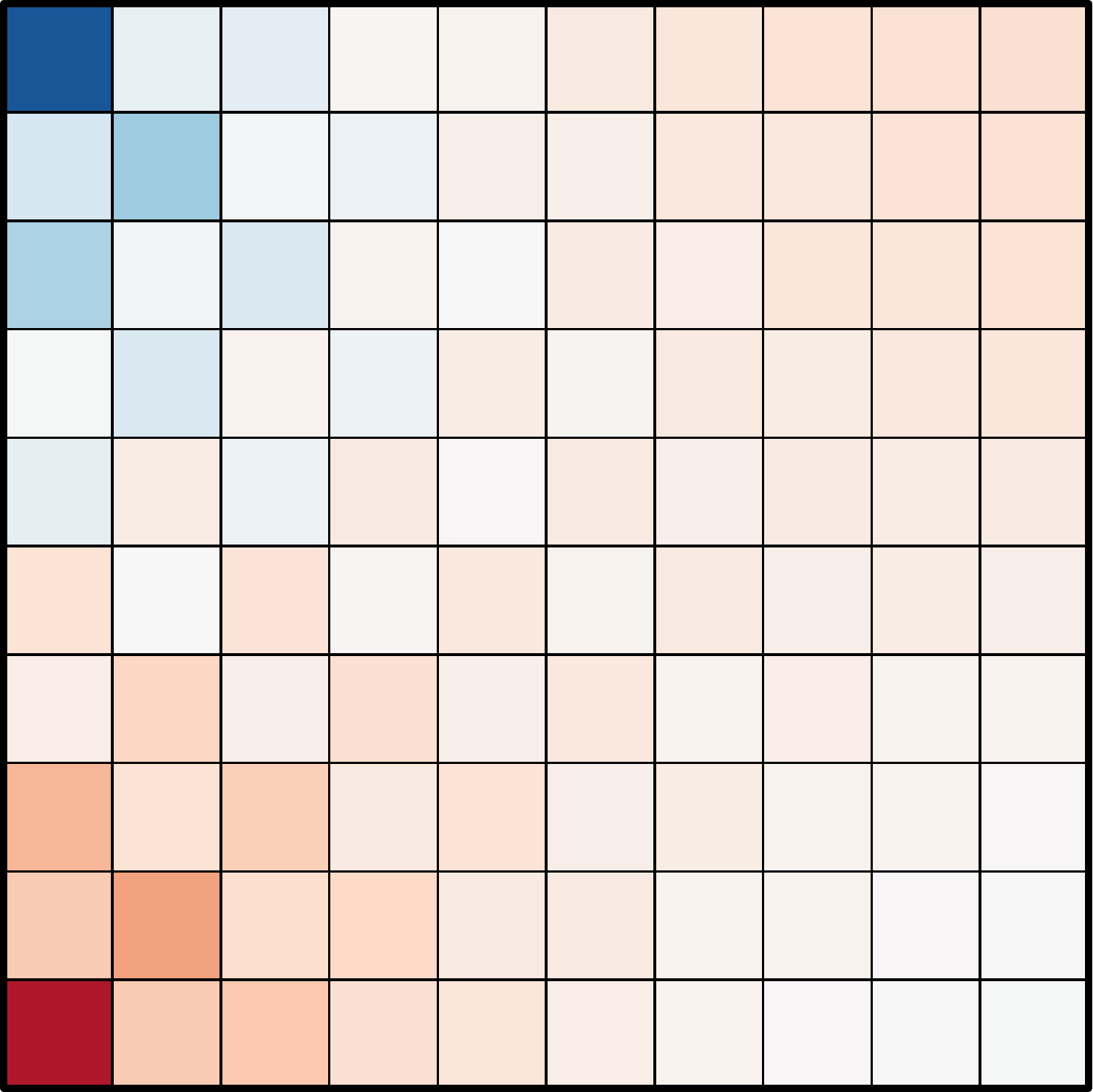}}
      \caption{{\footnotesize{\texttt{Cartesian}}}}
    \end{subfigure}
    \hfill    
    \begin{subfigure}[t]{.18\linewidth}
      \centerline{\includegraphics[width=\linewidth]{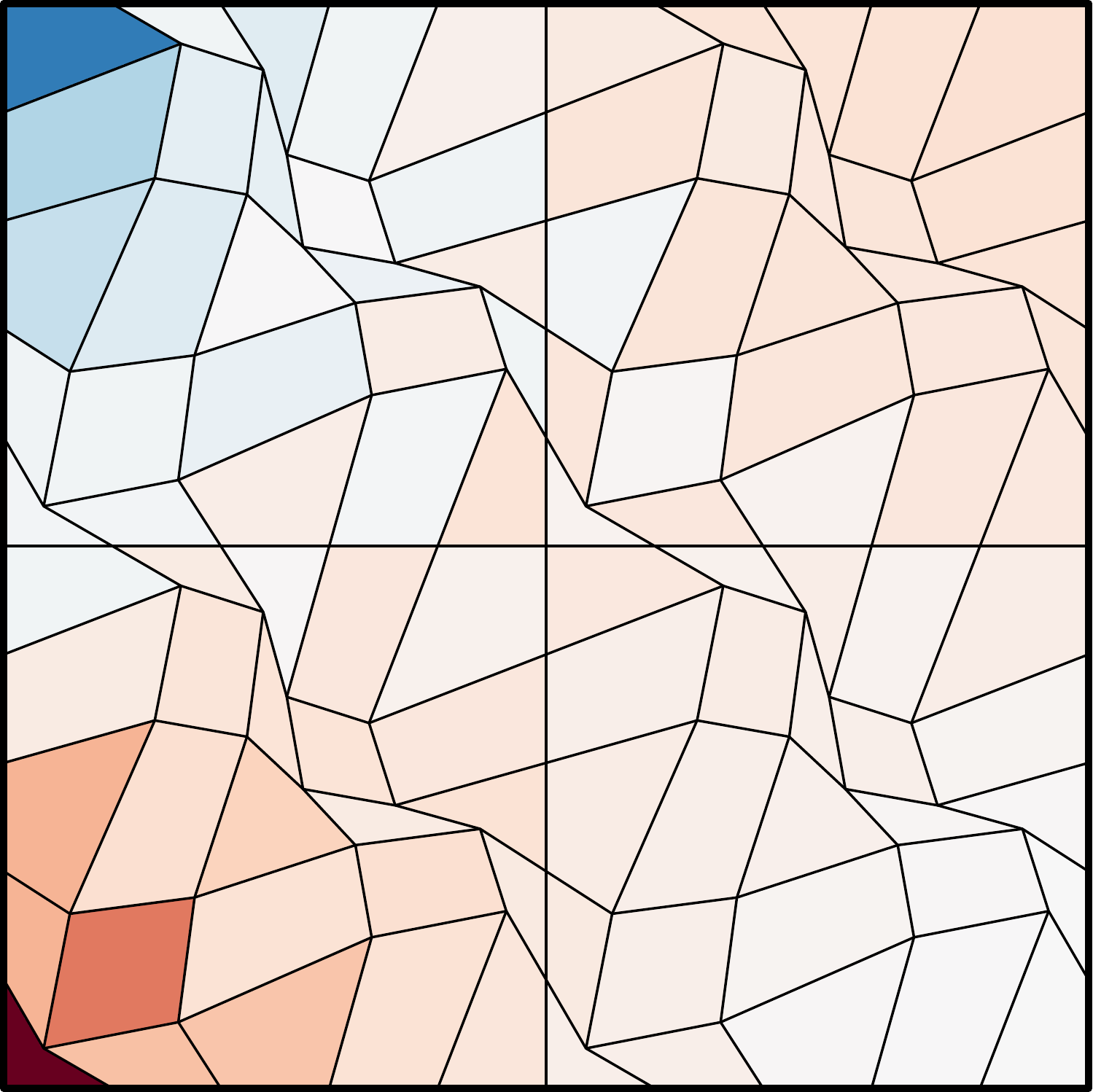}}
      \caption{{\footnotesize{\texttt{Skewed}}}}
    \end{subfigure}  
    \hfill    
    \begin{subfigure}[t]{.18\linewidth}
      \centerline{\includegraphics[width=\linewidth]{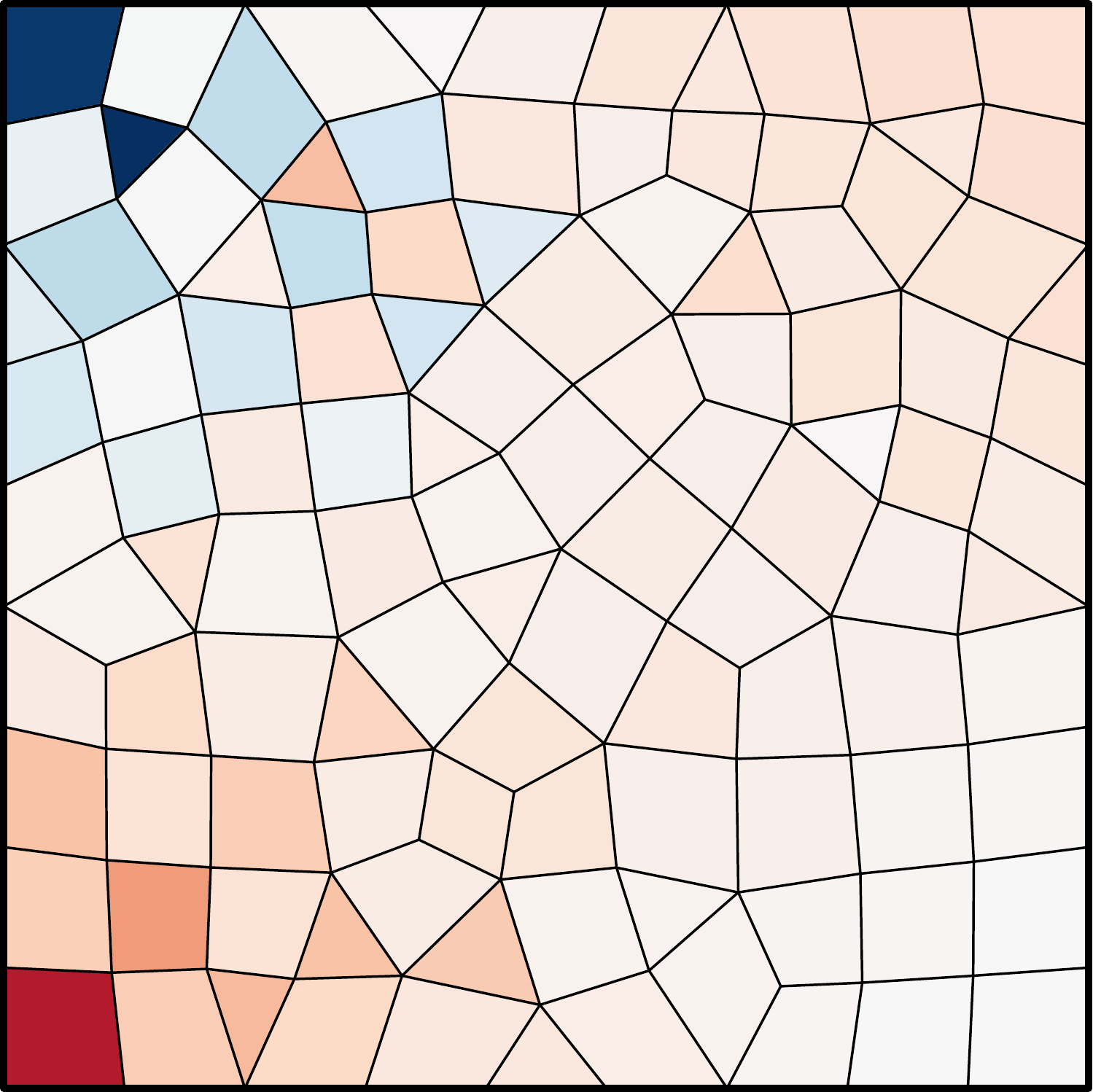}}
      \caption{{\footnotesize{\texttt{Hybrid}}}}
    \end{subfigure}
    \hfill      
    \begin{subfigure}[t]{.18\linewidth}
      \flushright
      \begin{tikzpicture}[scale=1.0]
        \node[anchor=south west,inner sep=0] (image) at (0,0) {\includegraphics[height=.8\linewidth]{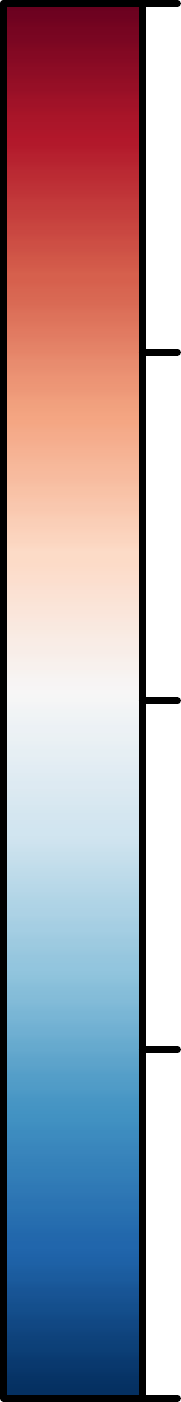}};
        \node [left] at (.5\linewidth,0) {$-3.200$};
        \node [left] at (.5\linewidth,0.2\linewidth) {$-1.600$};
        \node [left] at (.5\linewidth,0.4\linewidth) {$0.000$};
        \node [left] at (.5\linewidth,0.6\linewidth) {$1.600$};
        \node [left] at (.5\linewidth,0.8\linewidth) {$3.200$};
        \node [right] at (0,.9\linewidth) {Pressure [Pa]};
      \end{tikzpicture} 
    \end{subfigure}
    \hfill
    \begin{subfigure}[t]{.18\linewidth}
      \hspace{\linewidth}
    \end{subfigure}  
    \hfill\null  
  
  \caption{Cantilevered square block problem: pressure solution after a single time step $\Delta t= 1 \times 10^{-5}$ s for level 0 mesh of each family (see Table \ref{tab:grids}).   In panels a, b, and c, the star symbol ($\star$) superscript indicates that the local pressure-jump stabilization was introduced, with the unstructured macro-element mesh highlighted using thicker edges. The remaining panels (d to h) were obtained with the unstabilized formulation.}
  \label{fig:cartilever_cart_skew_hybr_dt_1e_5}
\end{figure}

\begin{figure}  
    \small
  
    \hfill
    \begin{subfigure}{.18\linewidth}
      \centerline{\includegraphics[width=\linewidth]{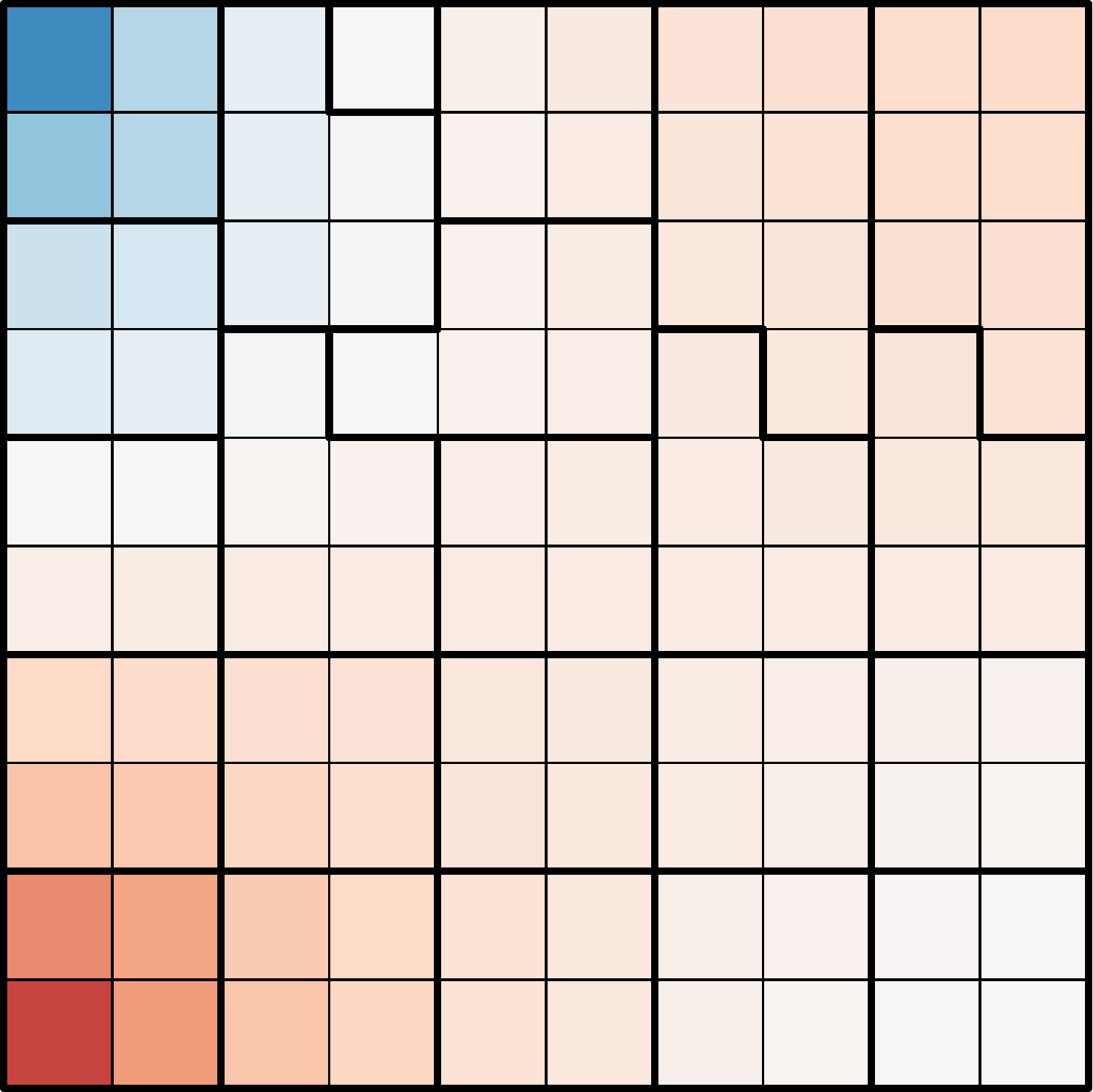}}
      \caption{{\footnotesize{\texttt{Cartesian}}}\textsuperscript{($\star$)}}
    \end{subfigure}
    \hfill    
    \begin{subfigure}{.18\linewidth}
      \centerline{\includegraphics[width=\linewidth]{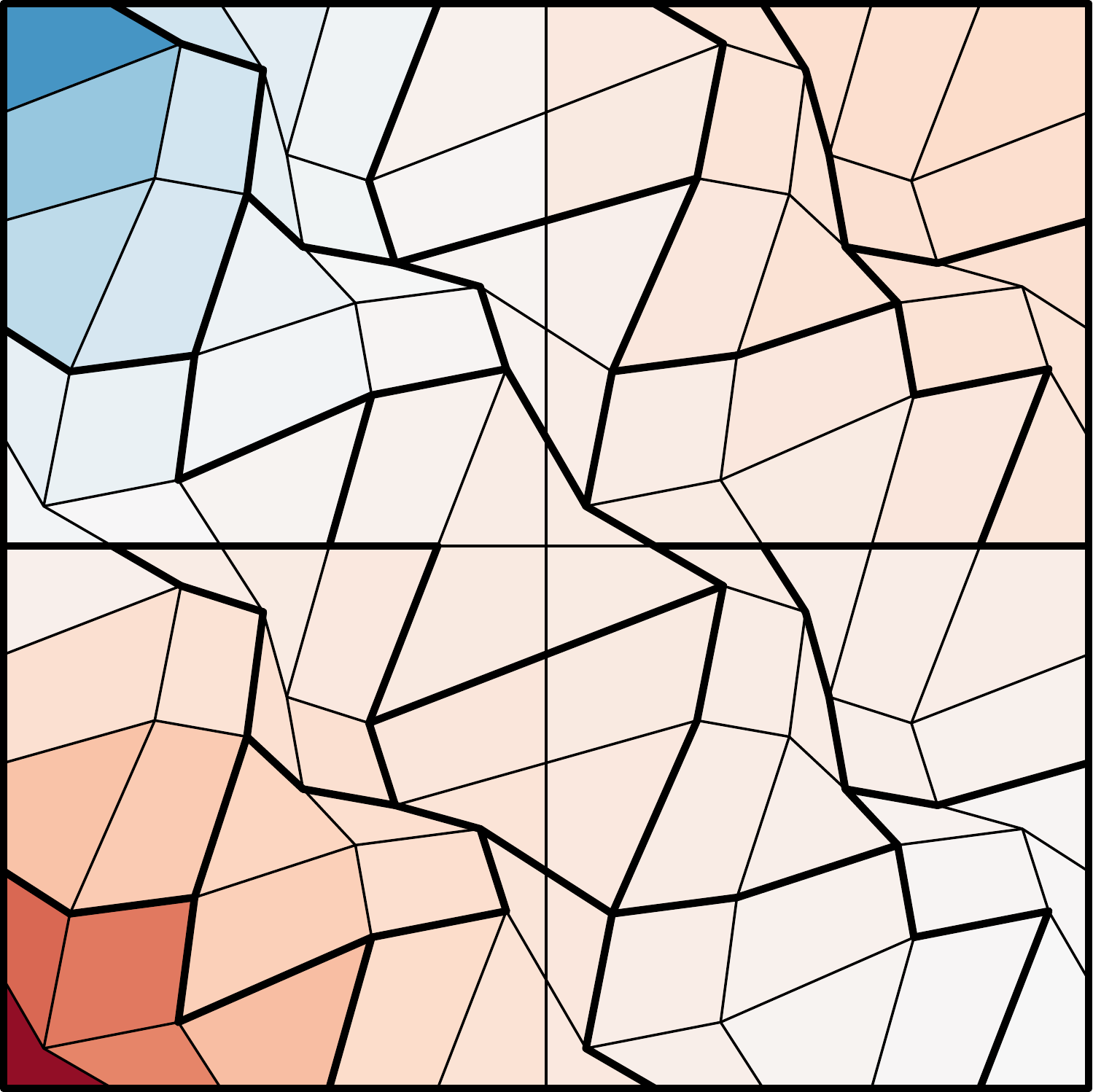}}
      \caption{{\footnotesize{\texttt{Skewed}}}\textsuperscript{($\star$)}}
    \end{subfigure}  
    \hfill    
    \begin{subfigure}{.18\linewidth}
      \centerline{\includegraphics[width=\linewidth]{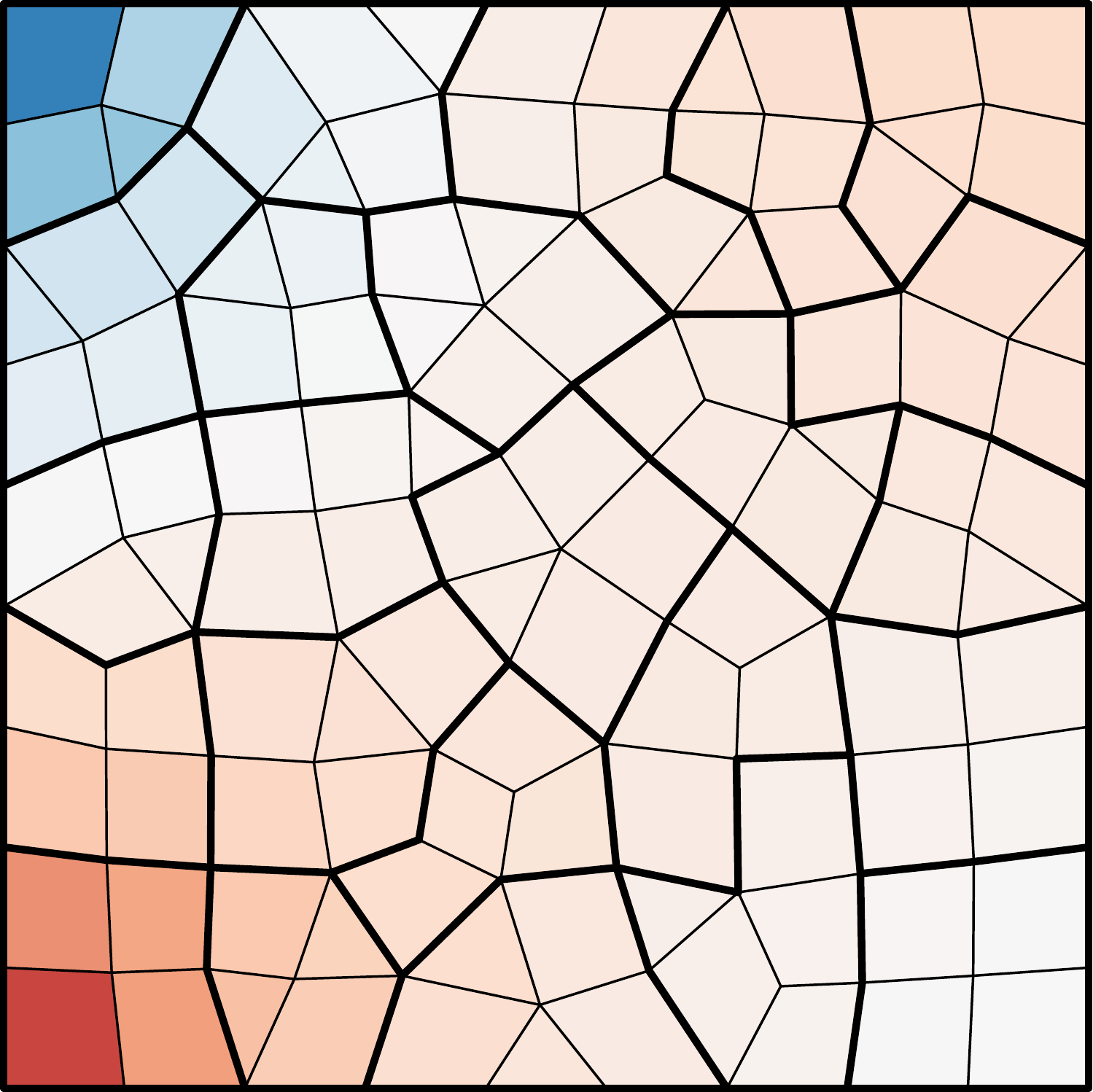}}
      \caption{{\footnotesize{\texttt{Hybrid}}}\textsuperscript{($\star$)}}
    \end{subfigure}
    \hfill    
    \begin{subfigure}{.18\linewidth}
      \centerline{\includegraphics[width=\linewidth]{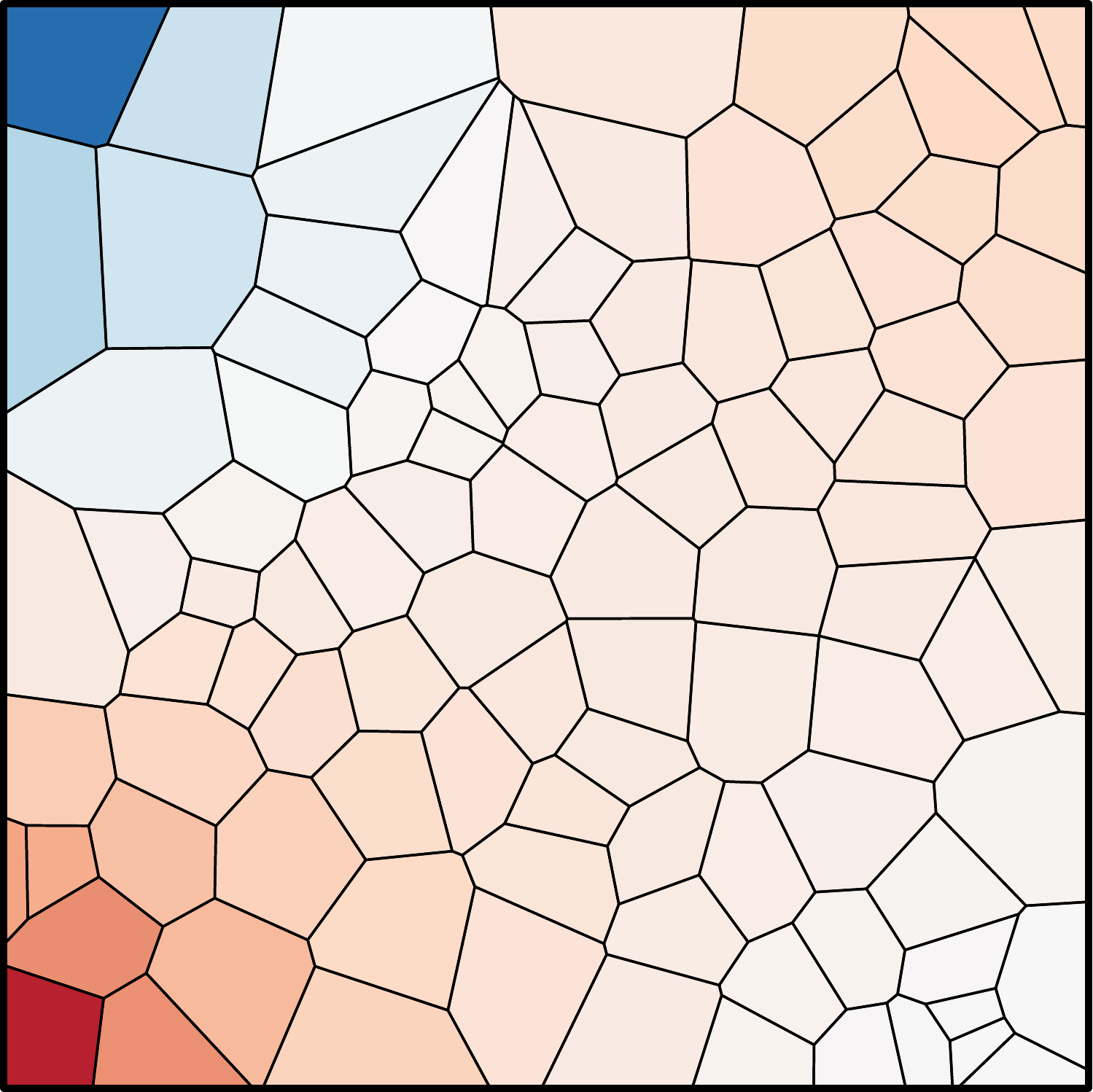}}
      \caption{{\footnotesize{\texttt{Polymesher1}}}}
    \end{subfigure}
    \hfill    
    \begin{subfigure}{.18\linewidth}
      \centerline{\includegraphics[width=\linewidth]{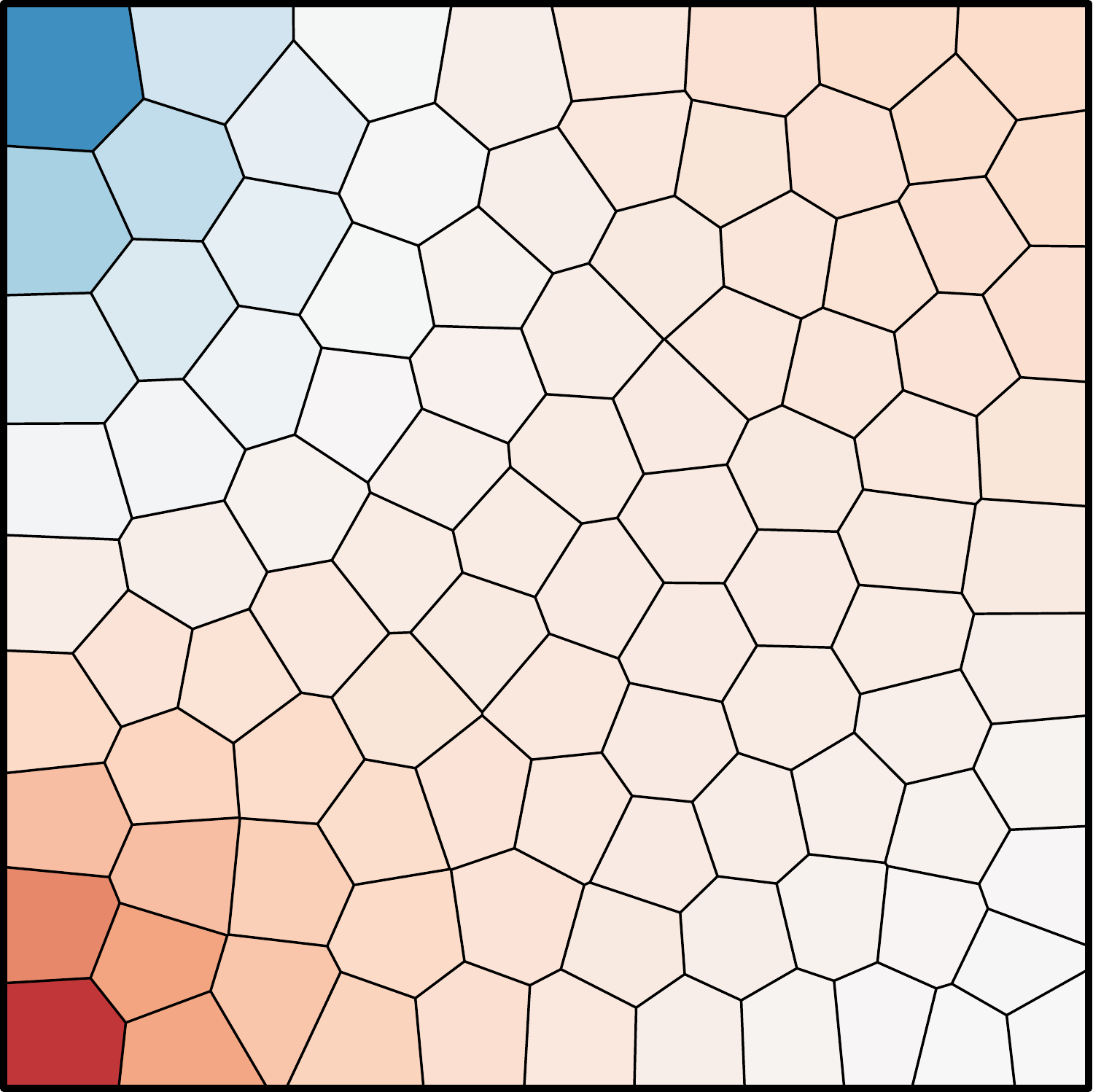}}
      \caption{{\footnotesize{\texttt{Polymesher20}}}}
    \end{subfigure}
    \hfill\null
  
    \bigskip
     
    \hfill
    \begin{subfigure}[t]{.18\linewidth}
      \centerline{\includegraphics[width=\linewidth]{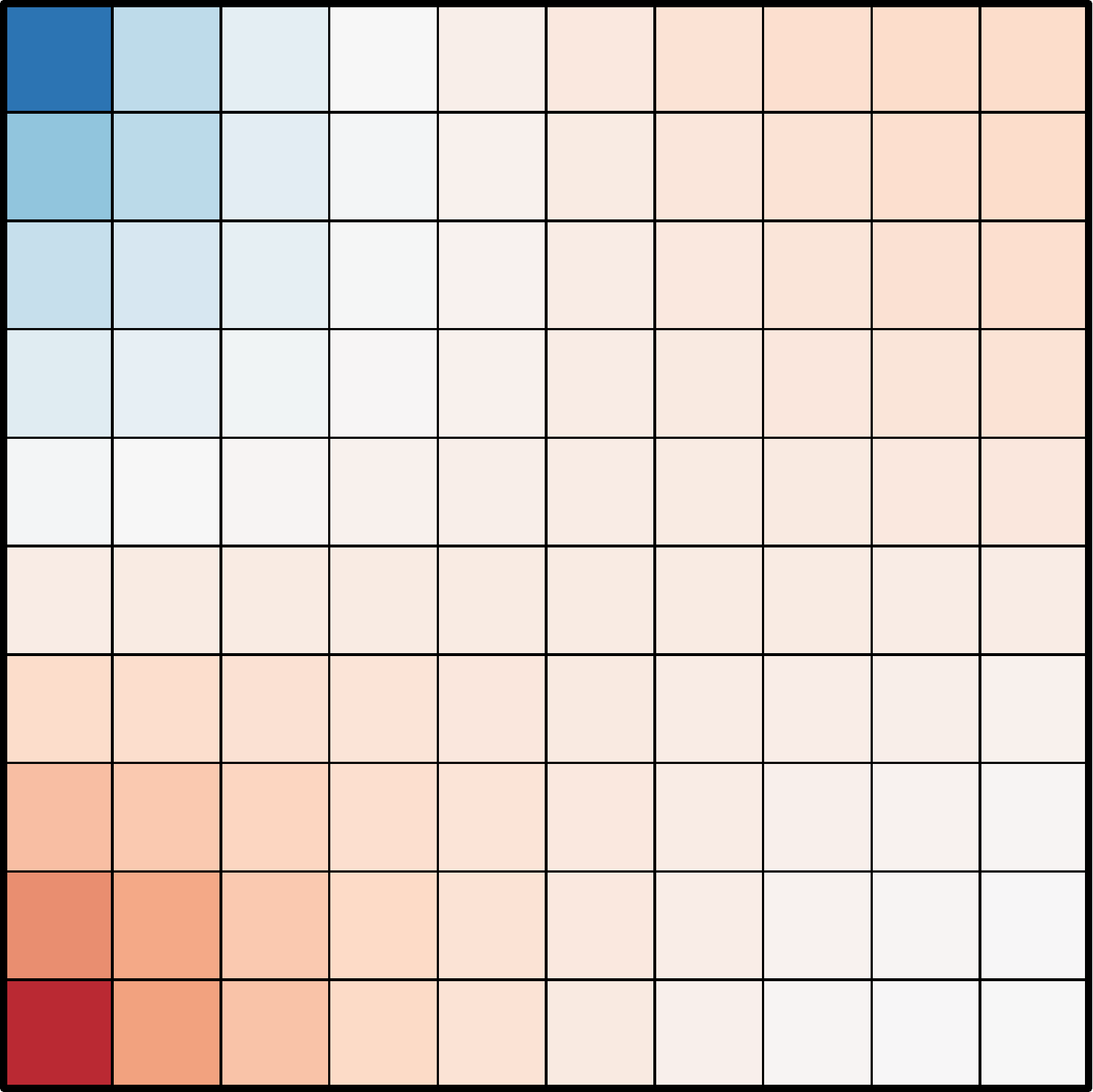}}
      \caption{{\footnotesize{\texttt{Cartesian}}}}
    \end{subfigure}
    \hfill    
    \begin{subfigure}[t]{.18\linewidth}
      \centerline{\includegraphics[width=\linewidth]{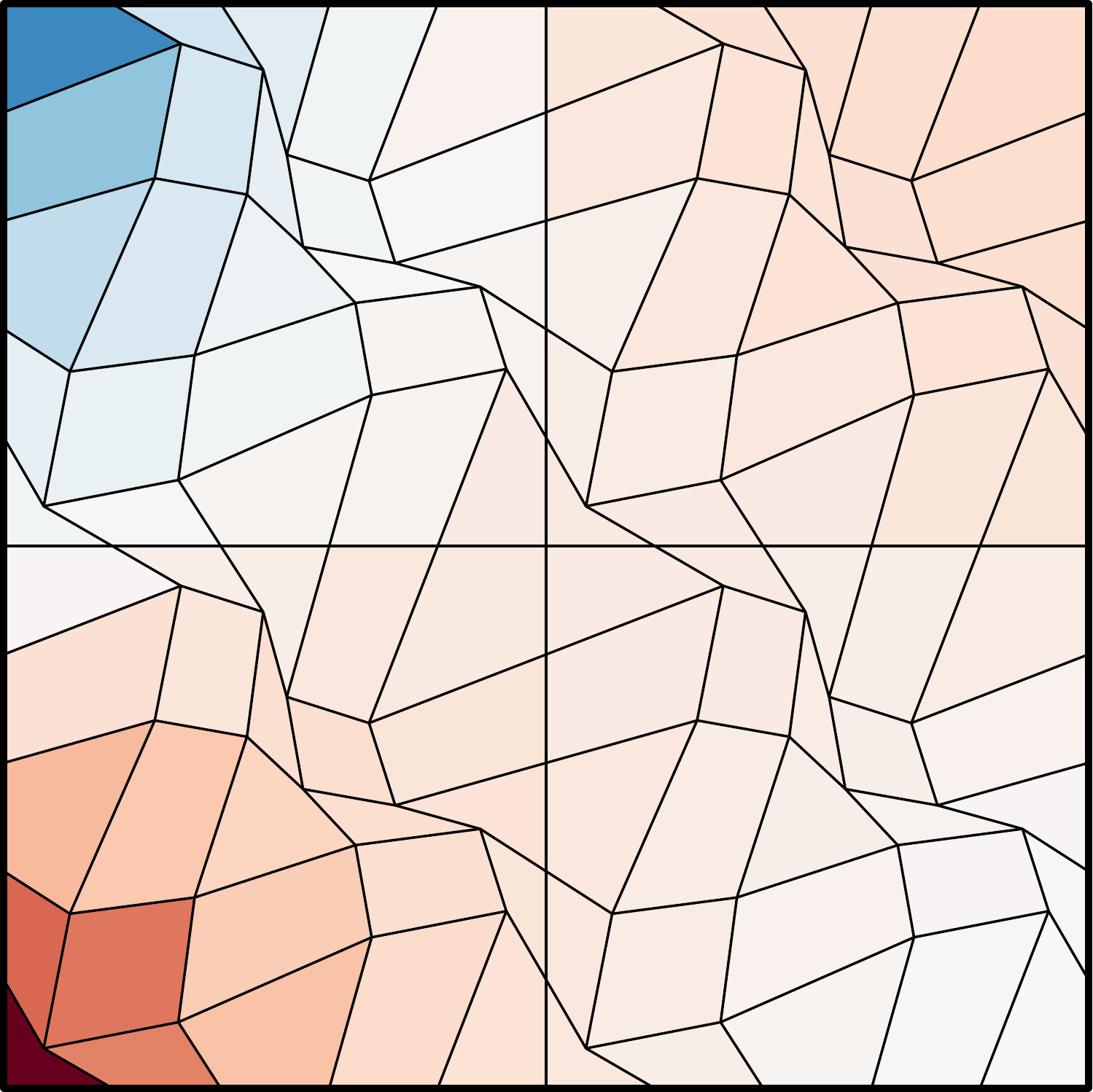}}
      \caption{{\footnotesize{\texttt{Skewed}}}}
    \end{subfigure}  
    \hfill    
    \begin{subfigure}[t]{.18\linewidth}
      \centerline{\includegraphics[width=\linewidth]{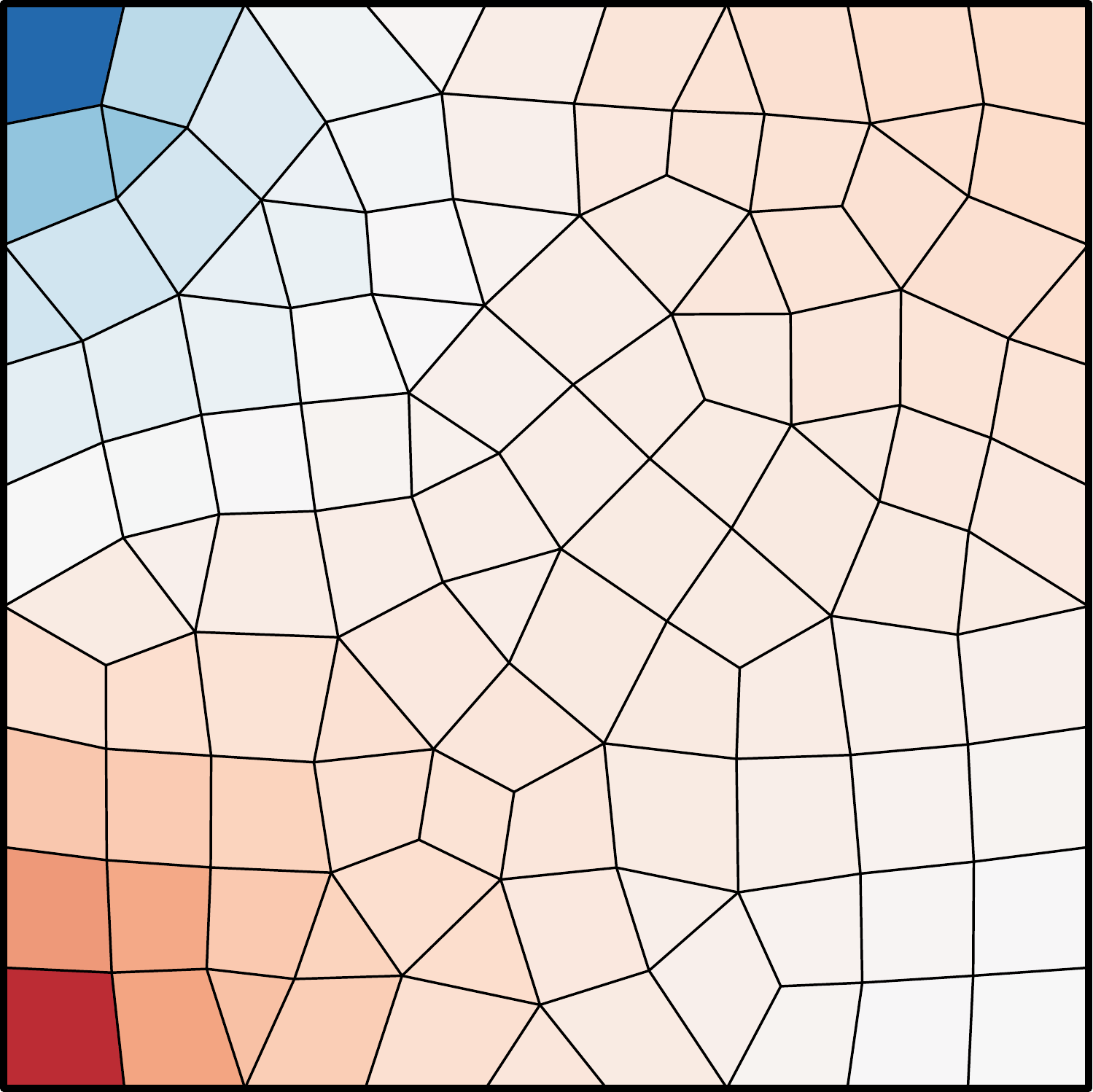}}
      \caption{{\footnotesize{\texttt{Hybrid}}}}
    \end{subfigure}
    \hfill      
    \begin{subfigure}[t]{.18\linewidth}
      \flushright
      \begin{tikzpicture}[scale=1.0]
        \node[anchor=south west,inner sep=0] (image) at (0,0) {\includegraphics[height=.8\linewidth]{./figs/cantilevered_legend}};
        \node [left] at (.5\linewidth,0) {$-2.600$};
        \node [left] at (.5\linewidth,0.2\linewidth) {$-1.300$};
        \node [left] at (.5\linewidth,0.4\linewidth) {$0.000$};
        \node [left] at (.5\linewidth,0.6\linewidth) {$1.300$};
        \node [left] at (.5\linewidth,0.8\linewidth) {$2.600$};
        \node [right] at (0,.9\linewidth) {Pressure [Pa]};
      \end{tikzpicture} 
      %\caption{•}{}
    \end{subfigure}
    \hfill
    \begin{subfigure}[t]{.18\linewidth}
      \hspace{\linewidth}
    \end{subfigure}  
    \hfill\null  
  
  \caption{Cantilevered square block: same as Fig.~\ref{fig:cartilever_cart_skew_hybr_dt_1e_1} for $\Delta t= 1 \times 10^{-1}$ s.}
  \label{fig:cartilever_cart_skew_hybr_dt_1e_1}
\end{figure}

First, we focus on the effectiveness of the proposed stabilization technique in preventing the occurrence of oscillations in the discrete pressure field.
A contour plot of the pressure MFD-VEM solution on the level 0 mesh of each family after a single timestep $\Delta t= 1 \times 10^{-5}$ s is given in Fig.~\ref{fig:cartilever_cart_skew_hybr_dt_1e_5}.
The use of an unstabilized formulation produces checkerboard oscillations for the \texttt{Cartesian}, \texttt{Skewed}, and \texttt{Hybrid} meshes, which are evident in the bottom panels (Fig.~\ref{fig:cartilever_cart_skew_hybr_dt_1e_5}f, g, h).
Such oscillations are successfully removed by the pressure-jump stabilization (Fig.~\ref{fig:cartilever_cart_skew_hybr_dt_1e_5}a, b, c).
As expected, stable results are obtained with meshes \texttt{Polymesher1} and \texttt{Polymesher20} without stabilization (Fig.~\ref{fig:cartilever_cart_skew_hybr_dt_1e_5}d, e) since all the vertices are connected to at most three faces \cite{da2010mimetic}.
Conversely, for larger timestep sizes the MFD-VEM formulation becomes intrinsically stable with no stabilization required  (Fig.~\ref{fig:cartilever_cart_skew_hybr_dt_1e_1}).

\begin{figure}[htbp]
  \small
  \centering
  {
\pgfplotsset{compat=1.11,
    /pgfplots/ybar legend/.style={
    /pgfplots/legend image code/.code={%
       \draw[##1,/tikz/.cd,yshift=-0.25em]
        (0cm,0cm) rectangle (3pt,0.8em);},
   }}
\begin{tikzpicture}
  \pgfplotsset{every axis legend/.append style={
at={(0.5,1.05)},
anchor=south}}
	\begin{axis}[
      width=\linewidth, height=.275\textwidth,
      xtick={1,...,6},
      xticklabels={0,1,2,3,4,5},
      ymin=0,ymax=70,ytick={0,10,...,70},
      xtick pos=bottom,
      xlabel near ticks,    
      ylabel near ticks,
      xlabel=Mesh level,
      ylabel=Iteration number,
      %legend pos=outer north east,
      legend columns = 3,
      transpose legend,
      legend cell align={left},
      legend style={font=\footnotesize},
      ymajorgrids=true,
      ybar=0pt,
      bar width=4.5pt ]
    \addplot [black, fill=mycolor1] table [x=level,y=nIter] {./plotData/cantileverPrecond/DIRECT/scaling_NONE/cart_dt_1_5_stab_0.dat};
    \addlegendentry{\texttt{Cartesian}-\textsc{direct}}; 
    \addplot [black, fill=mycolor1!66] table [x=level,y=nIter] {./plotData/cantileverPrecond/DIRECT/scaling_NONE/cart_dt_1_5_stab_1.dat};  
    \addlegendentry{\texttt{Cartesian}\textsuperscript{($\star$)}-\textsc{direct}};  
    \addplot [black, fill=mycolor1!33] table [x=level,y=nIter] {./plotData/cantileverPrecond/AMG/scaling_NONE/cart_dt_1_5_stab_1.dat}; 
    \addlegendentry{\texttt{Cartesian}\textsuperscript{($\star$)}-\textsc{amg}};

    \addplot [black, fill=mycolor2] table [x=level,y=nIter] {./plotData/cantileverPrecond/DIRECT/scaling_NONE/skew_dt_1_5_stab_0.dat}; 
    \addlegendentry{\texttt{Skewed}-\textsc{direct}}; 
    \addplot [black, fill=mycolor2!66] table [x=level,y=nIter] {./plotData/cantileverPrecond/DIRECT/scaling_NONE/skew_dt_1_5_stab_1.dat};    
    \addlegendentry{\texttt{Skewed}\textsuperscript{($\star$)}-\textsc{direct}};      
    \addplot [black, fill=mycolor2!33] table [x=level,y=nIter] {./plotData/cantileverPrecond/AMG/scaling_NONE/skew_dt_1_5_stab_1.dat};         
    \addlegendentry{\texttt{Skewed}\textsuperscript{($\star$)}-\textsc{amg}};

    \addplot [black, fill=mycolor3] table [x=level,y=nIter] {./plotData/cantileverPrecond/DIRECT/scaling_NONE/hybrid_dt_1_5_stab_0.dat};     
    \addlegendentry{\texttt{Hybrid}-\textsc{direct}};  
    \addplot [black, fill=mycolor3!66] table [x=level,y=nIter] {./plotData/cantileverPrecond/DIRECT/scaling_NONE/hybrid_dt_1_5_stab_1.dat};     
    \addlegendentry{\texttt{Hybrid}\textsuperscript{($\star$)}-\textsc{direct}};   
    \addplot [black, fill=mycolor3!33] table [x=level,y=nIter] {./plotData/cantileverPrecond/AMG/scaling_NONE/hybrid_dt_1_5_stab_1.dat};  
    \addlegendentry{\texttt{Hybrid}\textsuperscript{($\star$)}-\textsc{amg}};

    \addplot [black, fill=mycolor4] table [x=level,y=nIter] {./plotData/cantileverPrecond/DIRECT/scaling_NONE/poly01_dt_1_5_stab_0.dat};
    \addlegendentry{\texttt{Polymesher1}-\textsc{direct}};
    \addlegendimage{empty legend};
    \addlegendentry{};
    \addplot [black, fill=mycolor4!33] table [x=level,y=nIter] {./plotData/cantileverPrecond/AMG/scaling_NONE/poly01_dt_1_5_stab_0.dat};
    \addlegendentry{\texttt{Polymesher1}-\textsc{amg}}; 
    
    \addplot [black, fill=mycolor5] table [x=level,y=nIter] {./plotData/cantileverPrecond/DIRECT/scaling_NONE/poly20_dt_1_5_stab_0.dat};
    \addlegendentry{\texttt{Polymesher20}-\textsc{direct}};
    \addlegendimage{empty legend};
    \addlegendentry{};
    \addplot [black, fill=mycolor5!33] table [x=level,y=nIter] {./plotData/cantileverPrecond/AMG/scaling_NONE/poly20_dt_1_5_stab_0.dat};
    \addlegendentry{\texttt{Polymesher20}-\textsc{amg}};       

  \end{axis}
\end{tikzpicture}
}
  \caption{Cantilevered square block: right-preconditioned GMRES iteration number to solve system \eqref{eq:3x3_system} at the first timestep with $\Delta t = 1 \times 10^{-5}$ s using either nested direct solvers (\textsc{direct}) or algebraic multigrid (\textsc{amg}) to apply both $\widetilde{\Mat{A}}\sub{uu}^{-1}$ and $\widetilde{\Mat{C}}\sub{\pi\pi}^{-1}$. In the legend, the star symbol ($\star$) superscript indicates that the local pressure-jump stabilization was introduced.}
  \label{fig:cantilever_small_dt}
\end{figure}
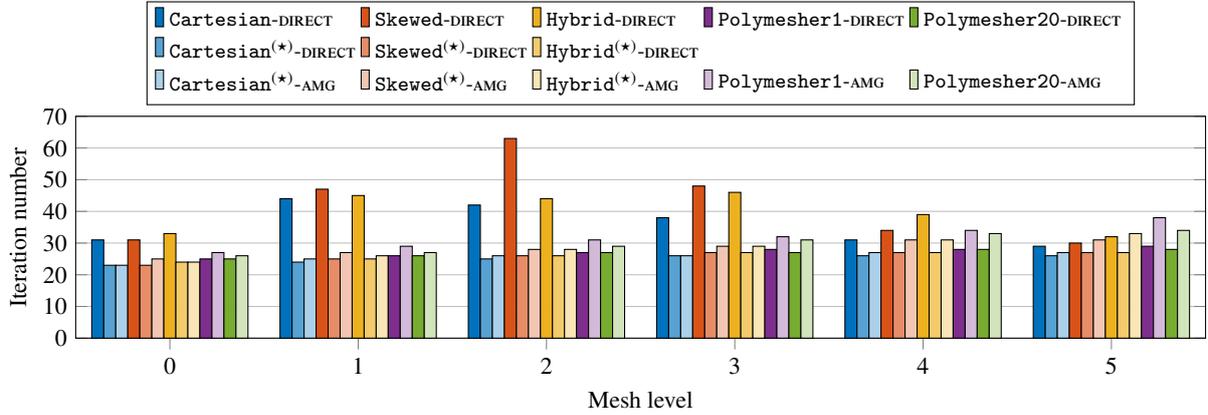

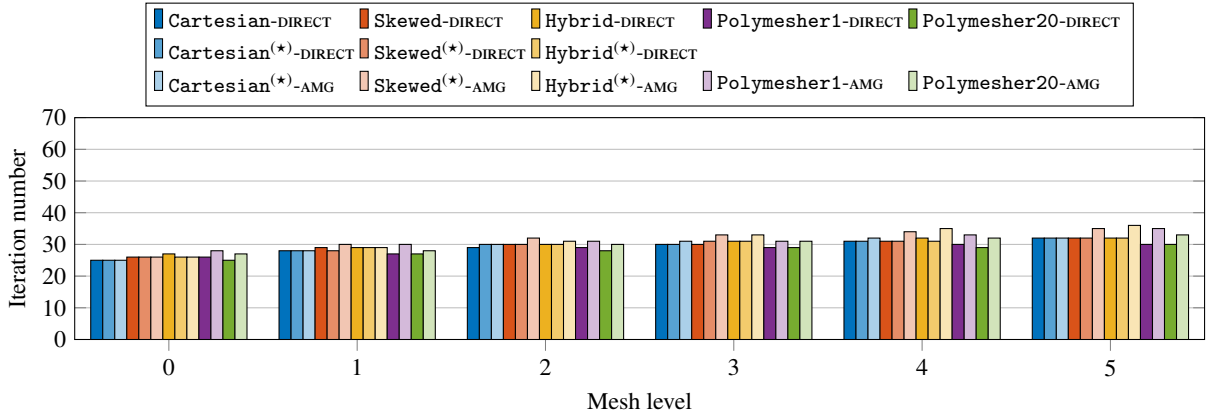
\begin{figure}[htbp]
  \small
  \centering
  {
\pgfplotsset{compat=1.11,
    /pgfplots/ybar legend/.style={
    /pgfplots/legend image code/.code={%
       \draw[##1,/tikz/.cd,yshift=-0.25em]
        (0cm,0cm) rectangle (3pt,0.8em);},
   }}
\begin{tikzpicture}
  \pgfplotsset{every axis legend/.append style={
at={(0.5,1.05)},
anchor=south}}
	\begin{axis}[
      width=\linewidth, height=.275\textwidth,
      xtick={1,...,6},
      xticklabels={0,1,2,3,4,5},
      ymin=0,ymax=70,ytick={0,10,...,70},
      xtick pos=bottom,
      xlabel near ticks,    
      ylabel near ticks,
      xlabel=Mesh level,
      ylabel=Iteration number,
      %legend pos=outer north east,
      legend columns = 3,
      transpose legend,
      legend cell align={left},
      legend style={font=\footnotesize},
      ymajorgrids=true,
      ybar=0pt,
      bar width=4.5pt ]
    \addplot [black, fill=mycolor1] table [x=level,y=nIter] {./plotData/cantileverPrecond/DIRECT/scaling_NONE/cart_dt_1_1_stab_0.dat};
    \addlegendentry{\texttt{Cartesian}-\textsc{direct}}; 
    \addplot [black, fill=mycolor1!66] table [x=level,y=nIter] {./plotData/cantileverPrecond/DIRECT/scaling_NONE/cart_dt_1_1_stab_1.dat};  
    \addlegendentry{\texttt{Cartesian}\textsuperscript{($\star$)}-\textsc{direct}};  
    \addplot [black, fill=mycolor1!33] table [x=level,y=nIter] {./plotData/cantileverPrecond/AMG/scaling_NONE/cart_dt_1_1_stab_1.dat}; 
    \addlegendentry{\texttt{Cartesian}\textsuperscript{($\star$)}-\textsc{amg}};

    \addplot [black, fill=mycolor2] table [x=level,y=nIter] {./plotData/cantileverPrecond/DIRECT/scaling_NONE/skew_dt_1_1_stab_0.dat}; 
    \addlegendentry{\texttt{Skewed}-\textsc{direct}}; 
    \addplot [black, fill=mycolor2!66] table [x=level,y=nIter] {./plotData/cantileverPrecond/DIRECT/scaling_NONE/skew_dt_1_1_stab_1.dat};    
    \addlegendentry{\texttt{Skewed}\textsuperscript{($\star$)}-\textsc{direct}};      
    \addplot [black, fill=mycolor2!33] table [x=level,y=nIter] {./plotData/cantileverPrecond/AMG/scaling_NONE/skew_dt_1_1_stab_1.dat};         
    \addlegendentry{\texttt{Skewed}\textsuperscript{($\star$)}-\textsc{amg}};

    \addplot [black, fill=mycolor3] table [x=level,y=nIter] {./plotData/cantileverPrecond/DIRECT/scaling_NONE/hybrid_dt_1_1_stab_0.dat};     
    \addlegendentry{\texttt{Hybrid}-\textsc{direct}};  
    \addplot [black, fill=mycolor3!66] table [x=level,y=nIter] {./plotData/cantileverPrecond/DIRECT/scaling_NONE/hybrid_dt_1_1_stab_1.dat};     
    \addlegendentry{\texttt{Hybrid}\textsuperscript{($\star$)}-\textsc{direct}};   
    \addplot [black, fill=mycolor3!33] table [x=level,y=nIter] {./plotData/cantileverPrecond/AMG/scaling_NONE/hybrid_dt_1_1_stab_1.dat};  
    \addlegendentry{\texttt{Hybrid}\textsuperscript{($\star$)}-\textsc{amg}};

    \addplot [black, fill=mycolor4] table [x=level,y=nIter] {./plotData/cantileverPrecond/DIRECT/scaling_NONE/poly01_dt_1_1_stab_0.dat};
    \addlegendentry{\texttt{Polymesher1}-\textsc{direct}};
    \addlegendimage{empty legend};
    \addlegendentry{};
    \addplot [black, fill=mycolor4!33] table [x=level,y=nIter] {./plotData/cantileverPrecond/AMG/scaling_NONE/poly01_dt_1_1_stab_0.dat};
    \addlegendentry{\texttt{Polymesher1}-\textsc{amg}}; 
    
    \addplot [black, fill=mycolor5] table [x=level,y=nIter] {./plotData/cantileverPrecond/DIRECT/scaling_NONE/poly20_dt_1_1_stab_0.dat};
    \addlegendentry{\texttt{Polymesher20}-\textsc{direct}};
    \addlegendimage{empty legend};
    \addlegendentry{};
    \addplot [black, fill=mycolor5!33] table [x=level,y=nIter] {./plotData/cantileverPrecond/AMG/scaling_NONE/poly20_dt_1_1_stab_0.dat};
    \addlegendentry{\texttt{Polymesher20}-\textsc{amg}};       

  \end{axis}
\end{tikzpicture}
}
  \caption{Cantilevered square block: same as Fig.~\ref{fig:cantilever_small_dt} for $\Delta t = 1 \times 10^{-1}$ s.}
  \label{fig:cantilever_large_dt}
\end{figure}

Finally, the algorithmic weak scalability of the iterative linear solver strategy described in Section \ref{sec:solution_strategy} is numerically investigated.
We solve \eqref{eq:3x3_system} only once for two different timestep sizes.
Figs.~\ref{fig:cantilever_small_dt}-\ref{fig:cantilever_large_dt} display the number of GMRES iterations needed in the sequence of refined problems detailed in Table \ref{tab:grids} 
When nested direct solvers are used to approximate the action of $\widetilde{\Mat{A}}\sub{uu}^{-1}$ and $\widetilde{\Mat{C}}\sub{\pi\pi}^{-1}$ (Fig.~\ref{fig:cantilever_small_dt}), two behaviors may be observed in case of a small timestep relative to the characteristic consolidation time, $\Delta t = 1 \times 10^{-5}$ s.
For a stable formulation---stabilized (\texttt{Cartesian}\textsuperscript{($\star$)}, \texttt{Skewed}\textsuperscript{($\star$)}, \texttt{Hybrid}\textsuperscript{($\star$)}) or intrinsically stable ((\texttt{Polymesher1}, \texttt{Polymesher20})---the preconditioner $\blkMat{P}^{-1}$ performs very well, with only a modest increase in the iteration count between levels 0 and 5.
Such almost-optimal behavior of $\blkMat{P}^{-1}$ with respect to mesh size is lost in case of an unstable formulation.
Indeed, GMRES shows erratic iteration counts that reflect the sensitivity of Krylov-based solvers to the presence of near-singular pressure modes.
For larger timestep size, i.e. $\Delta t = 1 \times 10^{-1}$ s, the system is far from the incompressibility limit (Fig.~\ref{fig:cantilever_large_dt}).
Hence, the linear solver always exhibits robust convergence and the stabilization effects become irrelevant.
Replacing nested direct solvers by AMG preconditioning for $\widetilde{\Mat{A}}\sub{uu}^{-1}$ and $\widetilde{\Mat{C}}\sub{\pi\pi}^{-1}$ preserves the desired behavior at a substantially smaller cost.

\section{Closure}
\label{sec:concluding_remarks}

In this work, we have presented a numerical scheme coupling hybrid MFD and VEM to discretize Biot's equations of poroelasticity on arbitrary polygonal meshes.
A key feature of the discretization is that it remains convergent in the presence of highly distorted cells with arbitrary shapes, as demonstrated numerically in Section~\ref{sec:numerical_examples}.
For incompressible problems approaching undrained conditions, the discretization is stabilized with a local pressure-jump technique applicable to unstructured meshes to prevent the development of spurious pressure modes (checkerboarding).
The proposed simulation framework also includes a fully coupled linear solution strategy, based on a block-triangular preconditioner, with excellent scalability.

Future work will focus on the extension of these results to three-dimensional polyhedral meshes to further demonstrate the potential of the numerical framework and confirm scalability  in large-scale problems.
A key goal is to address both fully unstructured meshes and stratigraphic corner-point grids \cite{lie2019introduction}.
The latter remain the industry-standard approach in reservoir engineering studies, but contain deformed (and sometimes degenerate) hexahedral cells with non-matching faces that cannot be handled by the standard FEM.

\section*{Acknowledgements}
\label{sec::acknow}
AB acknowledges financial support by the Italian MIUR
project ``Dipartimenti di Eccellenza 2018 - 2022 (CUP
E11G18000350001)'', by INdAM - GNCS, and by Politecnico di Torino
through project ``Starting grant RTD''.
Funding for FPH, NC, JAW, RRS, and HAT was provided by TOTAL S.A. through the FC-MAELSTROM project.
Portions of this work were performed under the auspices of the U.S. Department of
Energy by Lawrence Livermore National Laboratory under Contract DE-AC52-07-NA27344.
FPH completed this work during a visiting scientist appointment at LLNL.
The authors thank Prof. Hamdi Tchelepi for encouragements and fruitful discussions.
\appendix

\section{Unstructured macro-element construction for local pressure-jump stabilization}
\label{sec:macro_element_construction}

Here, we review the partitioning of the mesh involved in the local pressure-jump stabilization.
This is done using the mesh connectivity in two steps described in Algorithm~\ref{alg:macro_element_construction}.
We obtain non-overlapping macro-elements containing at least three cells.

To initialize the algorithm, we mark all the vertices as \textit{not visited yet}.
In the first step, we loop over the internal vertices of the mesh.
If a vertex is marked as \textit{visited}, it is skipped and we proceed to the next vertex in the list. 
If not, we form a macro-element made of the cells adjacent to this vertex, and we mark this vertex as \textit{visited}, as well as all the vertices of the cells forming the macro-element.
At the end of this first step, some cells are, in general, still unassigned to a macro-element.

In the second step, we collect these unassigned cells in a list, and for each cell in the list, we proceed as follows.
If an unassigned cell is adjacent (through faces) to at least one macro-element, this cell is assigned to its neighboring macro-element with the smallest number of cells, and is then removed from the list.
If an unassigned cell is surrounded by unassigned cells, it is placed at the back of the list, and will be processed later.
This procedure is applied until the list is empty, at which point all cells are assigned to a macro-element.

\begin{algorithm}[t]
\caption{Unstructured macro-element construction algorithm}
\label{alg:macro_element_construction}
\begin{algorithmic}[1]
        \State Mark all vertices as \textit{not visited yet}.
        \For{$v$ in $\mathcal{V}_{\textit{int}}$}
           \If{$v$ is marked \textit{not visited yet}}
              \State Form a macro-element with the cells adjacent to $v$.
              \State Mark all the vertices of the cells forming the new macro-element as \textit{visited}.
           \EndIf      
        \EndFor
        \State Collect cells not yet assigned to a macro-element in list $\mathcal{C}$.
        \While{$\mathcal{C}$ is not empty}
           \State Pick cell $K$ at the top of list $\mathcal{C}$.
           \If{$K$ is adjacent (through faces) to at least one macro-element}
             \State Assign $K$ to the neighboring macro-element with the smallest number of cells.
             \State Remove $K$ from list $\mathcal{C}$.
           \Else
             \State Put $K$ at the back of list $\mathcal{C}$.  
           \EndIf
        \EndWhile
 \end{algorithmic}
\end{algorithm}

\section*{References}
\biboptions{sort&compress}
\bibliographystyle{elsarticle-num}
%\bibliography{biblio.bib}
\bibliography{authorList,journalAbbreviations,proceedingCollectionNames,publisherNames,references}

\end{document}